\documentstyle[amscd,amssymb,verbatim]{amsart}
\newcommand{\tot}{\operatorname{tot}}
\newcommand{\secc}{\operatorname{sec}}
\newcommand{\lan}{\langle}
\newcommand{\ran}{\rangle}
\newcommand{\eps}{\varepsilon}

\newcommand{\LL}{{\cal L}}
\newcommand{\RR}{{\cal R}}
\newcommand{\Perv}{\operatorname{Perv}}
\newcommand{\Id}{\operatorname{Id}}
\newcommand{\Fr}{\operatorname{Fr}}
\newcommand{\pr}{\operatorname{pr}}
\newcommand{\gr}{\operatorname{gr}}
\newcommand{\tr}{\operatorname{tr}}

\newcommand{\Mod}{\operatorname{Mod}}

\newcommand{\SL}{\operatorname{SL}}

\newcommand{\Tr}{\operatorname{Tr}}

\newcommand{\HH}{{\cal H}}
\newcommand{\KK}{{\cal K}}
\newcommand{\BB}{{\cal B}}

\newcommand{\CC}{{\cal C}}

\newcommand{\ga}{\gamma}
\newcommand{\si}{\sigma}

\newcommand{\all}{\forall}

\newcommand{\im}{\operatorname{im}}

\newcommand{\Hom}{\operatorname{Hom}}

\newcommand{\Supp}{\operatorname{Supp}}

\newcommand{\ed}{\qed\vspace{3mm}}

\newcommand{\Spec}{\operatorname{Spec}}

\newtheorem{thm}{Theorem}[subsection]
\newtheorem{prop}[thm]{Proposition}
\newtheorem{lem}[thm]{Lemma}
\newtheorem{cor}[thm]{Corollary}
\newenvironment{ex}{\vspace{3mm} \noindent {\bf
Example.}}{\vspace{3mm}}
\theoremstyle{definition}
\newtheorem{defi}[thm]{Definition}
\newenvironment{rem}{\vspace{3mm}
\noindent {\bf Remark.}}{\vspace{3mm}}
\newenvironment{rems}{\vspace{3mm}
\noindent {\bf Remarks.}}{\vspace{3mm}}

\numberwithin{equation}{subsection}

\newcommand{\Pf}{\noindent {\it Proof}}

\renewcommand{\a}{\alpha}
\renewcommand{\b}{\beta}

\newcommand{\wt}{\widetilde}
\renewcommand{\mod}{\operatorname{mod}}

\renewcommand{\ss}{{\bold s}}
\newcommand{\A}{{\Bbb A}}
\renewcommand{\AA}{{\cal A}}

\newcommand{\FF}{{\cal F}}

\newcommand{\coker}{\operatorname{coker}}
\newcommand{\F}{{\Bbb F}}
\newcommand{\Q}{{\Bbb Q}}

\newcommand{\C}{{\Bbb C}}
\newcommand{\Z}{{\Bbb Z}}

\newcommand{\la}{\lambda}

\newcommand{\Ob}{\operatorname{Ob}}

\newcommand{\ov}{\overline}

\newcommand{\ra}{\rightarrow}
\newcommand{\hra}{\hookrightarrow}
\newcommand{\s}{\sigma}
\newcommand{\id}{\operatorname{id}}
\newcommand{\G}{{\Bbb G}}

\newcommand{\ot}{\otimes}

\renewcommand{\P}{{\Bbb P}}
\newcommand{\R}{{\Bbb R}}

\newcommand{\Ext}{\operatorname{Ext}}

\newcommand{\de}{\delta}

\newcommand{\D}{{\cal D}}
\newcommand{\De}{\Delta}

\newcommand{\sub}{\subset}
        
\newcommand{\lrar}[1]{\begin{picture}(50,10)(-25,-5)
\put(-25,0){\vector(1,0){50}}
\put(0,5){\makebox(0,0)[b]{\mbox{$#1$}}}
\end{picture}}

\newcommand{\ldar}[1]{\begin{picture}(10,50)(-5,-25)
\put(0,25){\vector(0,-1){50}}
\put(5,0){\mbox{$#1$}}
\end{picture}}

\newcommand{\ldrar}[1]{\begin{picture}(50,50)(-25,-25)
\put(-25,25){\vector(1,-1){50}}
\put(5,0){\mbox{$#1$}}
\end{picture}}

\newcommand{\ldlar}[1]{\begin{picture}(50,50)(-25,-25)
\put(25,25){\vector(-1,-1){50}}
\put(5,-5){\mbox{$#1$}}
\end{picture}}

\title{Gluing of Perverse Sheaves on the Basic Affine Space}

\begin{document}
\maketitle
\centerline{\sc A.~Polishchuk}
\bigskip
\centerline{\sc with an appendix by R.~Bezrukavnikov and A.~Polishchuk}

%\bigskip

\vspace{6mm}

The goal of this paper is to present the foundations of the
gluing techniques introduced by D.~Kazhdan  and G.~Laumon  in
\cite{KL}. They used this techniques to construct an abelian category
$\AA$ attached to a split semi-simple group $G$ over a finite
field $\F_q$ equipped with actions of $G$, $W$ (the Weil group)
and the Frobenius automorphism. More precisely, the category $\AA$
is glued from $|W|$ copies of the category of perverse sheaves on
the basic affine space $G/U$
(where $U$ is the unipotent radical of the Borel subgroup of $G$).
Conjecturally the category $\AA$
contains information about all ``forms" (in the algebro-geometric sense) of
principal series representations of $G(\F_q)$.
In this paper we explain the Kazhdan---Laumon 
construction in details and study 
gluing of categories in a broader context.

The idea of gluing of categories is based on the geometric model of gluing
a variety from open pieces. One just has to translate geometric data
into the categorical framework by passing to categories of sheaves.
Loosely speaking, in order to glue several categories one has to provide
functors between them which generate a ring-like structure in the
monoidal category of functors. When we want to glue abelian categories
we have to impose in addition some exactness conditions on these
functors. Then the glued category will also be abelian. The principal
difficulty arises when one tries to study homological properties of
the glued category. The intuition coming from the geometric model
here is misleading: it turns out that
even if the original categories have finite cohomological dimension
the glued category can have an infinite dimension.
The original motivation for the present work was a conjecture
made in \cite{KL} that the category $\AA$ has
finite cohomological dimension. This conjecture turned out to be wrong:
a counterexample was found by R.~Bezrukavnikov (see Appendix).
However, for applications to representation theory
a weaker version of this conjecture would suffice. Namely, it would be
enough to prove that objects of finite cohomological dimension
generate the (suitably localized) Grothendieck group of
the category $\AA$ and of its twisted versions. In this paper
we check that this is true (for the category $\AA$ itself) in the cases
when $G$ is of type $A_n$, $n\le 3$, or $B_2$. Also, we have
recorded an important theorem of L.~Positselski asserting that the
individual Ext-groups in the glued category are finite-dimensional (this
was conjectured in \cite{KL}).

The functors used by Kazhdan and Laumon for their gluing 
are certain generalized partial Fourier transforms $F_w$ on
the category of sheaves on $G/U$ 
parametrized by elements $w$ of the Weil  group $W$.
An  important observation made in \cite{KL} is that the functors
$F_w$ produce an action of the generalized positive braid 
monoid $B^+$ associated with $W$ on the category of sheaves on $G/U$.
In this paper we axiomatize this situation in the data which we call
{\it $W$-gluing} (in particular, the categories to be glued
are numbered by elements of $W$).
By combining this abstract formalism with some
geometric considerations we give a proof of Theorem 2.6.1 of \cite{KL}
(the proof given by the authors is insufficient).
We also define ``restriction" and ``induction" functors
between glued categories associated with parabolic subgroups in $W$.
This allows us to deduce
an interesting property of simple objects in the category $\AA$.
Namely, every such object has a {\it support} which is
a subset of $W$. We prove that the support of a simple
object in $\AA$ is either the entire set $W$ or a convex
subset  of  $W$  (in  the  sense of the metric given by reduced
decompositions).
Another topic discussed in this paper is an analogue of the
Kazhdan-Laumon gluing when $W$ is replaced by its ``half":
a subset of elements $w\in W$ with $\ell(sw)=\ell(w)+1$ where
$s$ is a fixed simple reflection (the set of $w$ with
$\ell(sw)=\ell(w)-1$ is ``another half"). 
Applying parabolic induction functors we show that
the homological properties of an object in the corresponding
glued category are completely defined by its restrictions
to the glued categories for proper parabolic subgroups of $W$

The basic idea of Kazhdan and Laumon was that the category $\AA$ ``looks"
as if it were the category of perverse sheaves on a scheme of finite type.
To extend this analogy we introduce a mixed version
of the category $\AA$ and prove a bound on the weights of
Ext-groups between pure objects which is known to hold in
the case of perverse sheaves on a scheme -- this gives one
more reason to believe that the category $\AA$ is ``similar"
to such a category.

As was shown in \cite{BP}, for regular characters of a finite torus 
Kazhdan---Laumon gluing gives the
same representation as the one defined by Deligne and Lusztig in
\cite{DL}. It seems that the methods of the present paper 
can be applied to try to analyze the 
Kazhdan---Laumon construction in the case of 
characters that are ``not very singular".
                         
Here is the plan of the paper. In section 1 we
give  basic  definitions  concerning   gluing   of   abelian
categories. In section \ref{groth} we compute the Grothendieck group
of the glued category (under some mild assumptions).
In section 3 we introduce the notion of $W$-gluing where
$W$ is a Coxeter group by axiomatizing the situation of
gluing on $G/U$. Here we prove a general theorem in the spirit
of Deligne's paper \cite{D} on
{\it quasi-actions} of a Coxeter group $W$ on a category
(by this we mean a collection of functors $(F_w, w\in W)$
and of morphisms $F_w F_{w'}\ra F_{ww'}$ that satisfy
the natural associativity condition but need not be
isomorphisms). Section 4 is devoted to the Fourier transform
on a symplectic vector bundle and its
restriction  to  the  complement of the zero section. In section
\ref{glu} we consider gluing on $G/U$. In particular, we
prove  Theorem  2.6.1  of \cite{KL} which asserts that generalized Fourier
transforms define a quasi-action of $W$ on sheaves on $G/U$.
In section 6 we introduce the technical notion of the
cubic Hecke algebra $\HH^c$ and
prove that the action of the generalized braid group $B$
on $K_0(G/U)$ factors through $\HH^c$. In section \ref{complex}
we give explicit construction of adjoint functors between
glued categories associated with parabolic subgroups of $W$.
Here we prove the theorem on supports of simple objects in
the glued category. Also we construct some canonical complexes
of functors on the glued category. In particular, for objects
supported on ``half"  of  $W$,  we  construct canonical
resolutions consisting  of  objects  induced  from  parabolic
subgroups of $W$. We also prove that for a gluing on "half"
of $W$ our gluing procedure coincides with the one considered
by Beilinson and Drinfeld in \cite{Beil}.
In section 8 we study Ext-groups in the
glued category.
Section 9 is devoted to the gluing of mixed
perverse sheaves. In section \ref{induction} we describe
an  embedding of the category of representations of the braid
group corresponding to a parabolic subgroup $W_J\sub W$
into the category of representations of the braid group $B$.
We also consider an analogue of this picture for quasi-actions
of these groups on categories. Finally, in section 11 we
introduce the notion of a {\it good} representation of the generalized braid
group $B$ and prove that some representations of $B$ are good.
This  is  applied   to   constructing   quotients   of
$K_0(\AA^{\Fr})$ on which the hypothetical bilinear form induces
non-degenerate forms. In Appendix, written jointly with R.~Bezrukavnikov,
we show that the category $\AA$ has infinite cohomological dimension
for $G=\SL_3$.
 
\vspace{2mm}

\noindent
{\it Acknowledgments}. I am grateful to D.~Kazhdan  for
introducing me to the subject and for guidance at
various stages of work on this project. My understanding
of gluing has been clarified a lot due to numerous discussions with
A.~Braverman and L.~Positselski. Moreover, the latter has made an
important contribution by proving the theorem on finiteness of
dimensions of Ext-groups in the glued category, as  well  as
pointing out an error in my earlier optimistic proof of finiteness of
cohomological  dimension. I am grateful to G.~Laumon
and  M.~Rapoport  for  communicating  to  me  the  proof  of
commutativity  of the diagram  (\ref{comm_tr}) which plays an
important role in checking the associativity  condition  for
the morphism from the square of the Fourier transform to the
identity functor. This work was partially supported by
NSF grant DMS-9700458.

\bigskip

\noindent {\bf Notation.}
Let $k$ be a field that is either finite or  algebraically
closed. For a scheme $S$ of finite type over $k$ and a prime
$l$ different from the characteristic of $k$, we denote
by $\D^b_c(S,\ov{\Q}_l)$ the triangulated category of
constructible $l$-adic sheaves on $X$ defined by Deligne in
\cite{DeW}. We denote by $\Perv(S)\sub\D^b_c(S,\ov{\Q}_l)$
the abelian subcategory of perverse sheaves
defined in \cite{BBD}. For a morphism $f:S\ra S'$ between 
such schemes we denote by $f^*$, $f_!$ and $f_*$ the corresponding
derived functors between categories $\D^b_c(S,\ov{\Q}_l)$
and $\D^b_c(S',\ov{\Q}_l)$. For a Coxeter group $(W,S)$
we denote by $\ell:W\ra\Z_{\ge0}$ the corresponding length function,
by $B$ the corresponding generalized braid group and by $B^+$
its positive brais submonoid. We denote by $\tau:W\ra B$
the canonical section given by reduced decompositions.
To avoid confusion with the notation for the set of cosets 
we denote the set-theoretic difference by $X-Y$ (instead of
the customary $X\setminus Y$).

\section{Basic definitions}
\subsection{Gluing of abelian categories}
Let  $(\CC_i)$,  $i=1,\ldots,  n$ be a collection of abelian
categories.
        
\begin{defi}
A {\it (left) gluing  data}  for  $(\CC_i)$   is   a
collection of right-exact (covariant) functors
\begin{equation}\label{functors}
\Phi_{i,j}:\CC_j\ra\CC_i
\end{equation}
for all pairs $(i,j)$ such that  $\Phi_{i,i}=\Id_{\CC_i}$, and a
collection of morphisms of functors
\begin{equation}\label{morphisms}
\nu_{i,j,k}:\Phi_{i,j}\circ\Phi_{j,k}\ra\Phi_{i,k}
\end{equation}
for all triples $(i,j,k)$
such  that  $\nu_{i,i,k}=\id$,   $\nu_{i,j,j}=\id$   and   the
following associativity equation holds:
\begin{equation}\label{asso}
\nu_{i,j,l}\circ(\Phi_{i,j}\nu_{j,k,l})=
\nu_{i,k,l}\circ(\nu_{i,j,k}\Phi_{k,l})
\end{equation}
for all quadruples $(i,j,k,l)$.
\end{defi}
        
One  defines {\it right gluing data} similarly by inverting
arrows and  requiring  functors  to  be  left-exact.
Henceforth, the term ``gluing data" by default will refer
to the ``left gluing data".

For  a  gluing data $\Phi=(\Phi_{i,j};\nu_{i,j,k})$ we define
the category $\CC(\Phi)$  as  follows.  The   objects   of
$\CC(\Phi)$ are collections $(A_i;\a_{ij})$ where $A_i$ is an
object of $\CC_i$ ($i=1,\ldots, n$),
$\a_{ij}:\Phi_{i,j}A_j\ra  A_i$  is  a  morphism  in $\CC_i$
(for every pair $(i,j)$) such that the  following diagram
is commutative:
\begin{equation}\label{comp}
\setlength{\unitlength}{0.25mm}
\begin{array}{ccccc}
\Phi_{i,j}\Phi_{j,k}A_k &
\setlength{\unitlength}{0.50mm}
\lrar{\Phi_{i,j}\a_{jk}}
& \Phi_{i,j}A_j\\
\ldar{\nu_{i,j,k}} & & \ldar{\a_{ij}} \\
\Phi_{i,k}A_k &
\setlength{\unitlength}{0.50mm}
\lrar{\a_{ik}} & A_i
\end{array}
\end{equation}
for every triple $(i,j,k)$.
        
A   morphism   $f:(A_i;\a_{ij})\ra   (A'_i;\a'_{ij})$  is  a
collection of morphisms $f_i:A_i\ra A'_i$ such that
$f_i\circ\a_{ij}=\a'_{ij}\circ\Phi_{i,j}(f_j)$ for all $(i,j)$.
        
\begin{lem} The category $\CC(\Phi)$ is abelian.
\end{lem}
        
\Pf . For a morphism $f:(A_i;\a_{ij})\ra(A'_i,\a'_{ij})$  in
$\CC(\Phi)$  there  are  natural objects of
$\CC(\Phi)$ extending the collections $(\ker(f_i))$ and
$(\coker(f_i))$, which constitute the kernel and the cokernel of
$f$. Indeed, the composition of the natural morphism
$$\Phi_{i,j}\ker(f_j)\ra\Phi_{i,j}A_j\ra A_i$$
with $f_i:A_i\ra A'_i$ is zero; hence, it factors through
a morphism $\Phi_{i,j}\ker(f_j)\ra\ker(f_i)$, and we get
the structure of an object of $\CC(\Phi)$ on $(\ker(f_i))$.
Similarly, the composition of the natural morphism
$$\Phi_{i,j}A'_j\ra A'_i\ra\coker(f_i)$$
with $\Phi_{i,j}f_j:\Phi_{i,j}A_j\ra\Phi_{i,j}A'_j$ is zero;
hence, it factors through a morphism
$\coker(\Phi_{i,j}f_j)\ra\coker(f_i)$. However,
$\coker(\Phi_{i,j}f_j)\simeq\Phi_{i,j}\coker(f_j)$ since
the functor $\Phi_{i,j}$ is right-exact. Thus, we get
a structure of an object of $\CC(\Phi)$ on the collection
$(\coker(f_i))$.
It follows that both the cokernel of kernel and kernel of cokernel of $f$
coincide with the natural object extending   the   collection
$(\im(f_i))$, and therefore, $\CC(\Phi)$ is abelian.
\ed
        
\begin{rems} 

\noindent 1. Dualizing the above construction we obtain
the glued category for right gluing data.
        
\noindent 2. The more general ``gluing" procedure is obtained
by considering an abelian category $\CC$ with
a right-exact functor $\Phi:\CC\ra\CC$ together with
morphisms of functors $\Id\ra\Phi$ and $\Phi^2\ra\Phi$
such that $\Phi$ is a monoid
object in the category of functors from $\CC$ to itself.
In our situation $\CC=\oplus_i \CC_i$ and $\Phi$ has
components $\Phi_{i,j}$. Note that if $\CC$ is the category
of vector spaces over a field $k$, then any $k$-algebra $A$
induces a gluing data on $\CC$ in the above generalized sense.
The functor $\Phi$ in this case is tensoring with $A$ and
the glued category is just the category of $A$-modules.
\end{rems}
        
\subsection{Adjunctions}
Let $\Phi$ be a left gluing data for $(\CC_i)$.
The functor $j_k^*:\CC(\Phi)\ra\CC_k:
(A_i;\a_{ij})\mapsto   A_k$  has   the   left      adjoint
$j_{k,!}:\CC_k\ra\CC(\Phi)$     such     that    $j_k^*\circ
j_{k,!}=\Id_{\CC_k}$.
Namely,
$$j_{k,!}(A)=(\Phi_{i,k}(A);\nu_{i,j,k})$$
where
$$\nu_{i,j,k}:\Phi_{i,j}\Phi_{j,k}(A)\ra\Phi_{i,k}(A)$$
is the structural morphism of the gluing data.
        
The following theorem states that
the existence of such adjoint functors essentially characterizes
the glued category. 

\begin{thm}\label{braverman} 
Let $\CC_k$, $k=1,\ldots,n$ be a collection of abelian
categories, and $\CC$ an abelian category equipped with exact
functors $j_k^*:\CC\ra\CC_k$ for $k=1,\ldots, n$. Assume that
for every $k$ there exists the left adjoint functor $j_{k,!}:\CC_k\ra\CC$
such that $j_k^*\circ j_{k,!}=\Id_{\CC_k}$. Assume also that
for an object $A\in\CC$ the condition $j_k^*A=0$ for all $k$ implies
that $A=0$. Then $\CC$ is equivalent to $\CC(\Phi)$ where $\Phi$
is the gluing data with $\Phi_{ij}=j_i^*\circ j_{k,!}$.
\end{thm}
The proof is not difficult and we leave it to the reader.
A more general statement of this kind is Theorem 2.6 of \cite{BBP}.

Assume that the functor $\Phi_{i,j}$ has the right  adjoint
$\Psi_{j,i}:\CC_i\ra\CC_j$   for  every  $(i,j)$.  Then  the
functors $\Psi_{j,i}$ give rise to a right gluing  data  for
$(\CC_i)$. Namely, by adjunction we have a natural morphism of
functors $\Id_{\CC_k}\ra\Psi_{k,j}\circ\Psi_{j,i}\circ\Phi_{i,k}$
induced by $\nu_{i,j,k}$, so we can form the composition
$$\la_{k,j,i}:\Psi_{k,i}\ra\Psi_{k,j}\circ\Psi_{j,i}\circ
\Phi_{i,k}\circ\Psi_{k,i}\ra\Psi_{k,j}\circ\Psi_{j,i}$$
where $\Phi_{i,k}\circ\Psi_{k,i}\ra\Id_{\CC_i}$
is the canonical adjunction morphism. One can see that the
associativity condition analogous to (\ref{asso}) holds for
$\la_{i,j,k}$.
On the other hand, the morphism of functors
$\nu_{j,i,j}:\Phi_{j,i}\circ\Phi_{i,j}\ra\Id_{\CC_j}$
by adjunction gives rise to a morphism
\begin{equation}\label{mu}
\mu_{i,j}:\Phi_{i,j}\ra\Psi_{i,j}
\end{equation}
for every $(i,j)$.
        
By  construction  the  glued  categories   $\CC(\Phi)$   and
$\CC(\Psi)$ are canonically equivalent.
In particular, by duality we have right adjoint functors
$j_{k,*}:\CC_k\ra\CC(\Phi)$ to the restriction functor
$j_k^*:\CC(\Phi)\ra\CC_k$:
$$j_{k,*}(A_k)=(\Psi_{i,k}(A_k);\a'_{ij})$$
where the morphisms
$\a'_{ij}:\Phi_{i,j}\Psi_{j,k}A_l\ra\Psi_{i,k}A_k$ are deduced by
adjunction from the right gluing data $(\Psi_{i,j},\la_{i,j,k})$.
One has $j_k^*\circ j_{k,*}=\Id_{\CC_k}$.

\subsection{Middle extensions and simple objects}
The morphism  (\ref{mu})  gives rise to a morphism of
functors  $\mu_k:j_{k,!}\ra j_{k,*}$,  such that  for every  object
$A\in\CC(\Phi)$ the composition of the adjunction morphisms
$$j_{k,!}j^*_k A\ra A\ra j_{k,*}j^*_k A$$
coincides with $\mu_l(j^*_k A)$. Thus, we can define
the middle extension functor
$j_{k,!*}:\CC_k\ra\CC(\Phi)$
$$j_{k,!*}(A_k)=\im(j_{k,!}(A)\ra j_{k,*}(A)).$$
        
\begin{lem}\label{ressim}
With the above assumption
for every simple object $A\in\CC(\Phi)$ and any
$l$ the restriction $j^*_lA\in\CC_l$ is either simple or zero.
\end{lem}
        
\Pf . Assume that $j^*_lA\neq 0$ and that there exist an exact sequence
$$0\ra B_l\ra j^*_lA\ra C_l\ra 0$$
with non-zero $B_l$ and $C_l$. Then by adjunction we have
non-zero morphisms $f:j_{l,!}(B_l)\ra A$ and $g:A\ra j_{l,*}(C_l)$ such that
$g\circ f=0$. Since $A$ is simple, $f$ should be surjective,
hence, $g=0$ --- a contradiction. Therefore, $j^*_lA$ is simple.
\ed
        
\begin{lem}\label{gormac}
For every simple object $A_l\in\CC_l$ there is a unique
(up to an isomorphism)
simple object $A\in\CC(\Phi)$ such that $j_l^*A\simeq A_l$.
Namely, $A=j_{l,!*}A_l$.
\end{lem}
        
\Pf . The uniqueness is clear: if $A\in\CC(\Phi)$ is simple
and $j^*_lA\neq 0$,
then  the  adjunction  morphisms  $j_{l,!}j^*_lA\ra  A$  and $A\ra
j_{l,*}j^*_lA$ are surjective and injective, respectively, and hence
$$A\simeq\im(\mu_l(j^*_lA):j_{l,!}j^*_lA\ra j_{l,*}j^*_lA)
=j_{l,!*}(j_l^*A).$$
It remains to check that if $A_l\in\CC_l$  is  simple,  then
$j_{l,!*}(A_l)$  is simple.      Let
$B\sub j_{l,!*}(A_l)$ be a simple subobject. Since
we have  an  inclusion  $B\sub  j_{l,*} A_l$  it  follows  from
adjunction that $j^*_lB$ is a non-zero subobject of $A_l$.
Therefore, $j^*_lB=A_l$ and $B=j_{l,!*}(j^*_lB)=j_{l,!*}(A_l)$.
\ed
                                                
\subsection{Gluing  data for triangulated categories}
We refer to \cite{BBD} for definitions concerning $t$-structures
on triangulated categories.

\begin{defi}  Let  $(\D_i)$  be a collection of triangulated
categories with $t$-structures.
A  {\it gluing  data}  for  $(\D_i)$   is   a
collection of exact functors
\begin{equation}\label{functors2}
\D\Phi_{i,j}:\D_j\ra\D_i
\end{equation}
for all pairs $(i,j)$ that are $t$-exact from the right
(with respect to the given $t$-structures),
such that  $\D\Phi_{i,i}=\Id_{\D_i}$, and a
collection of morphisms of functors
\begin{equation}\label{morphisms2}
\D\nu_{i,j,k}:\D\Phi_{i,j}\circ \D\Phi_{j,k}\ra\D\Phi_{i,k}
\end{equation}
for all triples $(i,j,k)$
such  that  $\D\nu_{i,i,k}=\id$,   $\D\nu_{i,j,j}=\id$   and   the
analogue of the associativity equation (\ref{asso}) holds.
\end{defi}
        
Let $\CC_i$ be the heart of the $t$-structure on $\D_i$.
It is easy to see that the functors
$$H^0\D\Phi_{i,j}|_{\CC_i}=\tau_{\ge 0}\D\Phi_{i,j}|_{\CC_j}$$
(where $\tau_{\ge 0}$ is the truncation with  respect to the $t$-structure)
extend to a gluing data for $\CC_i$. Indeed, since
$\D\Phi_{i,j}$ commutes with $\tau_{\ge 0}$ we have natural morphisms
$$\tau_{\ge 0}\circ\D\Phi_{i,j}\circ
\tau_{\ge 0}\circ\D_{\Phi_{j,k}}|_{\CC_k}\simeq
\tau_{\ge 0}\circ\D\Phi_{i,j}\circ\D_{\Phi_{j,k}}|_{\CC_k}
\lrar{\tau_{\ge 0}\D\nu_{i,j,k}}
\tau_{\ge 0}\circ\D\Phi_{i,k}|_{\CC_k}.
$$
        
\section{Grothendieck groups}\label{groth}
                                 
\subsection{Formulation of the theorem}
Let $\D\Phi=(\D\Phi_{i,j},\D\nu_{i,j,k})$
be the gluing  data for the derived categories
$(\D^b(\CC_i))$ of the abelian categories
$(\CC_i)_{1\le i\le n}$, where $\D^b(\CC_i)$ are
equipped with standard  $t$-structures.
Let $(\Phi_{i,j}=H^0\D\Phi_{i,j}|_{\CC_i},\nu_{i,j,k})$ be
the induced gluing data for $(\CC_i)$.
Then for every $(i,j)$
there  is an induced homomorphism of Grothendieck groups
$$\phi_{ij}=K_0(\Phi_{i,j}):K_0(\CC_j)\simeq K_0(\D_j)\ra
K_0(\D_i)\simeq K_0(\CC_i).$$
        
Let us denote by $\CC_{i,j}$ the full subcategory in $\CC_i$
consisting of objects $A$ such that the morphism
$\nu_{i,j,i,A}:\Phi_{i,j}\Phi_{j,i}(A)\ra A$ is zero.
        
\begin{lem}\label{subc}
$\CC_{i,j}$ is an abelian subcategory in $\CC_i$
closed under passing to quotients and subobjects.
\end{lem}
        
\Pf  .  It is clear that $\CC_{i,j}$ is closed under passing
to subobjects. The statement about quotients follows
from the fact that the functor
$\Phi_{i,j}\Phi_{j,i}$ is right-exact.
\ed
        
Let  us  denote by $K_{i,j}\sub K_0(\CC_i)$ the image of the
natural homomorphism $K_0(\CC_{i,j})\ra K_0(\CC_i)$.
        
\begin{thm}\label{K0}
Let $\D\Phi$ be a gluing data for $(\D^b(\CC_i))$, and
$\Phi$ the corresponding gluing data for $(\CC_i)$.
Assume that all categories  $\CC_i$  are  artinian  and
noetherian. Then the image of the natural map
$K_0(\CC(\Phi))\ra\oplus_{i=1}^n K_0(\CC_i):
[(A_i,\a_{ij})]\mapsto ([A_i])$
coincides with the subgroup
$$K(\Phi):=\{ (c_i)\in\oplus_i K_0(\CC_i) \ | \
\phi_{i,j}c_j-c_i\in K_{i,j} \}.$$
\end{thm}
                   
We need several lemmas for the proof.
                  
\subsection{Homological lemmas}
\begin{lem}\label{inters}
Assume  that the category $\CC_i$ is artinian
and  noetherian.  Then  for  every  collection  of   indices
$j_1$,..., $j_k$ the image of the natural map
\begin{equation}\label{inthom}
K_0(\CC_{i,j_1}\cap\ldots\cap\CC_{i,j_k})\ra K_0(\CC_i)
\end{equation}
coincides with $K_{i,j_1}\cap\ldots\cap K_{i,j_k}$.
\end{lem}
        
\Pf . Since $\CC_i$ is artinian and noetherian,
by the Jordan---H\"older theorem $K_0(\CC_i)$
is a free abelian group with the natural basis corresponding
to isomorphism classes of simple objects in $\CC_i$.
Now the categories $\CC_{i,j}$ are closed under passing to sub- and
quotient-objects by Lemma \ref{subc}.
Hence $K_{i,j}$ is spanned by the classes of simple objects
that belong to $\CC_{i,j}$, while
the  image  of  the  homomorphism
(\ref{inthom}) is spanned by the classes of simple objects
lying   in  $\CC_{i,j_1}\cap\ldots\cap\CC_{i,j_k}$  and  the
assertion follows.
\ed
        
\begin{lem}  Let  $\D$  be  a  triangulated  category with a
$t$-structure such that $\D=\cup_{n} \D^{\le n}$,
$F:\D\ra\D$  is an  exact  functor  which  is
$t$-exact from the right, i.e., $F(\D^{\le0})\sub\D^{\le0}$.
Let $\nu:F\ra\Id$ be a morphism of exact functors. Then for any
object $X\in\D$ such that the morphism $\nu_X:FX\ra X$ is zero the
natural  morphism  $\tau_{\ge 0}F(H^n(X))\ra H^n(X)$ is zero
for any $n$, where $H^n(X)=\tau_{\le n}\tau_{\ge n}(X)[n]$,
$\tau_{\cdot}$ are the truncation functors  associated  with
the $t$-structure.
\end{lem}
        
\Pf  . We may assume that $X\in\D^{\le k}$ for some $k$.
Consider the following morphism of exact triangles:
\begin{equation}
\begin{array}{cccccc}
F(\tau_{\le k-1}X) &\lrar{} & FX &\lrar{} & F(\tau_{\ge k}X)
&\lrar{}\ldots\\
\ldar{\nu} & & \ldar{0} & & \ldar{\nu} \\
\tau_{\le k-1}X &\lrar{} & X &\lrar{} & \tau_{\ge k}X
&\lrar{}\ldots
\end{array}
\end{equation}
Since $F$ is $t$-exact from the right we have
$F(\tau_{\le k-1} X)\in\D^{\le k-1}$, hence
$\Hom^{-1}(F(\tau_{\le k-1}X),\tau_{\ge k}X)=0$.
This implies that all vertical arrows in the above diagram
are zero, hence the conclusion holds for $n=k$ and the
assumption holds for $\tau_{\le k-1}X$, so we may proceed by
induction.
\ed
         
\begin{lem}\label{higherder}
For any object $A_j$ of $\CC_j$ we have
$[H^n\D\Phi_{i,j}A_j]\in K_{i,j}$ for $n\le -1$.
In paricular,
$[\Phi_{i,j}A_j]-\phi_{i,j}[A_j]\in K_{i,j}$.
\end{lem}
        
\Pf . Let us denote $A_i=\Phi_{i,j}A_j\in\CC_i$,
so that we have the morphism $\b:\D\Phi_{j,i}A_i\ra A_j$
induced by $\D\nu_{j,i,j}$.
Consider the following morphism of exact triangles:
\begin{equation}
\setlength{\unitlength}{0.15mm}
\begin{array}{cccccc}
\D\Phi_{i,j}\D\Phi_{j,i}(\tau_{\le       -1}\D\Phi_{i,j}A_j)
&\lrar{} &
\D\Phi_{i,j}\D\Phi_{j,i}(\D\Phi_{i,j}A_j) &\lrar{\pi} &
\D\Phi_{i,j}\D\Phi_{j,i}(A_i) &
\setlength{\unitlength}{0.12mm}
\lrar{}\ldots\\
\ldar{\ga'} & & \ldar{\ga} & & \ldar{} \\
\tau_{\le -1}\D\Phi_{i,j}A_j &\lrar{} & \D\Phi_{i,j}A_j &\lrar{} &
A_i &
\setlength{\unitlength}{0.12mm}
\lrar{}\ldots
\end{array}
\end{equation}
where the vertical arrows are induced by $\D\nu_{i,j,i}$.
The associativity equation for $\D\nu$ implies that
$\ga=\de\circ\pi$ where $\de$ is the morphism
$$\de=\D\Phi_{i,j}\b:\D\Phi_{i,j}\D\Phi_{j,i}A_i\ra
\D\Phi_{i,j}A_j.$$
It  follows  that the composition
$$\D\Phi_{i,j}\D\Phi_{j,i}\tau_{\le -1}\D\Phi_{i,j}A_j
\stackrel{\ga'}{\ra}\tau_{\le -1}\D\Phi_{i,j}A_j\ra
\D\Phi_{i,j}A_j$$
is     zero.     Since    $\D\Phi_{i,j}\D\Phi_{j,i}\tau_{\le
-1}\D\Phi_{i,j}A_j\in\D^{\le -1}$ we get that $\ga'=0$.
By the previous Lemma this implies that
$[H^n\D\Phi_{i,j}A_j]\in K_{i,j}$ for $n\le -1$ as required.
\ed
                        
\begin{lem}\label{devis}
Let   $\CC^n(\Phi)$  be  the  full
subcategory    of    $\CC(\Phi)$   consisting   of   objects
$(A_i;\a_{ij})$ with $A_n=0$. Then $\CC^n(\Phi)$
is equivalent to $\CC(\Phi')$ for some new gluing data $\Phi'$ for $n-1$
categories $\CC_{1,n}$,..., $\CC_{n-1,n}$.
\end{lem}
        
\Pf . If $(A_i;\a_{ij})$ is an object of $\CC^n(\Phi)$,
then the condition (\ref{comp}) implies that
the composition
$$\Phi_{i,n}\Phi_{n,j}A_j\stackrel{\nu_{i,n,j}}{\ra}
\Phi_{i,j}A_j\stackrel{\a_{ij}}{\ra} A_i$$
is zero for every pair $(i,j)$, in particular,
$A_i\in\CC_{i,n}$  for  every  $i$. For every pair $(i,j)$
such that $i,j\le n-1$ and $i\neq j$,
let us denote by $\Phi'_{i,j}$ the cokernel of the morphism
of functors $\nu_{i,n,j}:\Phi_{i,n}\Phi_{n,j}\ra\Phi_{i,j}$.
Then we have a canonical morphism of functors $\Phi_{i,j}\ra
\Phi'_{i,j}$ and
the above condition means that $\a_{ij}$ factors through a
morphism $\a'_{ij}:\Phi'_{i,j}A_j\ra A_i$. It is easy  to
see that the functor $\Phi'_{i,j}$ is  right-exact  (as  the
cokernel of right-exact functors). Let us check that
$\Phi'_{i,j}$ sends $\CC_j$ to $\CC_{i,n}$.
Indeed, for any $A_j\in\CC_j$ the morphism
\begin{equation}\label{phi'}
\nu_{i,n,i}:\Phi_{i,n}\Phi_{n,i}\Phi_{i,j}A_j\ra\Phi_{i,j}A_j
\end{equation}
factors through a morphism
$\Phi_{i,n}\Phi_{n,i}\Phi_{i,j}A_j\ra\Phi_{i,n}\Phi_{n,j}A_j$.
Hence, the composition of (\ref{phi'}) with the projection
$\Phi_{i,j}A_j\ra\Phi'_{i,j}A_j$ is zero which implies that
$\Phi'_{i,j}A_j\in\CC_{i,n}$.
Furthermore, there are unique
morphisms of functors $\Phi'_{i,j}\Phi'_{j,k}\ra\Phi'_{i,k}$
compatible with $\nu_{i,j,k}$, so that we get a gluing data
$\Phi'$ for $\CC_{1,n}$,...,$\CC_{n-1,n}$ such that
$(A_i,i=1,\ldots,n-1;\a'_{ij})$ is an object of $\CC(\Phi')$.
Clearly, this gives an equivalence of the category $\CC^n(\Phi)$
with $\CC(\Phi')$.
\ed
        
\subsection{Proof of theorem \ref{K0}}
Let  us  check  first  that the image of $K_0(\CC(\Phi))$ is
contained in $K(\Phi)$. It suffices to check that for every $i$
and every $A\in\CC(\Phi)$ one has $\phi_{i,n}[A_n]-[A_i]\in K_{i,n}$.  
We note that for every $A\in\CC(\Phi)$
the kernel and the cokernel
of  the natural morphism $j_{n,!}j^*_nA\ra A$ belong
to $\CC^n(\Phi)$. Hence, $K_0(\CC(\Phi))$ is generated by
the image of the map $K_0(\CC^n(\Phi))\ra K_0(\CC(\Phi))$ and
by  the  classes of  $j_{n,!}(A_n)$,  $A_n\in\CC_n$.
Let us check the above condition for these two classes of elements
separately. If
$(A_i;\a_{ij})\in\CC^n(\Phi)$,      then     by     definition
$A_i\in\CC_{i,n}$ for any $i$, while $A_n=0$; hence
$\phi_{i,n}[A_n]-[A_i]=-[A_i]\in K_{i,n}$ for all $i$. On
the other hand, if $(A_i;\a_{ij})=j_{n,!}(A_n)$, then
$A_i=\Phi_{i,n}A_n$ and
$$\phi_{i,n}[A_n]-[A_i]=\phi_{i,n}[A_n]-[\Phi_{i,n}A_n]\in
K_{i,n}$$
by Lemma \ref{higherder}. Thus, the condition
$\phi_{i,n}[A_n]-[A_i]\in  K_{i,n}$  is  satisfied  for  all
objects of $\CC(\Phi)$.

It remains to check that the map $K_0(\CC(\Phi))\ra K(\Phi)$
is surjective.
Let us denote  by  $\CC^l(\Phi)$  the  full  subcategory  of
$\CC(\Phi)$ consisting of all objects $(A_i;\a_{ij})$ with
$A_i=0$ for $i\ge l$. Similarly, let
$K^l(\Phi)$ be the subgroup of elements $(c_i)$
in $K(\Phi)$ with $c_i=0$ for $i\ge l$. Note that the image
of $K_0(\CC^l(\Phi))$ is contained in $K^l(\Phi)$.
Let  $p_j$ be the projection of $\oplus_i K_0(\CC_i)$ on its $j$-th
factor. Then $p_{l-1}(K^l(\Phi))\sub \cap_{j\ge l}K_{l-1,j}$
since for $(c_i)\in K^l(\Phi)$ and $j\ge l$ we have
$c_{l-1}=c_{l-1}-\phi_{l-1,j}c_j\in K_{l-1,j}$. Thus,
we have an exact sequence
$$0\ra K^{l-1}(\Phi)\ra K^l(\Phi)\ra \cap_{j\ge l} K_{l-1,j}.$$
Hence, to
prove the surjectivity of the map $K_0(\CC(\Phi))\ra K(\Phi)$
it is sufficient to check that the map
\begin{equation}\label{proj}
K_0(\CC^l(\Phi))\ra \cap_{j\ge l} K_{l-1,j}:
[(A_i;\a_{ij})]\mapsto [A_{l-1}]
\end{equation}
is surjective for each $l$.
Note that  for  every
gluing data $\Phi$ (not necessarily extending to the derived
categories) the natural  homomorphism   $K_0(\CC(\Phi))\ra
K_0(\CC_n)$ is surjective ($[j_{n,!}(A_n)]$ maps to $[A_n]$).
Now the iterated  application  of  Lemma  \ref{devis}  gives  an
equivalence  of  $\CC^l(\Phi)$  with  $\CC(\Phi')$  for some
gluing data $\Phi'$ on $l-1$ categories
$\cap_{j\ge  l}\CC_{1,j}$,...,  $\cap_{j\ge  l}\CC_{l-1,j}$,
such  that  the  map (\ref{proj}) is identified with the
homomorphism
$$K_0(\CC(\Phi'))\ra K_0(\cap_{j\ge l}\CC_{l-1,j})\ra
\cap_{j\ge l} K_{l-1,j}$$
which is surjective by Lemma \ref{inters}.
\ed
                                           
\subsection{Gluing of finite type}
\begin{prop}\label{inj}
Assume that every functor $\Phi_{i,j}$ has
the right adjoint $\Psi_{j,i}$ and  that  the  categories  $\CC_i$ are
artinian and noetherian. Then the natural homomorphism
$K_0(\CC(\Phi))\ra\oplus_i K_0(\CC_i)$ is
injective.
\end{prop}
        
\Pf  .  Since  the  category  $\CC(\Phi)$  is  artinian  and
noetherian,  it  follows  that  $K_0(\CC(\Phi))$ is the free
abelian group with the basis $[A]$, where $A$ runs through
the isomorphism classes of simple objects in $\CC(\Phi)$.
Now   the   assertion   follows   immediately   from  Lemmas
\ref{ressim} and \ref{gormac}.
\ed
                                                 
\begin{defi}  We  say that $\Phi$ is a gluing data of {\it finite
type} if all the categories $\CC_i$ are
artinian and noetherian, and
every functor $\Phi_{ij}$ has the right adjoint
$\Psi_{ji}$.
\end{defi}
        
Combining Proposition \ref{inj}  with  Theorem  \ref{K0}  we
obtain the following result.
        
\begin{thm}\label{main}
Let $\D\Phi$ be a gluing data for $(\D^b(\CC_i))$,
and $\Phi$ the corresponding gluing data for $(\CC_i)$.
Assume that $\Phi$ is of finite type.
Then the natural map
$K_0(\CC(\Phi))\ra\oplus_{i=1}^n K_0(\CC_i):
[(A_i,\a_{ij})]\mapsto ([A_i])$
induces an isomorphism $K_0(\CC(\Phi))\simeq K(\Phi)$.
\end{thm}

\section{Gluing for Coxeter groups}
                                                    
\subsection{$W$-gluing}
Let $(W,S)$ be a finite Coxeter group,
and $\ell:W\ra\Z_{\ge 0}$ be the length function.
Let $(\CC_w; w\in W)$ be a collection of abelian categories.
        
\begin{defi} A $W$-{\it gluing data} is a gluing data
$\Phi_w:\CC_{w'}\ra\CC_{ww'}$, $\nu_{w,w'}:\Phi_w\circ\Phi_{w'}\ra
\Phi_{ww'}$  for  $(\CC_w)$  such  that
$\nu_{w,w'}$ is an isomorphism for every pair
$w,w'\in W$ such that $\ell(ww')=\ell(w)+\ell(w')$.
\end{defi}
       
Note that the condition of finiteness of $W$ is imposed only
because we usually consider the gluing of a finite number
of categories. One can similarly treat the case of infinite number of
categories and define $W$-gluing data in this context.
       
\subsection{Quasi-actions of Coxeter groups}
\label{quasiactcox}
       
\begin{defi}  A  {\it  quasi-action}  of a monoid $M$ on a
category $\CC$ is a collection of functors
$T(f)$, $f\in M$, from $\CC$ to itself and of morphisms of functors
$c_{f,g}:T(f)\circ T(g)\ra T(fg)$, $f,g\in M$ satisfying the
associativity condition.
\end{defi}
        
An {\it action} of a monoid is a quasi-action such  that  all
morphisms $c_{f,g}$ are isomorphisms.
Let $(W,S)$ be a Coxeter system. The above definition of
$W$-gluing data can be reformulated as follows: there is
a quasi-action of $W$ on $\oplus_{w\in W}\CC_w$ given by
functors     $\Phi_w:\CC_{w'}\ra\CC_{ww'}$     such     that
$\Phi_1=\Id$   and   $c_{w,w'}$   are   isomorphisms    when
$\ell(ww')=\ell(w)+\ell(w')$.
   
Let us denote by $B$ the generalized
braid group corresponding to $(W,S)$ and by
$B^+\sub B$ the positive braids submonoid (see \cite{D}).
Let us denote by $b\mapsto\ov{b}$ the natural homomorphism
$B\ra W$ and by $\tau:W\ra B^+$ its canonical section
such that $\tau$ is the identity on $S$ and
$\tau(ww')=\tau(w)\tau(w')$ whenever
$\ell(ww')=\ell(w)+\ell(w')$.
According to the main theorem of \cite{D}, to give an action
of $B^+$ on a category $\CC$, it is sufficient to have functors
$T(w)$ corresponding to elements $\tau(w)\in B^+$
and isomorphisms $c_{w,w'}:T(w)\circ T(w')\wt{\ra} T(ww')$
for pairs $w,w'\in W$ such that $\ell(ww')=\ell(w)+\ell(w')$,
satisfying the associativity condition for
triples $w,w',w''\in W$, such that
$\ell(ww'w'')=\ell(w)+\ell(w')+\ell(w'')$.
It follows that a $W$-gluing   data   induces   an   action  of  $B^+$  on
$\oplus_{w}\CC_w$. Conversely, assume that we are given the
functors  $\CC_{w}$  and  isomorphisms  $c_{w,w'}$  only  for
$(w,w')$ with $\ell(ww')=\ell(w)+\ell(w')$ so that we have an action of
$B^+$ on $\oplus_w\CC_w$. To get a quasi-action of $W$ we have
to      give      in      addition      some       morphisms
$c_{s,s}:\Phi_s\Phi_s\ra\Id$.
The natural question arises as to what compatibility conditions
relating $c_{s,s}$ and the action of $B^+$ one should impose
to obtain a quasi-action of $W$. The answer is given by
the following theorem.
        
\begin{thm}\label{quasiact}
Assume that we have an action of  $B^+$  on  a
category  $\CC$  given by the collection of functors $T(w)$,
$w\in W$ (such that $T(1)=\Id$),
and isomorphisms $c_{w,w'}:T(w)T(w')\wt{\ra}T(ww')$
for $w,w'\in W$ such that $\ell(ww')=\ell(w)+\ell(w')$. Let
$c_{s,s}:T(s)T(s)\ra\Id$, $s\in S$, be a collection of  morphisms
satisfying the following two conditions:
\begin{enumerate}
\item
For every $s\in S$ the following associativity equation holds:
\begin{equation}\label{asso_simple}
T(s) c_{s,s,X}=c_{s,s, T(s)X}:T(s)^3X\ra T(s)X.
\end{equation}
\item
For every $w\in W$ and $s,s'\in S$ such that $sw=ws'$ and
$\ell(sw)=\ell(w)+1$ the following diagram is
commutative:
\begin{equation}\label{sizig}
\begin{array}{ccccc}
T(s)T(w)T(s') &\lrar{} & T(w)T(s')^2\\
\ldar{} & & \ldar{c_{s',s'}}\\
T(s)^2T(w) &\lrar{c_{s,s}} & T(w)
\end{array}
\end{equation}
where  the  unmarked  arrows  are  induced  by $c_{s,w}$ and
$c_{w,s'}$.
\end{enumerate}
Then there is a canonical quasi-action of $W$ on $\CC$.
\end{thm}
        
\Pf . Following \cite{B} we denote by $P_s$ the
set of $w\in W$ such that $\ell(sw)=\ell(w)+1$.
Let us first construct the canonical morphisms
$c_{s,w}:T(s)T(w)\ra  T(sw)$
for every $s\in S$, $w\in W$.
When $w\in P_s$ they
are given by the structure of $B^+$-action. Otherwise
$sw\in P_s$ and we have the morphism
$c_{s,w}:T(s)T(w)\wt{\ra} T(s)T(s)T(sw)\ra T(sw)$
induced  by  $c_{s,sw}^{-1}$  and
$c_{s,s}$. It is easy to check using (\ref{asso_simple})
that the following triangle
is commutative for every $s\in S$, $w\in W$:
\begin{equation}\label{asso0}
\begin{array}{ccc}
T(s)^2T(w) & &\\
\ldar{c_{s,w}} & \ldrar{c_{s,s}} &\\
T(s)T(sw) &\lrar{c_{s,ws}}& T_{w}
\end{array}
\end{equation}
Similarly, one defines morphisms
$c_{w,s}:T(w)T(s)\ra T(ws)$ for every $s\in S$, $w\in W$.
     
We claim that these  morphisms  satisfy  the
following  associativity  condition:  for  any   $w\in   W$,
$s,s'\in S$ the diagram
\begin{equation}\label{asso1}
\begin{array}{ccccc}
T(s)T(w)T(s')&\lrar{c_{w,s'}}& T(s)T(ws')\\
\ldar{c_{s,w}}& &\ldar{c_{s,ws'}}\\
T(sw)T(s')&\lrar{c_{sw,s'}}& T(sws')
\end{array}
\end{equation}
is  commutative. Assume at first that $w\in P_s$.
Consider two cases.
\begin{enumerate}
\item $ws'\in P_s$.
If  $\ell(ws')=\ell(w)+1$  then  $\ell(sws')=\ell(w)+2$ and the required
associativity holds by definition of the $B^+$-action.
Otherwise, $\ell(w)=\ell(ws')+1$,   and hence
$\ell(sws')=\ell(ws')+1=\ell(w)=\ell(sw)-1$, $T(w)\simeq T(ws')T(s')$,
$T(sw)\simeq T(sws')T(s')$, and we are reduced to the
commutativity of the diagram
\begin{equation}\label{asso1_1}
\begin{array}{ccccc}
T(s)T(ws')T(s')^2&\lrar{c_{s',s'}}& T(s)T(ws')\\
\ldar{c_{s,ws'}}& &\ldar{c_{s,ws'}}\\
T(sws')T(s')^2&\lrar{c_{s',s'}}& T(sws')
\end{array}
\end{equation}
which is clear.
\item $w\in P_s$, $ws'\not\in P_s$. In this case
according to \cite{B}, IV, 1.7 we have $sw=ws'$, so the
required associativity follows from (\ref{sizig}).
\end{enumerate}
     
Thus, the diagram (\ref{asso1}) is commutative for $w\in P_s$.
Now assume that $w=sw'$ with $w'\in P_s$.
Consider the following diagram:
\begin{equation}\label{asso1_2}
\setlength{\unitlength}{0.17mm}
\begin{array}{ccccccc}
T(s)T(sw')T(s')&\lrar{c_{s,w'}^{-1}}&T(s)^2T(w')T(s')&
\lrar{c_{w',s'}}& T(s)^2T(w's') &
\lrar{c_{s,w's'}}&T(s)T(sw's')\\
&\ldrar{c_{s,sw'}}&\ldar{c_{s,s}}&
&\ldar{c_{s,s}}&\ldlar{c_{s,sw's'}}&\\
&&T(w')T(s')&\lrar{c_{w',s'}}& T(w's')
\end{array}
\end{equation}
In this diagram the square is commutative,
the left triangle is commutative by
definition of $c{s,sw'}$ and the right triangle
is commutative by (\ref{asso0}). Also the commutativity
of (\ref{asso1}) for $w'$ implies that the composition
of the top arrows coincides with the top arrow in
the diagram (\ref{asso1}) for $w$ and hence, the commutativity
of (\ref{asso1}) for $w$.
     
Now for any $w\in W$ with a reduced decomposition
$w=s_{i_1}s_{i_2}\ldots s_{i_l}$ which we denote by $\ss$
and any $w'\in W$,
we define the morphism $c_{\ss,w'}:T(w)T(w')\ra T(ww')$
inductively as the composition
$$c_{\ss,w'}:T(w)T(w')\simeq
T(s_{i_1})T(s_{i_2}\ldots s_{i_l})T(w')\ra
T(s_{i_1})T(s_{i_2}\ldots s_{i_l}w')\ra T(ww')$$
where   the   latter   arrow   is  $c_{s_{i_1},s_{i_2}\ldots
s_{i_l}w'}$. We are going to prove that $c_{\ss,w'}$ does not
depend on a choice of the reduced decomposition $\ss$ of $w$.
First we note that the following diagram
is  commutative for every reduced decomposition $\ss$ of $w$
and every $s\in S$, $w'\in W$:
\begin{equation}\label{asso2}
\begin{array}{ccccc}
T(w)T(w')T(s)&\lrar{c_{w,s}}& T(w)T(w's)\\
\ldar{c_{\ss,w'}}& &\ldar{c_{\ss,w's}}\\
T(ww')T(s)&\lrar{c_{ww',s}}& T(ww's)
\end{array}
\end{equation}
Indeed, this follows from (\ref{asso1}) by induction in the
length of $w$. Now the required independence of $c_{\ss,w'}$
on a choice of $\ss$ follows from (\ref{asso2}) by induction
in the length of $w'$ (the base of the induction is the case
$w'=1$ when the assertion is obvious). Let us denote
$c_{w,w'}=c_{\ss,w'}$ for any reduced decomposition $\ss$ of
$w$. Using (\ref{asso2}) it is easy to show that one would
obtain  the  same   morphism   starting   with   a   reduced
decomposition of $w'$. It remains to check that the following
diagram is commutative for any $w,w_1,w_2\in W$:
\begin{equation}\label{asso3}
\begin{array}{ccccc}
T(w_1)T(w)T(w_2)&\lrar{c_{w,w_2}}& T(w_1)T(ww_2)\\
\ldar{c_{w_1,w}}& &\ldar{c_{w_1,ww_2}}\\
T(w_1w)T(w_2)&\lrar{c_{w_1w,w_2}}& T(w_1ww_2)
\end{array}
\end{equation}
When $\ell(w_2)=1$ this reduces to (\ref{asso2}).
The general case follows easily by induction in $\ell(w_2)$.
\ed
        
\subsection{Reduction to rank-2 subgroups}
The condition (2) of Theorem \ref{quasiact}
can be reduced to rank-2 subgroups of $W$ using the following result.
For every pair $(s_1,s_2)$ of elements of $S$, such that
the order of $s_1s_2$ is $2m+\eps$ with $\eps\in\{0,1\}$,
let us denote $s=s(s_1,s_2)=s_1$, $s'=s'(s_1,s_2)=s_1$ if $\eps=0$,
and $s'=s_2$ if $\eps=1$.
Let also $w=w(s_1,s_2)=(s_2s_1)^{n-\eps}s_2^{\eps}$.
Then $sw=ws'$ and $\ell(sw)=\ell(w)+1$.
        
\begin{prop}\label{sizig2} We keep the notation and assumptions
of Theorem \ref{quasiact}.
Assume that the diagram (\ref{sizig}) is commutative
for all the triples $s(s_1,s_2),s'(s_1,s_2),w(s_1,s_2)$
associated with pairs $(s_1,s_2)$ of elements of $S$ as above.
Then it is commutative for all $(w,s,s')$ such that
$sw=ws'$ and $\ell(sw)=\ell(w)+1$.
\end{prop}
                                                
%First, an easy lemma on Coxeter groups.
       
\begin{lem}\label{sizig3} Let $(W,S)$ be a Coxeter system.
Assume that $sw=ws'$, where
$w\in W$, $s,s'\in S$, $\ell(sw)=\ell(w)+1$, $w\neq 1$.
Then there exists an element $s_2\in S$, such that
for the element $w(s_1,s_2)$ associated with the pair
$(s_1=s,s_2)$, one has
$\ell(w)=\ell(w(s_1,s_2))+\ell(w(s_1,s_2)^{-1}w)$.
\end{lem}
                                                                  
\Pf. Consider a reduced decomposition $(s_2,\ldots, s_n)$ of
$w$. Then we have two reduced decompositions of $sw=ws'$:
${\bf s}=(s_1=s,s_2,\ldots,s_n)$ and ${\bf s'}=(s_2,\ldots,s_n,s')$. Now
the second sequence is obtained from the first one by a series
of standard moves associated with braid relations for couples of elements of
$S$. If the first move in such a series
does not touch the first member of ${\bf s}$, then it just
changes a reduced decomposition of $w$. Thus, we can choose
an  initial reduced decomposition for $W$ in such a way that
the first move does touch $s_1$. Then we have $w=w(s_1,s_2)w'$
where $\ell(w)=\ell(w(s_1,s_2))+\ell(w')$ as required.
\ed
                                              
\noindent
{\it Proof of Proposition \ref{sizig2}}.
Induction on the length of $w$ and Lemma \ref{sizig3}
show that it is sufficient to prove the following:
if $w=w_1w_2$ with $\ell(w)=\ell(w_1)+\ell(w_2)$ and one
has $s_1w_1=w_1s_2$, $s_2w_2=w_2s_3$ where $\ell(s_1w)=\ell(w)+1$,
then the commutativity of the diagram (\ref{sizig})
for the triples $(w_1,s_1,s_2)$ and $(w_2,s_2,s_3)$
implies its commutativity for $(w,s_1,s_3)$.
This can be easily checked using the fact that we have a
$B^+$-action.
\ed
       
\subsection{Grothendieck groups}
As before we can define a notion of $W$-gluing data $\D\Phi$ for
the derived categories $(\D^b(\CC_w))$, which induce the
$W$-gluing data $\Phi$ for $(\CC_w)$.
Let us denote by $\phi_w:K_0(\CC_w')\ra
K_0(\CC_{ww'})$   the   corresponding    homomorphisms    of
Grothendieck  groups.  If in addition the categories $\CC_w$
are artinian and noetherian, then according to Theorem \ref{K0}
the image  of  the  homomorphism  $K_0(\CC(\Phi))\ra\oplus_w
K_0(\CC_w)$ coincides with the subgroup
$$K(\Phi)=\{ (c_w)\in\oplus_{w\in W} K_0(\CC_w) \ | \
\phi_w c_{w'}-c_{ww'}\in K_{ww',w'},\ w,w'\in W \}.$$
        
The property that $\nu_{w,w'}$ is an isomorphism
when $\ell(ww')=\ell(w)+\ell(w')$ allows us to give an alternative
definition for the category $\CC(\Phi)$. Namely,
the morphisms $\a_{w',w}:\Phi_{w'}A_w\ra A_{w'w}$
for all $w'\in W$
can be recovered uniquely from the morphisms
$\a_{s,w}:\Phi_sA_w\ra A_{sw}$ provided the latter morphisms
respect the relations in $W$ in the obvious sense.
On the level of Grothendieck groups  this is reflected  in  the
following result.
        
\begin{prop}\label{simple}
Let $\D\Phi$ be a $W$-gluing data for
$(\D^b(\CC_w))$. Then
\begin{equation}\label{newK0}
K(\Phi)=\{ (c_w)\in\oplus_{w\in W} K_0(\CC_w) \ | \
\phi_s c_w-c_{sw}\in K_{sw,w}, \ w\in W, s\in S \}.
\end{equation}
\end{prop}
        
\Pf . It is clear that the left-hand side of (\ref{newK0})
is contained in the right-hand side, so we have to check the
inverse inclusion.
If $\ell(sw)=\ell(w)+1$, then $\phi_{sw}=\phi_s\phi_w$. Hence,
$$\phi_{sw} c_{w'}-c_{sww'}=(\phi_s c_{ww'}- c_{sww'})+
\phi_s (\phi_w c_{w'}-c_{ww'}).$$
Thus, it is sufficient to prove that
$$K_{sww',ww'}+\phi_s(K_{ww',w'})\sub K_{sww',w'}$$
provided that $\ell(sw)=\ell(w)+1$.
By    definition   $\CC_{sww',w'}$   consists   of   objects
$X\in\CC_{sww'}$ for which the morphism
\begin{equation}\label{mor}
\Phi_s\Phi_{w}\Phi_{w^{-1}}\Phi_s X\ra X
\end{equation}
is zero. Since this morphism factors  through  $\Phi_s\Phi_s
X\ra X$ we have an obvious inclusion $\CC_{sww',ww'}\sub
\CC_{sww',w'}$. By Lemma \ref{higherder} for $Y\in\CC_{ww'}$
we have $[H^n\D\Phi_s Y]\in K_{sww',ww'}$ for $n\le -1$. Thus,
it is sufficient to prove the inclusion
$\Phi_s(\CC_{ww',w'})\sub\CC_{sww',w'}$.       But       for
$X=\Phi_s Y$ the morphism (\ref{mor}) factorizes as follows:
$$\Phi_s\Phi_w\Phi_{w^{-1}}\Phi_s\Phi_sY\ra
\Phi_s\Phi_w\Phi_{w^{-1}}Y\ra\Phi_s Y.$$
If $Y\in\CC_{ww',w'}$, then the latter arrow is zero, hence,
$\Phi_s Y\in\CC_{sww',w'}$.
\ed
       
\section{Symplectic Fourier transform}
                                                
\subsection{Functors and distinguished triangles}\label{foursq}
Let $k$ be a field of characteristic $p>0$ that is
either finite or algebraically closed, and $S$ a
scheme of  finite type over $k$.
Let  $\pi:V\ra S$ be a symplectic vector bundle of rank $2n$ over
$S$, and  $\langle  ,  \rangle:  V\times_S  V\ra\G_a$ the corresponding
symplectic   pairing.   Let us fix a non-trivial additive character
$\psi:\F_p\ra\ov{\Q}_l^*$.
The Fourier---Deligne transform
$\FF=\FF_\psi$ is the  involution  of  $\D^b_c(V,\ov{\Q}_l)$
defined by 
$$\FF(K)=p_{2!}(\LL\ot p_1^*(K))[2n](n).$$
where $p_i$ are the projections of the product $V\times_S V$ on its
factors, and $\LL=\LL_{\psi}(\langle , \rangle)$ is a smooth rank-1
$\ov{\Q}_l$-sheaf on $V\times_S V$ which
is the pullback of the Artin---Schreier sheaf $\LL_{\psi}$ on $\G_a$ under
the morphism  $\langle  ,  \rangle$.
Let $s:S\ra V$ be  the  zero  section,  $j:U\ra  V$  the
complementary open subset to $s(S)$, and $p=\pi\circ j:U\ra S$
the projection of $U$ to $S$. Let us denote
$$\FF_!=j^*\FF j_!:\D^b_c(U,\ov{\Q}_l)\ra \D^b_c(U,\ov{\Q}_l).$$
                                                                     
\begin{prop}\label{Fouriersq}
For every $K\in \D^b_c(U,\ov{\Q_l})$
there    is    a   canonical   distinguished   triangle   in
$\D^b_c(U,\ov{\Q}_l)$:
\begin{equation}\label{triangle}
\FF_!^2(K)\ra K\ra p^*p_! K[4n](2n)\ra\ldots
\end{equation}
\end{prop}
        
\Pf . Note that $\FF_!(K)\simeq p_{2!}(\LL'\ot p_1^*K)[2n](n)$
where $\LL'=\LL|_{U^2}$;
hence,
$$\FF_!^2(K)\simeq p_{2!}((\LL'\circ\LL')\ot p_1^*K)[4n](2n)$$
where                  $\LL'\circ\LL'=p_{13!}(p_{12}^*\LL'\ot
p_{23}^*\LL')\in\D^b_c(U^2,\ov{\Q}_l)$.
Let $k=\id\times j\times\id: V\times_S U\times_S V\ra V^3$ be
the open embedding. Then
$$\LL'\circ\LL'\simeq (j\times j)^*p_{13!}
k_!k^*(p_{12}^*\LL\ot p_{23}^*\LL).$$
Note that we have
a canonical isomorphism
$$\LL_{\psi}(\langle    x_1, x_2\rangle)\ot\LL_{\psi}(\langle
x_2, x_3\rangle)\simeq \LL_{\psi}(\langle x_2, x_3-x_1\rangle)$$
on $V^3$, and the latter rank-1 local system is  trivial  on
the  complement  to $k$; hence we have the following exact
triangle:
$$\LL'\circ\LL'\ra (j\times j)^*(\LL\circ\LL)\ra
\ov{\Q}_{l,U^2}\ra\ldots,$$
where $\LL\circ\LL=p_{13!}(p_{12}^*\LL\ot p_{23}^*\LL)\simeq
\De_*\ov{\Q}_{l,V}[-4n](-2n)$, which induces the exact triangle
(\ref{triangle}).
\ed
                                 
Let us denote by
$\sideset{^p}{_!}{\FF}=\sideset{^p}{^0}{H}\FF_!:\Perv(U)\ra\Perv(U)$
the  right-exact  functor  induced  by $\FF_!$. Then we have a
canonical morphism of functors
$\nu:(\sideset{^p}{_!}{\FF})^2\ra\Id$.
        
\begin{cor}\label{zer}
The full subcategory of objects  $K\in\Perv(U)$  such
that $\nu_K:(\sideset{^p}{_!}{\FF})^2K\ra K$ is zero coincides with
the subcategory $p^*[2n](\Perv(S))\sub\Perv(U)$.
\end{cor}
        
\Pf  .  For $K\in\Perv(U)$ the exact
triangle (\ref{triangle}) induces the long exact sequence
$$\ldots\ra(\sideset{^p}{_!}{F})^2 K\ra K\ra
p^*[2n](\sideset{^p}{^{2n}}{H}p_!K)(2n)\ra 0.$$
Since $\sideset{^p}{^{4n}}{H}p_!p^*(2n)\simeq\Id_{\Perv(S)}$,
the assertion follows.
\ed
       
\subsection{Associativity}\label{assofourier}
To check the associativity condition (\ref{asso_simple})
for the above morphism $\FF_{!}^2\ra\Id$,
it is sufficient to check this condition for the
morphism  $\FF^2\ra\Id$,  which  is  done  in  the  following
proposition.
       
\begin{prop} Let $\FF$ be the Fourier  transform  associated
with  a  symplectic  vector  bundle  $V$ of rank $2n$. Then the
canonical morphism of functors $c:\FF^2\ra\Id$ satisfies the
associativity condition $\FF\circ c=c\circ\FF:\FF^3\ra\FF$.
\end{prop}
        
\Pf . Recall that $\FF$ is given by the kernel $\LL[2n](n)$ on
$V^2$ and $c$ is induced by the canonical morphisms of kernels
\begin{align*}
&p_{13!}(p_{12}^*\LL\ot p_{23}^*\LL[4n](2n))\ra
p_{13!}(\De_{13*}\De_{13}^*(p_{12}^*\LL\ot
p_{23}^*\LL[4n](2n)))\simeq\\
&\simeq p_{13!}\De_{13*}\ov{\Q}_{l,V^2}[4n](2n)
\stackrel{\Tr}{\ra}\De_*\ov{\Q}_{l,V},
\end{align*}
where $p_{ij}$ are the projections of $V^3$  on  the  double
products, $\De:V\ra V^2$ and $\De_{13}:V^2\ra V^3$ are the diagonals,
$\Tr$ is induced by the trace morphism (see \cite{SGA4})
(the triviality of $\De_{13}^*(p_{12}^*\LL\ot p_{23}^*\LL)$
follows from the skew-symmetry of $\LL$).
Now $\FF^3$ is given by the kernel
$p_{14!}(L[6n](3n))$ on
$V^2$  where $L=p_{12}^*\LL\ot p_{23}^*\LL\ot p_{34}^*\LL$
$p_{ij}$  are  the projections of $V^4$ on the
double products, and the associativity boils down to the
commutativity of the diagram
\begin{equation}
\begin{array}{ccccc}
& &p_{14!}(L[6n](3n))& &\\
&\ldlar{} & &\ldrar{}&\\
p_{14!}(\De_{13*}\De_{13}^*L[6n](3n))& &
& & p_{14!}(\De_{24*}\De_{24}^*L[6n](3n))\\
&\ldrar{\Tr} & &\ldlar{\Tr}& \\
& &\LL[2n](n)& &
\end{array}
\end{equation}
where $\De_{ij}:V^3\ra V^4$ are the diagonals,
and the lower diagonal arrows are induced by the isomorphisms
$\De_{13}^*L\simeq\De_{24}^*L\simeq p_{13}^*\LL$ on $V^3$
and by the trace morphisms.
Changing  the  coordinates  on  $V^4$  by  $(x,y,z,t)\mapsto
(x,y-t,z-x,t)$   one  can  see  that  this  reduces  to  the
commutativity of the diagram
\begin{equation}\label{comm_tr}
\begin{array}{ccccc}
p_!(\LL[4n](2n)) &\lrar{} p_!(i_{1*}i_1^*\LL[4n](2n))\simeq &
p_!(\ov{\Q}_{l,V}[4n](2n))\\
\ldar{} & &\ldar{\Tr}\\
p_!(i_{2*}i_2^*\LL[4n](2n))&\simeq p_!(\ov{\Q}_{l,V}[4n](2n))
\lrar{\Tr} & \ov{\Q}_l
\end{array}
\end{equation}
where $i_1:V\times 0\ra V^2$, $i_2:0\times  V\ra  V^2$  are  the
natural  closed  embeddings,   and $p$  denotes the projection to
$\Spec(k)$.
Here  is the argument  due  to  M.~Rapoport  verifying  the
commutativity of (\ref{comm_tr}). Obviously, we may assume
that the base is a point and $n=1$, so that $V=\A^2$.
Let   $i:I\hookrightarrow   \P^1\times   \A^2$   be   the
tautological  line  bundle  over  the  projective  line (the
incidence correspondence), and $p_2^*\LL$ be  the  pullback  of
$\LL$ under the projection $p_2:\P^1\times\A^2\ra\A^2$.
Then  $i$  is  an embedding of a lagrangian subbundle in the
trivial rank-2 symplectic bundle over $\P^1$; hence,
we have a sequence of canonical morphisms
$$q^*p_!(\LL[4](2))\wt{\ra}
\pi_!(p_2^*\LL[4](2))\ra
\pi_!(i_*i^*p_2^*\LL[4](2))\wt{\ra}
\pi_!(\ov{\Q}_{l,I}[4](2))\wt{\ra} q^*\ov{\Q}_l$$
where $\pi$ denotes the projection to  $\P^1$, and  $q$  is  the
projection of $\P^1$ to  $\Spec(k)$. Let
$\phi:q^*p_!(\LL[4](2))\ra q^*\ov{\Q}_l$ be the composed
morphism.   Then   $\phi=q^*\wt{\phi}$   for  some  morphism
$\wt{\phi}:p_!(\LL[4](2))\ra\ov{\Q}_l$   (since   in    fact
$p_!(\LL[4](2))\simeq\ov{\Q}_l$). On the other hand, the two morphisms
$p_!(\LL[4](2))\ra\ov{\Q}_l$ in (\ref{comm_tr})
are the restrictions of $\phi$ to the points $x,y\in\P^1$
corresponding to the coordinate lines in $\A^2$; hence, they
are both equal to $\wt{\phi}$.
\ed

\section{Gluing on the basic affine space}\label{glu}
       
\subsection{Setup}
Let  $k$ be a field of characteristic $p>2$, which is either
finite or algebraically closed, let
$G$ be a connected, simply-connected, semisimple algebraic group over
$k$. Assume that $G$ is split over $k$
and fix a split maximal torus $T\sub G$ and a Borel subgroup
$B$ containing $T$. Also we denote by $W=N(T)/T$ the Weil group of $G$,
and by $S\sub W$ the set of simple reflections.
Let $X=G/U$ be the corresponding
basic affine space, where
$U$  is  the unipotent radical of $B$.
Following \cite{KL} we are going to construct a $W$-gluing
data such that $\CC_w$ is
the category $\Perv(X)$ of perverse sheaves on $X$ for every
$w$.
To define the gluing functors we need some additional data.
First, we fix a nontrivial additive character
$\psi:\F_p\ra\ov{\Q_l}^*$.  We  denote  by  $\LL_{\psi}$ the
corresponding Artin---Schreier sheaf on $\G_{a,\F_p}$.
Second, for every simple root $\a_s$ we fix an isomorphism
of the corresponding 1-parameter subgroup $U_s\sub U$ with
the additive  group  $\G_{a,k}$.  This  defines  uniquely  a
homomorphism $\rho_s:\SL_{2,k}\ra G$ such that
$\varphi_s$ induces the given isomorphism of $\G_{a,k}$
(embedded in $\SL_{2,k}$ as upper-triangular matrices)
with $U_s$ (see \cite{KL}). Let
$$n_s=\rho_s\left(\matrix 0 & 1 \\ -1 & 0 \endmatrix\right).$$
For every $w\in W$ with a reduced decomposition
$w=s_1\ldots s_l$ we set $n_w=n_{s_1}\ldots n_{s_l}\in G$.
Then  $n_w$  does not  depend  on  a  choice  of  a   reduced
decomposition, so we get a canonical system of representatives
for $W$ in $N(T)$.
       
\subsection{Gluing functors}  For  every  element  $w\in W$ we consider the
subtorus
$$T_w=\prod_{\a\in R(w)}{\a^{\vee}(\G_m)}\sub T$$
where $R(w)\sub R^+$ is the set of positive roots $\a\in R^+$
such that $w(\a)\in -R^+$, $\a^{\vee}$ is the corresponding
coroot. By definition we have a surjective homomorphism
$$\prod_{\a\in R(w)}{\a^{\vee}}:\G_{m,k}^{R(w)}\ra T_w.$$
It is easy to see that
$$T_w=\prod_{s\in S_w} T_s$$
where $S_w\sub S$ is the set of simple reflections $s$ such
that  $s\leq  w$  with  respect  to  the  Bruhat  order (see
\cite{KL}, 2.2.1).
Now we define $X(w)\sub X\times_k X$ as the subvariety
of pairs $(gU,g'U)\sub X\times_k X$ such that
$g^{-1}g'\in Un_wT_wT$. There is a canonical projection
$\pr_w:X(w)\ra T_w$ sending $(gU,g'U)$ to the unique
$t_w\in T_w$ such that $g^{-1}g'\in Un_wt_wU$.
The morphism $\pr_w$ is smooth of relative dimension
$\dim X +\ell(w)$, surjective, with connected geometric fibers.
The last ingredient in the definition of the gluing functors is
the morphism
$$\si_w:\G_{m,k}^{R(w)}\ra\G_{a,k}:
(z_{\a})_{\a\in R(w)}\mapsto -\sum_{\a\in R(w)} z_{\a}.$$
Now we set
$$K(w)=K_{\psi}(w)=
\pr_w^*(\prod_{\a\in R(w)}\a^{\vee})_!\si_w^*\LL_{\psi}
[2\ell(w)](\ell(w)).$$
As shown in (\cite{KL}, 2.2.8) this is, up to shift by $\dim X$,
an irreducible
perverse sheaf on $X(w)$. Finally, one defines
$\ov{K(w)}$ to be the Goresky---MacPherson extension of
$K(w)$ to the closure $\ov{X(w)}$ of $X(w)$ in $X\times X$.
The gluing functors
$F_{w,!}:\D^b_c(X,\ov{\Q}_l)\ra\D^b_c(X,\ov{\Q}_l)$
are defined by 
$$F_{w,!}(A)=p_{2,!}(p_1^*(A)\ot \ov{K(w)}).$$
     
In the case when $w=s$ is a simple reflection
the morphism $\pr_s:X(s)\ra T_s\simeq\G_{m,k}$ extends to
$\ov{\pr}_s:\ov{X(s)}\ra\G_{a,k}$ and we have
$$\ov{K(s)}\simeq (-\ov{\pr}_s)^*\LL_{\psi}.$$
This leads to the alternative construction of the functor
$F_{s,!}$ using the embedding of
$G/U$ in a rank-2 vector bundle and  the  corresponding
partial Fourier transform. Namely,
let us denote $M_s=\rho_s(\SL_{2,k})\sub G$ and
consider   the   projection   $p_s:X=G/U\ra    G/Q_s$,    where
$Q_s=M_sU\sub G$ (note that $U$ normalizes $M_s$).
Now  $p_s$ is the complementary open subset
to the zero section in a $G$-equivariant rank-2 symplectic vector bundle
$\pi_s:V_s\ra G/Q_s$ (see \cite{KL}). Furthermore, we have
$\ov{X(s)}\simeq X\times_{G/Q_s} X$ and the morphism
$-\ov{\pr}_s$ from $\ov{X(s)}$ to $\G_{a,k}$ coincides  with
the restriction of the symplectic pairing on $V_s$. It follows
that $F_{s,!}=j^*\FF j_!$ where $j:X\hra V_s$ is the embedding,
$\FF$ is the (symplectic) Fourier transform for  $V_s$  (see
the previous section),  so   $\FF^2\simeq\Id$.  Since  $j_!$  is right
$t$-exact with respect to  perverse  $t$-structures,  so  is
$F_s$. As is shown in \cite{KL} (see also \ref{recKL} below)
for every reduced decomposition $w=s_1\ldots s_l$ one has  a
canonical isomorphism of functors
$$F_{w,!}\simeq F_{s_1,!}\circ\ldots\circ F_{s_l,!}.$$
It  follows  that  all  the  functors  $F_{w,!}$  are  right
$t$-exact.
       
\begin{thm} The functors $F_{w,!}$ define a quasi-action of
$W$ on $\D^b_c(X,\ov{\Q}_l)$.
\end{thm}
       
It follows that the functors $\sideset{^p}{^0}{H} F_{w,!}$ can be used
to define a $W$-gluing for $|W|$ copies of $\Perv(X)$.
The resulting glued category is denoted by $\AA$.
This theorem was stated as Theorem 2.6.1 in \cite{KL}.
However, the proof presented in loc. cit. is insufficient.
Namely, it is proved in \cite{KL} that the functors
$F_{w,!}$ generate the action of the positive braid monoid
on  $\D^b_c(X,\ov{\Q}_l)$,  and  the  morphisms  $F_{s,!}^2\ra\Id$   are
constructed. To finish the proof one has to show
that the other conditions of Theorem \ref{quasiact}
are  satisfied.   The   condition   (1)   was   checked   in
\ref{assofourier},  and  the  condition  (2) will be checked
below in \ref{assocomp}.
       
\subsection{Isomorphism}\label{recKL} 
Let us recall from \cite{KL} the construction of
the canonical isomorphism
$$F_{w,!}\simeq F_{w_1,!}\circ F_{w_2,!}$$
associated with a decomposition $w=w_1w_2$ such that
$\ell(w)=\ell(w_1)+\ell(w_2)$.
One starts with a commutative diagram
\begin{equation}\label{Xw_1w_2cart}
\begin{array}{ccccc}
X(w_1)\times_X X(w_2) &\lrar{} &T_{w_1}\times T_{w_2}\\
\ldar{} & & \ldar{m_{w_1,w_2}}\\
X(w) &\lrar{} &T_w
\end{array}
\end{equation}
where the map
$m_{w_1,w_2}:T_{w_1}\times T_{w_2}\ra T_{w_1w_2}$ is given
by $(t_1,t_2)\mapsto w_2^{-1}(t_1)t_2$.
It is easy to check that this diagram is cartesian. Moreover,
since $R(w_1w_2)$ is the disjoint union of  $w_2^{-1}(R(w_1))$
and $R(w_2)$, by the base change we get an isomorphism
\begin{equation}\label{w_1w_2}
K(w)\simeq p_{13,!}(p_{12}^*K(w_1)\otimes p_{23}^*K(w_2)),
\end{equation}
where $p_{ij}$ are the projections from
$X(w_1)\times_X X(w_2)\sub X^3$.
To derive from this the isomorphism
\begin{equation}\label{cw_1w_2}
\ov{K(w)}\simeq \ov{K(w_1)}\circ \ov{K(w_2)}
\end{equation}
it remains to check that the right-hand side is an irreducible
perverse sheaf up to shift by $\dim X$.
This is done in \cite{KL}, (2.4)---(2.5).
                                           
\subsection{Geometric information}\label{Xw_1w_2}
We need some more information about varieties
$X(w)$.
For every $w_1, w_2\in W$ such that $l(w_1w_2)=l(w_1)+l(w_2)$
consider the natural maps
\begin{align*}
&i_{w_1,w_2}:  X(w_1)\times_X   X(w_2)\ra   X(w)\times
T_{w_1}:
(g_1U, g_2U, g_3U)\mapsto ((g_1U, g_3U), \pr_{w_1}(g_1U,g_2U),\\
&i'_{w_1,w_2}:  X(w_1)\times_X   X(w_2)\ra   X(w)\times
T_{w_2}:
(g_1U, g_2U, g_3U)\mapsto ((g_1U, g_3U), \pr_{w_2}(g_2U,g_3U),
\end{align*}
where $w=w_1w_2$.
     
\begin{lem}\label{iw_1w_2}
The morphism $i_{w_1,w_2}$ (resp. $i'_{w_1,w_2}$) is an isomorphism onto
the  locally  closed subvariety of pairs
$(x,t_1)\in X(w)\times
T_{w_1}$ such that $\pr_{w}(x)t_1^{-1}\in T_{w_2}$
(resp. $(x,t_2)\in X(w)\times
T_{w_2}$ such that $\pr_{w}(x)t_2^{-1}\in
w_2^{-1}(T_{w_1})$).
\end{lem}
     
\Pf . The cartesian diagram (\ref{Xw_1w_2cart}) allows one to
identify $X(w_1)\times_X X(w_2)$ with the subvariety
in $X(w)\times T_{w_1}\times T_{w_2}$ consisting
of $(x,t_1,t_2)$ such that $\pr_w(x)=w_2^{-1}(t_1)t_2$.
Together with the fact that $w_2(t)t^{-1}\in T_{w_2}$
for any $t\in T$, this implies the assertion.
\ed
     
Let    $\BB=G/B$    be    the    flag    variety   of   $G$,
and $O(w)\sub\BB\times\BB$ the $G$-orbit corresponding to
$w\in  W$.  The  canonical  projection  $X(w)\ra  O(w)$ can
be considered as a $T\times T_w$-torsor. Namely, we have the
natural action of $T\times T$ on $X\times X$ such that
$(t,t')(gU,g'U)=(gtU,g't'U)$. Now the subgroup
$\{(t,t')\  |\  t^{-1}t'\in  T_w  \}$  (which  is  naturally
isomorphic to $T\times T_w$) preserves  $X(w)$  and  induces
the  above torsor structure. Let us denote by $\ov{X(w)}\sub
X\times  X$  and  $\ov{O(w)}\sub\BB\times\BB$  the   Zariski
closures.
       
\begin{lem}\label{tors} For every $s\in S$ the morphism
$\ov{X(s)}\ra \ov{O(s)}$ is a $T\times T_s$-torsor.
\end{lem}
       
This follows essentially from the rank-1 case when
$\ov{X(s)}=(\A^2-0)^2$, $\ov{O(s)}=\P^1\times\P^1$.
       
\subsection{Associativity}\label{assocomp}
Now we will check the second condition of Theorem \ref{quasiact}
for the functors $F_{w,!}$. Clearly, it is sufficient to
check the commutativity of the corresponding diagram of
morphisms between kernels on $X\times X$.
Let us recall the construction of the morphism
$$c_{s,s}:\ov{K(s)}\circ\ov{K(s)}\ra\De_*\ov{\Q}_{X,l}.$$
Consider the natural embedding
\begin{equation}\label{De_s}
\De_s:\ov{X(s)}\hra\ov{X(s)}\times_X\ov{X(s)}:
(g_1U,g_2U)\mapsto (g_2U,g_1U,g_2U).
\end{equation}
We have the morphism
$$\ov{\pr}_{s,s}:\ov{X(s)}\times_X\ov{X(s)}\ra \G_{a,k}:
(x,x')\mapsto\ov{\pr}_s(x)+\ov{\pr}_s(x')$$
such that $p_{12}^*\ov{K(s)}\otimes p_{23}^*\ov{K(s)}=
(-\ov{\pr}_{s,s})^*\LL_{\psi}$. The composition
$\ov{\pr}_{s,s}\De_s$ is the constant map to
$\{0\}\in\G_{a,k}$, hence we get the canonical isomorphism
$$\De_s^*(p_{12}^*\ov{K(s)}\otimes p_{23}^*\ov{K(s)})
\simeq\ov{\Q}_{\ov{X(s)},l}[2](1).$$
Now $c_{s,s}$ corresponds by adjunction to the morphism
$$\De^*(\ov{K(s)}\circ\ov{K(s)})\simeq
p_{1,!}(\De_s^*(p_{12}^*\ov{K(s)}\otimes p_{23}^*\ov{K(s)})\simeq
p_{1,!}(\ov{\Q}_{\ov{X(s)},l}[2](1))\stackrel{\tr}{\ra}
\ov{\Q}_{X,l}$$
where $\tr$ is the relative trace morphism for $p_1$.
       
\begin{thm}\label{sizker}
Let $sw=ws'$, where
$w\in W$, $s,s'\in S$ and $l(sw)=l(w)+1$.
Then the following diagram in $\D^b_c(X\times X,\ov{\Q}_l)$ is
commutative:
\begin{equation}\label{sizigker}
\begin{array}{ccccc}
\ov{K(s)}\circ\ov{K(w)}\circ\ov{K(s')} &\lrar{} &
\ov{K(w)}\circ\ov{K(s')}\circ\ov{K(s')}\\
\ldar{} & & \ldar{c_{s',s'}}\\
\ov{K(s)}\circ\ov{K(s)}\circ\ov{K(w)} &\lrar{c_{s,s}} &
\ov{K(w)}
\end{array}
\end{equation}
where  unmarked  arrows  are  induced  by  the  isomorphisms
(\ref{cw_1w_2}).
\end{thm}
                                     
The main step in the proof is the following lemma.
       
\begin{lem}\label{mainsiz} Under the conditions of Theorem \ref{sizker}
there is a canonical isomorphism
$$\a:\ov{X(s)}\times_X X(w)\ra X(w)\times_X\ov{X(s')}$$
of schemes over $X\times X$ and a canonical isomorphism
$$\ov{p}_s^*\ov{K(s)}\otimes p_w^*K(w)\simeq
\a^*(p_w^*K(w)\otimes \ov{p}_{s'}^*\ov{K(s)})$$
which is compatible with the isomorphisms (\ref{w_1w_2}),
where e.g. $\ov{p}_s:\ov{X(s)}\times_X X(w)\ra \ov{X(s)}$
is the natural projection, etc.
\end{lem}
       
\Pf . Consider the morphism
$$\ov{i'}_{s,w}:\ov{X(s)}\times_X X(w)\ra
X\times X\times T_w:
(g_1U, g_2U, g_3U)\mapsto (g_1U, g_3U, \pr_w(g_2U,g_3U))$$
extending the morphism $i'_{s,w}$ defined in \ref{Xw_1w_2}.
It is easy to see that $\ov{i'}_{s,w}$ is an embedding.
Let us consider the locally closed subvariety
$Y(s,w)=O(sw)\cup  O(w)\sub\BB\times\BB$.
We have the following isomorphism of
$\BB\times\BB$-schemes:
$$\ov{O(s)}\times_{\BB} O(w)\simeq Y(s,w).$$
Using Lemma \ref{tors} we see that the natural projection
$$\ov{p}_{s,w}:\ov{X(s)}\times_X X(w)\ra Y(s,w)$$
is   a   $T\times   T_s\times  T_w$-torsor,  where
$T\times T_s\times T_w$ is identified with the subgroup
$$\{(t_1,t_2,t_3)\in T\times T\times T|\ t_1^{-1}t_2\in T_s,
t_2^{-1}t_3\in T_w\}$$
acting naturally on $\ov{X(s)}\times_X X(w)$.
Let $p:X\times X\times T_w\ra\BB\times\BB$  be  the projection.
Then $p$ is a $T\times T\times T_w$-torsor, and
via $\ov{i'}_{s,w}$ we can identify the morphism
$$p^{-1}(Y(s,w))\ra Y(s,w)$$
with the $T\times T\times T_w$-torsor over $Y(s,w)$ induced from
$\ov{p}_{s,w}$ by the embedding $T\times T_s\times T_w\hra
T\times   T\times   T_w$.   In   particular,  the image  of
$\ov{i'}_{s,w}$ is a closed subvariety in $p^{-1}(Y(s,w))$.
Thus, $\im(\ov{i'}_{s,w})$ is the closure of $\im(i'_{s,w})$
in $p^{-1}(Y(s,w))$.
       
Similarly, we have an embedding
$$\ov{i}_{w,s'}:X(w)\times_X \ov{X(s)}\ra
X\times X\times T_w$$
extending $i_{w,s'}$, such that $\im(\ov{i}_{w,s'})$
is the closure of $\im(i_{w,s'})$ in $p^{-1}(Y(w,s'))$, where
$Y(ws')=O(ws')\cup O(w)$. Now since $sw=ws'$ we have
$Y(s,w)=Y(w,s')$ and we claim that
$\im(i'_{s,w})=\im(i_{w,s'})$, which implies immediately
that $\im(\ov{i'}_{s,w})=\im(\ov{i}_{w,s'})$. Indeed,
according to Lemma \ref{iw_1w_2} the image of $i'_{s,w}$
consists of pairs $(x,t)\in X(sw)\times T_w$
such that $\pr_{sw}(x)t^{-1}\in w^{-1}(T_s)$, while
the image of $i_{w,s'}$ consists of
$(x,t)$ such that $\pr_{sw}(x)t^{-1}\in T_s'$. Now
the equality $w^{-1}sw=s'$ implies that $w^{-1}(T_s)=T_{s'}$
which proves our claim.
                                   
Let
$$\a:\ov{X(s)}\times_X X(w)\ra X(w)\times_X\ov{X(s')}$$
be the unique isomorphism compatible with embeddings
$i'_{s,w}$ and $i_{w,s'}$.
It     remains     to     check     that  the     sheaves
$\ov{p}_s^*\ov{K(s)}\otimes p_w^*K(w)$ and
$p_w^*K(w)\otimes\ov{p}_{s'}^*\ov{K(s')}$ correspond to each
other under $\a$. Since $K(w)$ is the inverse image of a sheaf
on $T_w$ it is sufficient to check that the sheaves
$\ov{p}_s^*\ov{K(s)}$ and $\ov{p}_{s'}^*\ov{K(s')}$ correspond
to each other under $\a$. Since both these sheaves are local
systems (up to  shift) it is sufficient to check that
$p_s^*K(s)$ and $p_{s'}^*K(s')$ correspond to each other
under the restriction of $\a$ to the open subset $X(s)\times_X X(w)$.
It is easy to check that the following diagram is commutative
\begin{equation}
\begin{array}{ccccc}
X(s)\times_X X(w) &\lrar{\a} &
X(w)\times_X X(s')\\
\ldar{\pr_s} & & \ldar{\pr_s'}\\
T_s &\lrar{w^{-1}} & T_{s'}
\end{array}
\end{equation}
Moreover, since $l(sw)=l(w)+1$ we have
$w^{-1}(\a_s^{\vee})=\a_{s'}^{\vee}$, hence
the bottom arrow becomes the identity under the
identification of $T_s$ (resp. $T_{s'}$) with $\G_m$
via $\a_s^{\vee}$ (resp. $\a_{s'}^{\vee}$),
and our assertrion follows immediately.
\ed
       
\begin{lem}\label{cor}
Let $Y$ be a scheme, $A$ a correspondence
over $Y\times Y$, $C$ and $C'$ symmetric
correspondences over $Y\times Y$ such that an isomorphism
of $Y\times Y$-schemes is given by
$$\a:C\times_Y A\wt{\ra} A\times_Y C'.$$
Let
\begin{align*}
&\De_C:C\hra C\times_Y C:(y_1,y_2)\mapsto (y_2,y_1,y_2),\\
&\De_{C'}\si:C'\hra C'\times_Y C':(y_1,y_2)\mapsto (y_1,y_2,y_1)
\end{align*}
be  the  natural  embeddings  (where  $\si:C'\ra  C'$ is the
permutation of factors in $Y\times Y$). Then the following diagram
is commutative:
\begin{equation}
\begin{array}{ccccc}
C\times_Y A &\lrar{\a'} &
A\times_Y C'\\
\ldar{\De_C} & & \ldar{\De_{C'}\si}\\
C\times_Y C\times_Y A &\lrar{} & A\times_Y C'\times_Y C'
\end{array}
\end{equation}
where the bottom arrow is induced by $\a$,
$\a'$ is an isomorphism given by
$$\a'(y_1,y_2,y_3)=(y_2,y_3,y_2'),$$
for $(y_1,y_2,y_3)\in C\times_Y A$ where
$\a(y_1,y_2,y_3)=(y_1,y_2',y_3)$.
\end{lem}
       
The  proof is straightforward. Note that $\a'$ commutes with
projections to $A$, and $\a^{-1}\circ\a'$ is an involution of
$C\times_Y A$.

\vspace{2mm}       

\noindent
{\it Proof of Theorem \ref{sizker}}.
First we note that
$\ov{K(s)}\circ\ov{K(w)}\circ\ov{K(s')}\in
\sideset{^p}{^{\leq\dim X}}{\D}(X\times X)$.
Indeed, the functor $K\mapsto \ov{K(s)}\circ K$
from $\D^b_c(X\times X,\ov{\Q}_l)$ to itself coincides
with the functor $\FF_!$ of (\ref{foursq}) for the
rank-2 symplectic bundle $V_s\times X\ra G/Q_s\times X$,
and hence it is right $t$-exact with  respect  to  the  perverse
$t$-structure. Similarly, the
functor $K\mapsto K\circ\ov{K(s')}$ is right $t$-exact,
hence our claim follows from the fact that
$\ov{K(w)}[\dim X]$ is a
perverse sheaf on $X\times X$. It follows that
we can replace all objects in the diagram (\ref{sizigker})
by their $\sideset{^p}{^{\dim X}}H$. Since
$\ov{K(w)}$ is the Goresky---MacPherson extension from
$X(w)\sub X\times_k X$, it is sufficient to check the
commutativity of the restriction of (\ref{sizigker}) to
$X(w)$. Applying Lemma \ref{cor} to the correspondences
$A=X(w)$, $C=\ov{X(s)}$ and $C'=\ov{X(s')}$ we obtain
the commutative diagram
\begin{equation}
\begin{array}{ccccc}
\ov{X(s)}\times_X X(w) &\lrar{\a'} &
X(w)\times_X \ov{X(s')}\\
\ldar{\De_s} & & \ldar{\De_{s'}\si}\\
\ov{X(s)}\times_X \ov{X(s)}\times_X X(w) &\lrar{} &
X(w)\times_X \ov{X(s')}\times_X \ov{X(s')}
\end{array}
\end{equation}
From the construction of the morphisms $c_{s,s}$ and
$c_{s',s'}$ and Lemma \ref{mainsiz}, we see that
it is sufficient to prove that the trivializations (up to
shift and twist)
of $\De_s^*(\ov{p}_{s,1}^*\ov{K(s)}\otimes
\ov{p}_{s,2}^*\otimes p_w^*K(w))$
and of
$\De_{s'}^*(\ov{p}_{s',1}^*\ov{K(s')}\otimes
\ov{p}_{s',2}\otimes p_w^*K(w))$
are compatible via the above commutative diagram with
$\a'$ and the isomorphism of Lemma \ref{mainsiz}
(here $\ov{p}_{s,1}$ and $\ov{p}_{s,2}$ are the projections onto the
first and the second factors $\ov{X(s)}$).
To this end we can replace $\ov{X(s)}$, $\ov{X(s')}$, $\ov{K(s)}$
and $\ov{K(s')}$ by $X(s)$, $X(s')$, $K(s)$ and $K(s')$, 
respectively. 
Now this follows immediately from the fact that the isomorphism
$\a:X(s)\times_X X(w)\simeq X(w)\times_X X(s')$ is compatible
with the projections to $T_s\times T_w\simeq T_w\times T_{s'}$
(see the proof of Lemma \ref{mainsiz}).
\ed
                                                   
\subsection{Grothendieck group of the glued category}
It is easy to check that every functor $F_{w,!}$ has the right
adjoint $F_{w,*}$. Indeed, $F_{w,!}$ is the composition
of the functors $F_{s,!}$ corresponding to simple reflections.
Now $F_{s,!}=j^*\FF j_!$ where $\FF$ is the Fourier transform,
and $j$ is an open embedding; hence it is left adjoint
to $j^*\FF j_*$. In fact, it is shown in \cite{KL} that
$F_{w,*}(A)=p_{2,*}(p_1^*A\otimes\ov{K(w)})$.
     
The Proposition \ref{simple} combined with Theorem \ref{main}
gives  a simple description of the Grothendieck group of the
abelian category $\AA$ resulting from gluing
on $G/U$. 
As we have seen in Corollary \ref{zer} the subgroup
$K_{sw,w}\sub K_0(\Perv(G/U))$ coincides with the image of
the natural (injective) homomorphism
$p_s^*:K_0(\Perv(G/Q_s))\ra K_0(\Perv(G/U))$. Hence,
we get the  following  description of $K_0(\AA)$.
        
\begin{thm}\label{K_0A} The group $K_0(\AA)$ is isomorphic
to the group
$$\KK=\{ (c_w)\in\oplus_{w\in W} K_0(\Perv(G/U)) \ | \
\phi_s c_w-c_{sw}\in  p_s^*(K_0(\Perv(G/Q_s))),  \  w\in W,
s\in S \},$$
where $\phi_s$ is an operator on $K_0(\Perv(G/U))$
induced by the partial Fourier transform $F_{s,!}$.
\end{thm}

\section{Cubic Hecke algebra}
        
\subsection{A property of the Fourier transform}
Let $\pi:V\ra S$ be a symplectic rank-2 bundle, let
$\FF:\D^b_c(V,\ov{\Q}_l)\ra\D_c^b(V,\ov{\Q}_l)$ be
the corresponding Fourier
transform, $j:U\ra V$ the complement to the zero section, and
$\FF_!=j^*\FF j_!:\D_c^b(U,\ov{\Q}_l)\ra\D_c^b(U,\ov{\Q}_l)$.
        
\begin{lem}\label{mainl}
There is a canonical isomorphism of functors
\begin{equation}\label{shift}
\FF_!\circ p^*\simeq p^*[1](1)
\end{equation}
where $p=\pi\circ j:U\ra S$.
\end{lem}
        
\Pf . We start with canonical isomorphisms
\begin{align*}
&\FF\circ s_*\simeq \pi^*[2](1),\\
&\FF\circ\pi^*\simeq s_*[-2](-1)
\end{align*}
where $s:S\ra V$ is the zero section (see \cite{Laumon}).
Now applying $\FF$ to the exact triangle
$$j_!p^*F\ra\pi^*F\ra s_*F\ra\ldots$$
and using the above isomorphisms  we  obtain  the 
exact triangle on $V$:
$$\FF j_! p^*F\ra s_*F[-2](-1)\ra\pi^*F[2](1)\ldots$$
Restricting to $U\sub V$ we get the required isomorphism.
\ed
                    
\subsection{Cubic relation}
For every scheme $Y$ we denote by $K_0(Y)=K_0(\D_c^b(Y,\ov{\Q}_l))$
the Grothendieck group of the category 
$\D_c^b(Y,\ov{\Q}_l)$. We define the action of the algebra
$\Z[u,u^{-1}]$ (where $u$ is an indeterminate) on $K_0(Y)$
by setting $u\cdot [F]=[F(-1)]$, where $F\mapsto F(1)$ is the Tate twist.
        
\begin{prop}  
Let $\phi:K_0(U)\ra K_0(U)$ be the operator induced by $\FF_!$.
Then $\phi$  satisfies the 
equation
\begin{equation}\label{cubic}
(\phi+u^{-1})(\phi^2-1)=0.
\end{equation}
\end{prop}
        
\Pf     .     From    the    previous    lemma    we    have
$(\phi+u^{-1})|_{\im(p^*)}=0$.
On the  other hand, from  Proposition \ref{Fouriersq} we have
$\im(\phi^2-1)\sub\im(p^*)$, hence the assertion.
\ed

\begin{cor}\label{s2-1} Let $R=\Z[u,u^{-1},(u^2-1)^{-1}]$. The submodule 
$K_0(G/Q_s)\otimes_{\Z[u,u^{-1}]} R\sub K_0(G/U)
\otimes_{\Z[u,u^{-1}]} R$ coincides with the image of the operator
$\phi^2-1$ on the latter $R$-module.
\end{cor}
        
Let $r:U\ra\P(V)$ be the projection to the projectivization
of $V$, $q:\P(V)\ra S$ the projection to the base, so that
$p=q\circ r$.
        
\begin{prop}\label{Fourproj} The following relation between operators
$K_0(\P(V))\ra K_0(U)$ holds:
\begin{equation}
\phi r^*= r^* - u^{-1}p^*q_!= r^* (\id - u^{-1}q^*q_!).
\end{equation}
\end{prop}
        
\Pf . By definition we have
$$\FF_!(r^*G)=p_{2!}(p_1^*r^*G\ot\LL|_{U^2})[2](1)\simeq
p'_{2!}(p_1^{\prime *}G\ot (r\times\id_U)_!(\LL|_{U^2}))[2](1)$$
where
$p_i$ (resp. $p'_i$) are the projections of $U^2$
(resp. $\P(V)\times U$) on its factors.
Note that $(r\times\id_U)_!(\LL|_{U^2})\simeq
((r\times\id_V)_!(\LL|_{U\times V}))|_{\P(V)\times U}$.
To compute the latter sheaf we decompose $r$ as follows:
$r=l\circ k$ where $k:U\hra I$ is the open embedding into
the incidence correspondence $I\sub\P(V)\times V$
($k:v\mapsto  (\langle  v\rangle,v)$),  $l:I\ra\P(V)$ is the
projection. Let $s:\P(V)\ra \P(V)\times 0\sub I$ be the zero
section. Then since  the  images  of  $k$  and  $s$  are
complementary to each other we have an exact triangle
\begin{equation}\label{trian}
(k\times\id_V)_!(\LL|_{U\times V})\ra\wt{\LL}\ra
(s\times\id_V)_*(s\times\id_V)^*\wt{\LL}\ra\ldots
\end{equation}
where $\wt{\LL}$ is the pullback of $\LL$ by the morphism
$I\times V\ra V\times V$ induced by the projection $I\ra V$.
Now we have $(l\times\id_V)_!(s_*s^*\wt{\LL})\simeq
\ov{\Q}_{l,\P(V)\times V}$ and
$(l\times\id_V)_!(\wt{\LL})\simeq \ov{\Q}_{l,I}[-2](-1)$.
Indeed, the latter isomorphism follows from the fact that
$(l\times\id_V)_!(\wt{\LL})$      is      supported       on
$I\sub\P(V)\times V$ and the trivialiaty of
$\wt{\LL}|_{(l\times\id_V)^{-1}(I)}$.
Applying the functor $(l\times\id_V)_!$ to the triangle
(\ref{trian}) and using these isomorphisms we obtain
the exact triangle
$$(r\times\id_V)_!(\LL|_{U\times V})\ra
\ov{\Q}_{l,I}[-2](-1)\ra \ov{\Q}_{l,\P(V)\times V}\ra\ldots$$
Note that the restriction of the projection $I\ra V$ to
$I\cap (\P(V)\times U)$ is an isomorphism. Hence, passing to
Grothendieck groups we obtain
$$[\FF_!r^*G]=[r^*G]-[p^*q_!G(1)]$$
as required.
\ed
                    
\subsection{Action of the cubic Hecke algebra}
Now we return to the situation of section \ref{glu}.
Let $\HH$ be the Hecke algebra defined as the quotient of
the group algebra with coefficients in $\Z[u,u^{-1}]$
where $u$ is indeterminate,
of the generalized braid group corresponding to $(W,S)$ by
the relations $(s+1)(s-u)=0$, $s\in S$.
Recall that there is an action of $\HH$ on  $K_0(G/B)$  such
that $\Z[u,u^{-1}]$ acts in the standard way (using Tate twist) and
the  action  of  $s$  is  given  by  the
correspondence $O(s)\sub (G/B)^2$:
$O(s)=\{ (gB,g'B)\ |\ g^{-1}g'\in BsB\}$. In  terms  of  the
projective bundle $q_s:G/B\ra G/M_sB$ associated with $s$
we have $T_s=q_s^*q_{s!}-\id$.
        
\begin{prop} The functors $F_{s,!}$, $s\in S$ extend to the action  on
$K_0(G/U)$ of the cubic Hecke algebra $\HH^c$ which is obtained
as the quotient of the group algebra $\Z[u,u^{-1}][B]$ of the braid group
$B$ corresponding to $(W,S)$ by the relations
$(s+u)(s^2-1)=0$, $s\in S$.
This action preserves $K_0(G/B)\sub K_0(G/U)$
and restricts to the standard action of the quadratic Hecke algebra
$\HH$ on $K_0(G/B)$ via the $\Z[u,u^{-1}]$-linear homomorphism
$$\HH^c\ra\HH:s\mapsto -s^{-1}.$$
\end{prop}
       
The proof reduces to a simple computation for the Fourier transform
on a symplectic bundle of rank 2.
                    
\section{Adjoint functors and canonical complexes for $W$-gluing}
\label{complex}
   
\subsection{Adjoint functors for parabolic subgroups}
Let $(W,S)$ be a finite Coxeter group, $J\sub S$  a subset,
and $W_J\sub W$ the subgroup generated
by simple reflections in $J$.
   
We will frequently use the following fact
(see  \cite{B},  IV,  Exercise   1.3): every left
(or right) $W_J$-coset
contains a unique element of minimal length.
Furthermore, an element $w\in W$ is the shortest element in
$W_Jw$ if and only if $l(sw)=l(w)+1$ for every $s\in J$,
and  in  this  case  we  have  $l(w'w)=l(w')+l(w)$ for every
$w'\in W_J$.
   
In particular, for any coset $W_Jx$ and
any $w\in W$ there exists the unique element $p_{W_Jx}(w)\in W_Jx$
such that $n_{W_Jx}(w):=wp_{W_Jx}(w)^{-1}$ has minimal possible length.
Namely, $n_{W_Jx}(w)$ is the shortest element in the coset
$wx^{-1}W_J$.
   
Let $(\CC_w,w\in W)$ be a collection of abelian
categories, $\Phi_w:\CC_{w'}\ra\CC_{ww'}$,
a $W$-gluing data, and $\CC(\Phi)$  the corresponding
glued category.
For  every  subset  $P\sub  W$ let us denote by $\Phi_P$ the
gluing data on categories $\CC_w$, $w\in P$, induced by $\Phi$.
We have natural restriction functors
$j_P^*:\CC(\Phi)\ra\CC(\Phi_P)$.
 
We claim that for every coset $W_Jx\sub W$
the functor $j^*_{W_Jx}$ has the left adjoint
functor $$j_{W_Jx,!}:\CC(\Phi_{W_Jx})\ra\CC(\Phi).$$
Indeed, let $A=(A_w,w\in W_Jx;\a_{w,w'})$
be an object of $\CC(\Phi_{W_Jx})$.
Set $j_{W_Jx,!}=(A_w,w\in W;\a_{w,w'})$
where $$A_w=\Phi_{n(w)}A_{p(w)},$$
$n(w):=n_{W_Jx}(w)$, $p(w):=p_{W_Jx}(w)$, and
the morphisms $\a_{w,w'}:\Phi_w(A_{w'})\ra A_{ww'}$
are defined as follows. Let us write
$p(ww')=w_1p(w')$ with $w_1\in W_J$. Then we have
$$wn(w')=n(ww')w_1$$
and $l(n(ww')w_1)=l(n(ww'))+l(w_1)$.
Hence, we can define $\a_{w,w'}$ as the composition
$$\Phi_w(A_{w'})=\Phi_w\Phi_{n(w')}(A_{p(w')})\ra
\Phi_{wn(w')}(A_{p(w')})
\simeq\Phi_{n(ww')}\Phi_{w_1}(A_{p(w')})
\setlength{\unitlength}{0.50mm}
\lrar{\Phi_{n(ww')}(\a_{w_1,p(w')})}
\Phi_{n(ww')}(A_{w_1p(w')})=A_{ww'}.
$$
One can easily check that $j_{W_Jx,!}$ is indeed an
object of $\CC(\Phi)$.
   
\begin{prop}\label{jWJx!}
The functor $j_{W_Jx,!}$ is left adjoint to
$j^*_{W_Jx}$.
\end{prop}
   
\Pf.  Let  $A=(A_w,w\in  W_Jx;\a_{w,w'})$  be  an  object of
$\CC(\Phi_{W_Jx})$, and $B=(B_w,w\in W;\b_{w,w'})$ an object
of $\CC(\Phi)$. We have an obvious map
\begin{equation}\label{adjmap}
\Hom_{\CC(\Phi)}(j_{W_Jx,!}A,B)\ra
\Hom_{\CC(\Phi_{W_Jx})}(A,j^*_{W_Jx}B).
\end{equation}
The inverse map is constructed as follows.
Assume that we are given a morphism
$f=(f_w,w\in W_Jx):A\ra j^*_{W_Jx}B$.
Then for every $w\in W$ we define the morphism
$\wt{f}_w:\Phi_{n(w)}A_{p(w)}\ra B$
as the composition
$$
\setlength{\unitlength}{0.35mm}
\Phi_{n(w)}A_{p(w)}\lrar{\Phi_{n(w)}(f_{p(w)})}
\Phi_{n(w)}B_{p(w)}\lrar{\b_{n(w),p(w)}} B_w.$$
It is easy to check that $\wt{f}=(\wt{f}_w,w\in W)$
is the morphism in $\CC(\Phi)$ between $j_{W_Jx,!}A$ and
$B$ and that the obtained map
$$\Hom_{\CC(\Phi_{W_Jx})}(A,j^*_{W_Jx}B)\ra
\Hom_{\CC(\Phi)}(j_{W_Jx,!}A,B):f\mapsto\wt{f}$$
is inverse to (\ref{adjmap}).
\ed
 
If we have an inclusion $W_Jx\sub P\sub W$, then the restriction
functor
$j_{W_Jx,P}^*:\CC(\Phi_P)\ra\CC(\Phi_{W_Jx})$ has the left adjoint
$$j_{W_Jx,P;!}:=j^*_Pj_{W_Jx,!}:\CC(\Phi_{W_Jx})\ra
\CC(\Phi_P)$$
(this follows from the proof of the above proposition).
Moreover, by construction the composition $j_{W_Jx,P}^*j_{W_Jx,P;!}$
is the identity functor on $\CC(\Phi_{W_Jx})$. Hence, 
we can apply Theorem \ref{braverman} to conclude that for any subset 
$P\sub W$ which is a union of subsets of the form $W_Jx$, the category 
$\CC(\Phi_P)$ is obtained by gluing from the categories 
$\CC(\Phi_{W_Jx})$ for $W_Jx\sub P$.
 
Let $J\sub K\sub S$. Then one has canonical isomorphisms
$$j^*_{W_Jx}\simeq j^*_{W_Jx,W_Kx}\circ j^*_{W_Kx},$$
$$j_{W_Jx,!}\simeq j_{W_Kx,!}\circ j_{W_Jx,W_Kx;!}.$$
Also one has the canonical morphism of functors
\begin{equation}\label{WJK}
j_{W_Jx,!}j^*_{W_Jx}\ra j_{W_Kx,!}j^*_{W_Kx}.
\end{equation}
               
Assume that every functor $\Phi_w$ has the left derived 
$L\Phi_w:\D^-(\CC_w')\ra\D^-(\CC_{ww'})$ and that these functors
satisfy $L\Phi_{w_1}\circ L\Phi_{w_2}\simeq L\Phi_{w_1w_2}$
when $\ell(w_1w_2)=\ell(w_1)+\ell(w_2)$.
The following proposition gives a sufficient condition for the existence of
the derived functor for $j_{W_Jx,!}$.
   
\begin{prop}\label{leftderived} Assume that for every $w\in W$
there is a family of objects
$\RR_w\sub \Ob\CC_w$ that are $\Phi_{w'}$-acyclic for every
$w'\in W$ and such that every object of $\CC_w$
can be covered by an object in $\RR_w$. Then the
functor $j_{W_Jx,!}$ has the left derived 
$Lj_{W_Jx,!}:\D^-(\CC(\Phi_{W_Jx}))\ra\D^-(\CC(\Phi))$ which
is left adjoint to the restriction functor
$j^*_{W_Jx}:\D^-(\CC(\Phi))\ra\D^-(\CC(\Phi_{W_Jx}))$.
Furthermore, one has an isomorphism of functors
$$j^*_w\circ Lj_{W_Jx,!}\simeq L\Phi_{n(w)}j^*_{p(w)}.$$
\end{prop}
   
\Pf . Let $\RR\sub\Ob\CC(\Phi_{W_Jx})$ be the family
of objects that are direct sums of objects of the form
$j_{y,W_Jx;!}R_y$ where $y\in W_Jx$, $R_y\in\RR_y$.
Clearly, every object of $\CC(\Phi_{W_Jx})$ can be covered by
an object in $\RR$.
We claim that $\RR$ is an adapted class of objects for the functor
$j_{W_Jx,!}$. Indeed, it suffices to prove that for every $w\in W$ the
object $\Phi_{p(w)y^{-1}}R_y$ is $\Phi_{n(w)}$-acyclic, where
$w=n(w)p(w)$ is the decomposition used in the definition of $j_{W_Jx,!}$.
Now we use the fact that
$$\ell(wy^{-1})=\ell(n(w)p(w)y^{-1})=\ell(n(w))+\ell(p(w)y^{-1}).$$
Therefore, by our assumption
$$\Phi_{wy^{-1}}R_y=L\Phi_{wy^{-1}}R_y\simeq L\Phi_{n(w)}\circ
L\Phi_{p(w)y^{-1}}R_y=L\Phi_{n(w)}(\Phi_{p(w)y^{-1}}R_y)$$
and our claim follows. Thus, the left derived
functor for $j_{W_Jx,!}$ exists and can be computed
using resolutions in $\RR$. The remaining assertions
can be easily checked using such resolutions.
\ed

Note that Theorem \ref{adapted} implies that the conditions
of the previous proposition are satisfied for gluing on the
basic affine space.

\subsection{Canonical complex}
Let us fix a complete order on $S$: $S=\{s_1,\ldots,s_n\}$.
Given an object $A\in\CC(\Phi)$ we construct a homological
coefficient system on the $(n-1)$-simplex $\De_{n-1}$ with values in
$\CC(\Phi)$. Namely, to a subset
$J=\{i_1<\ldots<i_k\}\sub [1,n]$ we assign the
object
$$A(J)=\oplus j_{W_{S-J}x,!}j^*_{W_{S-J}x}A$$
where the sum is taken over all right $W_{S-J}$-cosets.
For every inclusion $J\sub J'$ we have the canonical morphism
$A(J')\ra A(J)$ with components (\ref{WJK}).
Thus, we can consider the corresponding chain complex
\begin{equation}
C_{\cdot}(A):C_{n-1}=A([1,n])\ra\ldots\ra
C_1=\oplus_{|J|=2}A(J)\ra C_0=\oplus_{|J|=1}A(J)
\end{equation}
The sum of adjunction morphisms
$$C_0(A)=\oplus_{|J|=n-1,W_Jx\sub  W}j_{W_Jx,!}j^*_{W_Jx}A\ra
A$$
induces the morphism
\begin{equation}\label{H0}
H_0(C_{\cdot}(A))\ra A.
\end{equation}
It  turns  out that $C_{\cdot}(A)$ is almost a resolution of
$A$. To describe its homology we need to introduce the
functor $\iota:\CC(\Phi)\ra\CC(\Phi)$. For $A=(A_w;\a_{w,w'})$
we set $\iota A=(\Phi_{w_0}A_{w_0w},\wt{\a}_{w,w'})$
where $w_0$ is the longest element in $W$, the morphism
$\wt{\a}_{w,w'}:\Phi_w\Phi_{w_0}A_{w_0w'}\ra
\Phi_{w_0}A_{w_0ww'}$ is equal to the composition
$$\Phi_w\Phi_{w_0}A_{w_0w'}\simeq
\Phi_{w_0}\Phi_{w_0ww_0}A_{w_0w'}
\setlength{\unitlength}{0.50mm}
\lrar{\Phi_{w_0}(\a_{w_0ww_0,w_0w'})}\Phi_{w_0}A_{w_0ww'}.$$
Here we used the following identity in the braid group:
$$\tau(w)\tau(w_0)=\tau(w)\tau(w^{-1}w_0)\tau(w_0ww_0)=
\tau(w_0)\tau(w_0ww_0).$$
Note that for every $y\in W$ we have the natural morphism
$$\a_y:\iota A\ra j_{y,!}j^*_yA$$
with components
$$\Phi_{wy^{-1}}(\a_{yw^{-1}w_0,w_0w}):\Phi_{w_0}A_{w_0w}\ra
\Phi_{wy^{-1}}A_y.$$
It is easy to see that one has the canonical morphism
\begin{equation}\label{Hn-1}
\iota A\ra H_{n-1}(C_{\cdot}(A)),
\end{equation}
induced by the morphism
$$\iota A\lrar{((-1)^{l(y)}\a_y)} C_{n-1}(A)=\oplus_{y\in  W}
j_{y,!}j^*_yA.$$

\begin{thm}\label{homology}
One has $H_i(C_{\cdot}(A))=0$ for $i\neq 0,n-1$,
$H_0(C_{\cdot}(A))\simeq A$, and $H_{n-1}(C_{\cdot}(A))\simeq\iota A$.
\end{thm}
   
\Pf. Let us consider the complex
$\wt{C}_{\cdot}(A)$
obtained from $C_{\cdot}$ by adding the terms
$\wt{C}_n=\iota A$ and $\wt{C}_{-1}=A$ with additional differentials
induced by (\ref{H0}) and (\ref{Hn-1}).
By   definition   for  every  $w\in  W$  the  complex
$j_w^*\wt{C}_{\cdot}(A)$ looks as follows:
\begin{align*}
&\Phi_{w_0}A_{w_0w}\ra\oplus_{x\in W}\Phi_{wx^{-1}}A_x\ra
\oplus_{|J|=1,x\in W_J\backslash W}
\Phi_{n_{W_Jx}(w)}A_{p_{W_Jx}(w)}\ra\ldots\\
%\ra\oplus_{|J|=n-2,x\in W_J\backslash W}
%\Phi_{n_{W_Jx}(w)}A_{p_{W_Jx}(w)}
&\ra\oplus_{|J|=n-1,x\in W_J\backslash W}
\Phi_{n_{W_Jx}(w)}A_{p_{W_Jx}(w)}\ra A_w.
\end{align*}
Recall  that  as  $x$  runs  through  the  set   of   cosets
$W_J\backslash W$,
the element $n_{W_Jx}(w)$ runs through all the elements
$y\in W$ such that $l(ys)=l(y)+1$ for all $s\in J$.
Let us denote by $P_j\sub W$ the set of all $y\in W$
such that $l(ys_j)=l(y)+1$.
For $J\sub [1,n]$ we denote $P_J=\cap_{j\in J} P_j$.
Then we can rewrite the complex $j^*_w\wt{C}_{\cdot}(A)$
as 
\begin{align*}
&\Phi_{w_0}A_{w_0w}\ra\oplus_{y\in W}\Phi_{y}A_{y^{-1}w}\ra
\oplus_{j,y\in P_j}\Phi_yA_{y^{-1}w}\ra
\ldots\\
%\ra\oplus_{|J|=n-2,y\in P_J}\Phi_{y}A_{y^{-1}w}
&\ra\oplus_{|J|=n-1,y\in P_J}\Phi_{y}A_{y^{-1}w}\ra A_w.
\end{align*}
Consider the increasing filtration $F_0\sub F_1\sub\ldots$
on $j^*_w\wt{C}_{\cdot}(A)$ such that
$F_n$ contains only summands $\Phi_y A_{y^{-1}w}$ with
$l(y)\le n$. Then the differentials in
$\gr_{F}j^*_w\wt{C}_{\cdot}(A)$ are only $\pm\id$ or zeroes.
Note that $P_S=\{1\}$ while
$\cup_{j=1}^n P_j=W-\{w_0\}$.
Hence the complex $\gr_F j_w^*\wt{C}_{\cdot}(A)$ is acyclic,
and so is $\wt{C}_{\cdot}(A)$.
\ed
   
Note that we can attach the end of the complex
$C_{\cdot}(\iota A)$ to the beginning of the complex
$C_{\cdot}(A)$ via the map
$$C_0(\iota A)\ra \iota A\ra C_{n-1}(A)$$
to get the complex
$\wt{C}=\wt{C}_{\cdot}(A)$ with homologies
$H_0(\wt{C})=A$ and $H_{2n-1}(\wt{C})=\iota^2(A)$ (all other homologies
vahish).
The functor $\iota^2$  sends     the     object
$A=(A_w,\a_{w,w'})\in\CC(\Phi)$ to the object
$(\Phi_{\pi} A_w, \Phi_{\pi}(\a_{w,w'}))$
where $\pi=\tau(w_0)^2\in B$ is the canonical central element.
   
The importance of the above construction is that the members
of the complex $C_{\cdot}$ are direct sums of objects of the
form $j_{W_Jx,!}(\cdot)$ where $J\sub S$ is a proper subset.
Hence, it can be used for the induction process.
For example, suppose we know that for all
proper subsets $J\sub S$ the  categories  $\CC(\Phi_{W_Jx})$
have finite cohomological dimension (this is true for
the gluing on the basic affine space if the rank of $G$ is equal to $2$).
The derived category
version of the above construction gives a canonical
morphism of functors $\Id\ra\Phi_{\pi}[2n]$. Now
an object $A\in\CC(\Phi)$ has finite projective dimension
if and only if some power of the morphism
$A\ra\Phi_{\pi}^k(A)[2nk]$ vanishes.
       
\subsection{Canonical complex for a ``half" of $W$}
                                  
Let us fix an element $s_i\in S$. Recall that we denote
by $P_i$ the subset of $w\in W$ such that $l(ws_i)>l(w)$.
One has a decomposition of $W$ into the disjoint union
of $P_i$ and $P_is_i$. Let us consider gluing data
$\Phi_{P_i}$ and $\Phi_{P_is_i}$
corresponding to these two pieces.
We are going to construct an analogue of the complex
$C_{\cdot}$ for these partial gluing data.
As before, for every $A\in\CC(\Phi_{P_i})$
we can define a homological
coefficient system on $\De_{n-1}$ with values in
$\CC(\Phi_{P_i})$ as follows.
For every $j\in [1,n]$ let us denote
$W^{(j)}=W_{[1,n]-j}$.
Now to a nonempty subset
$J=\{i_1<\ldots<i_k\}\sub [1,n]$ our coefficient system assigns the
object
$$A(J)=\oplus j_{W_{S-J}x,P_i;!}
j^*_{W_{S-J}x,P_i}A$$
where  the  sum  is  taken  over  all
$x\in W_{S-J}\backslash W$ such that
$W^{(j)}x\sub P_i$ for every $j\in J$.
For every inclusion $J\sub J'$ we have the canonical morphism
$A(J')\ra A(J)$. Let us consider the corresponding chain complex
\begin{equation}\label{chaincomplex}
C_{\cdot}(P_i,A):C_{n-1}=A([1,n])\ra\ldots\ra
C_1=\oplus_{|J|=2}A(J)\ra C_0=\oplus_{|J|=1}A(J).
\end{equation}
  
\begin{thm}\label{complexPi} One has $H_0(C_{\cdot}(P_i,A))\simeq A$
and $H_j(C_{\cdot}(P_i,A))=0$ for $j\neq 0$.
\end{thm}
 
The proof of this theorem will be given in \ref{proofPi}.
 
\subsection{Homological lemma}
Let $T$ be a finite set, and $(T^{1}_i)$ and $(T^2_i)$
two families of subsets of $T$ indexed by $i\in [1,n]$.
Assume that we have a family $(B_t)$ of objects of some
additive category indexed by $t\in T$. Then we can construct
a homology coefficient system on $\De_{n-1}$ by setting
$B(J)=\oplus_{t\in T(J)}B_t$ where
$$T(J)=\cap_{j\in J}T^1_j\cap\cap_{j\in\ov{J}}T^2_j,$$
$J\sub [1,n]$ is a subset, $\ov{J}$ is the complementary
subset. For $J\sub K$ we have the natural map
$B(K)\ra B(J)$ which is the following composition of the projection
and the embedding:
$$\oplus_{t\in T(K)}B_t\ra\oplus_{t\in T(J)\cap T(K)}B_t\ra
\oplus_{t\in T(J)}B_t.$$
Let  $D_{\cdot}=C_{\cdot}(\De_{n-1},B(\cdot))$  be  the  corresponding
chain complex.
  
\begin{lem}\label{homlem}  Assume that for every $t\in T$ the sets
$I^1(t)=\{i\ |\ t\in T^1_i\}$ and $I^2(t)=\{i\ |\ t\not\in T^2_i\}$
are different. Then $H_i(D_{\cdot})=0$ for $i>0$,
$H_0(D_{\cdot})\simeq B(\emptyset)=\oplus_{t\in T(\emptyset)}B_t$
where $T(\emptyset)=\cap_{j=1}^nT^2_j$.
\end{lem}
  
\Pf . Let $\wt{D}_{\cdot}$ be the complex
$$B([1,n])\ra\oplus_{|J|=n-1}B(J)\ra\ldots\ra
\oplus_{|J|=1}B(J)\ra B(\emptyset)$$
obtained from $D_{\cdot}$ by attaching one more term
$B(\emptyset)$. Then $\wt{\D}_{\cdot}$ is the direct sum
over $t\in T$ of the complexes  $\wt{D}(t)$ where
$\wt{D}_i(t)=\oplus_{J:t\in T(J),|J|=i}B_t$.
Note that $t\in T(J)$ if and only $I^2(t)\sub J\sub I^1(t)$.
Hence, the condition $I^1(t)\neq I^2(t)$ implies
that the complex $\wt{D}_i(t)$ is exact.
\ed
 
\subsection{Convexity}
 
\begin{lem}\label{W'Pi}
Let $W'=W_J\sub W$ for some $J\sub S$,
and let $y\in P_i$ be an element. Then $W'y\sub P_i$ if and only if
$ys_iy^{-1}\not\in W'$.
\end{lem}
 
\Pf. The second condition
is equivalent to the requirement that the cosets $W'y$ and
$W'ys_i$ are different. This is in turn equivalent
to the condition that the double coset
$W'y\lan s_i\ran$, where $\lan s_i\ran=\{1,s_i\}$,
has $2|W'|$ elements. Let $y_0$ be the shortest
element in this double coset (see \cite{B}, IV, Exercise 1.3).
Then every element  in  $W'y\lan  s_i\ran$  can  be  written
uniquely in  the  form  $w_1y_0w_2$  where  $w_1\in W'$, $w_2\in \lan
s_i\ran$ and $l(w_1y_0w_2)=l(w_1)+l(y_0)+l(w_2)$. Since there
are $2|W'|$ such expressions this gives a bijection between
$W'\times\lan s_i\ran$ and $W'y\lan s_i\ran$.
The condition $y\in P_i$ implies that
$y=w'y_0$ for some $w'\in W'$, i.~e. $W'y=W'y_0$. But for
every $w_1\in W'$ one has $l(w_1y_0s_i)=l(w_1y_0)+1$, hence
$W'y_0\sub P_i$.
\ed
    
\begin{lem}\label{appear}
For every $y\in P_i$ the set of $j\in [1,n]$ such that
$W^{(j)}y\sub P_i$ coincides with the set of $j$ such that
$s_j$ appears in a reduced decomposition of $ys_iy^{-1}$.
\end{lem}
 
\Pf . Applying Lemma \ref{W'Pi}
to  $W'=W^{(j)}$  we  obtain that $W^{(j)}y\sub P_i$ if and
only if $ys_iy^{-1}\not\in W^{(j)}$.  The latter  condition
means that $s_j$ appears in every (or some) reduced decomposition
of $ys_iy^{-1}$.
\ed
               
Consider the graph with vertices $W$ and
edges between $w$ and $sw$ for every $w\in W$, $s\in S$.
To give a path from $w_1$ to $w_2$ in this graph is
the same as giving a decomposition of $w_2w_1^{-1}$
into a product of simple reflections. The shortest paths,
{\it geodesics}, correspond to reduced decompositions.
Let us call a subset $P\sub W$ convex if every geodesic
between vertices in $P$ lies entirely in $P$.
 
\begin{prop}\label{convex}
The subset $P_i$ is convex.
\end{prop}
 
\Pf . First let us show that for every pair of elements
$y,w\in P_i$ there exists a geodesic from $y$ to $w$ which
lies entirely in $P_i$. Arguing by induction in $l(wy^{-1})$,
we see that it is sufficient to check that there exists
$s_k\in S$ such that either $s_ky\in    P_i$    and
$l(wy^{-1}s_k)<l(wy^{-1})$, or $ws_k\in P_i$ and
$l(s_kwy^{-1})<l(wy^{-1})$. Assume that  there  is  no  such
$s_k$. Note that by Lemma \ref{W'Pi} the condition $s_ky\in P_i$
is equivalent to $ys_iy^{-1}\neq s_k$. In particular,
this condition fails to be true for at most  one $k$.
Thus, for $w\neq y$ our assumption implies that
$ys_iy^{-1}=s_k$ for some $k$, $ws_iw^{-1}=s_l$ for some $l$
and $l(wy^{-1}s_j)>l(wy^{-1})$ for all $j\neq k$,
$l(s_jwy^{-1})>l(wy^{-1})$ for all $j\neq l$. Let us set
$x=wy^{-1}\in W$.  Then  $xs_kx^{-1}=s_l$  and  $x$  is  the
shortest element in $W^{(l)}xW^{(k)}$. Using Lemma \ref{sizig3}
one can easily see that this is impossible.
                                                               
Note that every two geodesics with common ends are connected
by a sequence of transformations of the following kind:
replace the segment
\begin{equation}\label{seg1}
y\ra s_j y\ra s_k s_j y\ra\ldots\ra w_0(s_j,s_k)y
\end{equation}
by the segment
\begin{equation}\label{seg2}
y\ra s_k y\ra s_j s_k y\ra\ldots\ra w_0(s_j,s_k)y
\end{equation}
where $w_0(s_i,s_j)$ is the longest element in the subgroup
generated by $s_i$ and $s_j$. Thus, it is sufficient to check
that if segment (\ref{seg1}) lies in $P_i$, then the corresponding
segment  (\ref{seg2}) does as well. Applying Lemma \ref{W'Pi} we see
that  for $y\in P_i$ the segment (\ref{seg1}) lies in $P_i$
if and only if
\begin{equation}\label{ysi}
ys_iy^{-1}\not\in\{s_j, s_ks_js_k, s_js_ks_js_ks_j,
\ldots\}.
\end{equation}
Now if $l(w_0(s_i,s_j))$ is odd, then the latter set
remains the same if we switch $s_j$ and $s_k$. Otherwise
(if $l(w_0(s_i,s_j))$ is even)
the path inverse to segment (\ref{seg2}) has the form
$$w_0(s_j,s_k)\ra s_j w_0(s_j,s_k)\ra\ldots\s_k y\ra y.$$
Hence, the conditions $w_0(s_j,s_k)\in P_i$ and (\ref{ysi})
imply that this segment lies in $P_i$.
\ed
 
\subsection{Proof of Theorem \ref{complexPi}}
\label{proofPi}
Let us denote by $Q_{i,j}$ the set of $y\in P_i$ such
that $W^{(j)}y\sub P_i$. Then for any $w\in P_i$ the complex
$j^*_wC_{\cdot}(P_i,A)$ can be written as 
\begin{equation}
\oplus_{x\in Q(w,[1,n])}\Phi_{wx^{-1}}A_x\ra\ldots
\ra\oplus_{|J|=2, x\in Q(w,J)}
\Phi_{wx^{-1}}A_x
\ra\oplus_{|J|=1, x\in Q(w,J)}
\Phi_{wx^{-1}}A_x
\end{equation}
where
$Q(w,J)=\cap_{j\in J} Q_{i,j}\cap\cap_{k\in \ov{J}}P_k^{-1}w$,
$A_x=j_x^*A$.
We can filter this complex by the length of $wx^{-1}$ as in  the
Proof of Theorem \ref{homology}. Now the associated graded factor
$\gr_Fj^*_wC_{\cdot}(P_i,A)$ is the complex
of Lemma \ref{homlem} with $T=P_i$,
$T^1_j=Q_{i,j}$ and $T^2_j=P_i\cap P_j^{-1}w$.
Since $\cap_{j=1}^n P_k^{-1}w=\{ w\}$  it  remains  to  show
that the conditions of Lemma \ref{homlem} are satisfied in
our situation. Indeed, assume that the set of $j$
such that $W^{(j)}y\sub P_i$ coincides with the set of $j$
such that $wy^{-1}\not\in P_j$. Let us denote this set by
$S_1\sub [1,n]$. Then for every $j\in S_1$ we have
$l(wy^{-1}s_j)<l(wy^{-1})$. Therefore,
$l(w_0wy^{-1}s_j)>l(w_0wy^{-1})$ for $j\in S_1$, where
$w_0\in W$ is the longest element. This means that
$w_0wy^{-1}$ is the shortest element in the coset
$w_0wy^{-1}W_{S_1}$. Hence, $l(w_0wy^{-1}w_1)=
l(w_0wy^{-1})+l(w_1)$ for every $w_1\in W_{S_1}$.
This can be rewritten as $l(wy^{-1}w_1)=l(wy^{-1})-l(w_1)$
for every $w_1\in W_{S_1}$. In other words, for every
$w_1\in W_{S_1}$ there exists a geodisic from $y$ to $w$
passing through $w_1y$. According to Lemma \ref{convex}
this implies that $w_1y\in P_i$ for every $w_1\in W_{S_1}$.
Now by Lemma \ref{appear}
we have $ys_iy^{-1}\in W_{S_1}$. Thus, taking
$w_1=ys_iy^{-1}$  we  obtain  that  $ys_i\in  P_i$  ---  
contradiction.
\ed

\subsection{Comparison with Beilinson---Drinfeld gluing} 
For every nonempty subset $J\sub [1,n]$
let us denote
$\CC_J=\oplus_\CC(\Phi_{W_{S-J}x})$
where   the   sum   is   taken   over
$x\in W_{S-J}x\backslash W$ such that
$W^{(j)}x\sub P_i$ for every $j\in J$. Then
$(\CC_J)$ is a family of abelian categories cofibered
over the category of nonempty subsets $J\sub [1,n]$
and their embeddings. Namely, if $J\sub K\sub [1,n]$,
then we have obvious restriction functors
$j^*_{J,K}:\CC_J\ra\CC_K$. Following Beilinson and Drinfeld \cite{Beil} 
we denote by $\CC_{\tot}$ the category of
cocartesian sections of $(\CC_J)$. An object of
$\CC_{\tot}$ is a collection of objects $A_J\in\CC_J$
and isomorphisms $\a_{J,K}:A_K\wt{\ra}j^*_{J,K}A_J$ for
$J\sub K$, such that
for $J\sub K\sub L$ one has
\begin{equation}\label{aJKL}
\a_{J,L}=j^*_{K,L}\a_{J,K}\circ \a_{K,L}.
\end{equation}
             
For every $J\sub [1,n]$ we have an obvious restriction functor
$j^*_J:\CC(\Phi_{P_i})\ra\CC_J$ and these functors fit
together into the functor $(j^*_{\cdot}):\CC(\Phi_{P_i})\ra\CC_{\tot}$.

\begin{thm} The functor
$(j^*_{\cdot}):\CC(\Phi_{P_i})\ra\CC_{\tot}$ is an
equivalence of categories.
\end{thm}

\Pf . For every $J\sub [1,n]$ we have the left adjoint
functor $j_{J,!}:\CC_J\ra\CC(\Phi_{P_i})$ to
$j^*_J$. Namely, $j_{J,!}$
is equal to $j_{W_Jx,P_i;!}$ on
$\CC(\Phi_{W_{S-J}x})\sub \CC_J$.
Now for every $(A_J,\a_{J,K})\in \CC_{\tot}$
the objects $j_{J,!}A_J$ form a homological coefficient system
on $\De_{n-1}$, so we can consider the corresponding
chain complex $\CC_{\cdot}(A_{\cdot})$. We claim that the
functor
$$(A_J)\mapsto H_0(\CC_{\cdot}(A_{\cdot}))$$
is quasi-inverse to $(j^*_{\cdot})$.
Notice that for $A\in\CC(\Phi_{P_i})$ the complex
$\CC_{\cdot}(j^*_{\cdot}A)$ coincides with the complex
(\ref{chaincomplex}), hence by
Theorem \ref{complexPi} it is a resolution of $A$.
It remains to show that for $(A_{\cdot})\in\CC_{\tot}$ and
for every $J\sub [1,n]$ there is a system of compatible
isomorphisms
$$H_0(j^*_J\CC_{\cdot}(A_{\cdot}))\simeq A_J.$$
It is sufficient to construct
canonical isomorphisms
$\a_k:H_0(j^*_{\{k\}}(\CC_{\cdot}(A_{\cdot})))\wt{\ra} A_{\{k\}}$
for every $k\in [1,n]$.
These morphisms are induced by the canonical
projections $\CC_0(A_{\cdot})\ra A_{\{k\}}$. To check that
$\a_k$ are isomorphisms it is sufficient to restrict everything
by some functor $j^*_w$ and to apply arguments from the proof of Theorem
\ref{complexPi}.
\ed

\subsection{Cohomological dimension}  
Following  \cite{Beil}  let  us  consider the
category $\secc_{-}=\secc_{-}(\CC_{\cdot})$ whose objects are
collections $A_J\in\CC_J$, $J\sub [1,n]$,
$\a_{J,K}:A_K\ra j_{J,K}^*A_J$ for $J\sub K$
satisfying (\ref{aJKL}). We  consider  $\CC_{\tot}$  as  the
subcategory in $\secc_{-}$ and denote by $\D^b_{\tot}$ the
full subcategory in the derived category
$\D^b(\secc_{-})$ consisting of complexes $C^{\cdot}$
with $H^i(C^{\cdot})\in\CC_{\tot}$.
The standard $t$-structure on $\D^b(\secc_{-})$ induces
a $t$-structure on $\D^b_{\tot}$ with core $\CC_{\tot}$.
As shown in \cite{Beil} for every $M,N\in\CC_{\tot}$
there    is    a    spectral    sequence    converging    to
$\Ext^{p+q}_{\D_{\tot}}(M,N)$ with
$E_1^{p,q}=\oplus_{|J|=p+1}\Ext^q_{\CC_J}(M_J,N_J)$.

\begin{thm}\label{anBe}
Assume that every object  of $\CC_w$ can be
covered by an object which is acyclic with respect to
all the functors $\Phi_{w'}$.
Then for every $M,N\in\CC_{\tot}$ the natural
map $\Ext^i_{\CC_{\tot}}(M,N)\ra\Ext^i_{\D_{\tot}}(M,N)$
is an isomorphism.
\end{thm}

\Pf . It is sufficient to prove that for every element
$e\in\Ext^i_{\D_{\tot}}(M,N)$ there exists a surjection
$M'\ra M$ in $\CC_{\tot}$ such that
the corresponding homomorphism $\Ext^i_{\D_{\tot}}(M,N)\ra
\Ext^i_{\D_{\tot}}(M',N)$ sends $e$ to zero.
Note that for every $w\in P_i$ the functor
$$j_{w,P_i;!}:\CC_w\ra\CC(\Phi_{P_i})\simeq\CC_{\tot}
\ra\secc_{-}$$
is left adjoint to the restriction functor
$j_w^*:\secc_{-}\ra\CC_w$. Furthermore, if an object
$P\in\CC_w$ is acyclic with respect to all the functors
$\Phi_{w'}$, then it is also $j_{w,P_i;!}$-acyclic.
Thus, for such an object $P$ we have
$\Ext^i_{\D_{\tot}}(j_{w,P_i;!}P,B)\simeq
\Ext^i_{\CC_w}(P,j^*_wB)$. Thus, we can start by
choosing surjections $P_w\ra j^*_wA$ such that
$P_w$ is acyclic with respect to all $\Phi_{w'}$
Then for every $w$ we can choose a surjection
$P'_w\ra P_w$ killing the image of $e$
in the space $\Ext^i_{\CC_w}(P_w,j^*_wB)$,
and take
$M'=\oplus_w j_{w,P_i;!}P'_w$ with the natural morphism
$M'\ra \oplus_w j_{w,P_i;!}j^*_wM\ra M$.
\ed

\begin{cor} Assume that the conditions of Theorem \ref{anBe}
are satisfied and in addition all categories $\CC(\Phi_{W_Jx})$,
where $J\sub [1,n]$ is a proper subset, have
finite cohomological dimension (this is true for the gluing
on the basic affine space if the rank of $G$ is equal to $2$).
Then the category $\CC(\Phi_{P_i})$ also
has finite cohomological dimension.
\end{cor}

\subsection{Gluing from two halves} 
Theorem \ref{complexPi} also implies that
the restriction functor
$j_{P_i}^*:\CC(\Phi)\ra\CC(\Phi_{P_i})$ has the left adjoint.
Namely,
for every $A\in\CC(\Phi_{P_i})$ we have
$$A=\coker
(\oplus_{|J|=2}A(J)\ra \oplus_j A(j)).$$
Hence, for any $B\in\CC(\Phi)$ we have
$$\Hom(A,j^*_{P_i}B)\simeq
\ker(\oplus_j\Hom(A(j),j^*_{P_i}B)\ra\oplus_{|J|=2}
\Hom(A(J),j^*_{P_i}B)).$$
Now by adjunction we have
$$\Hom(A(J),j^*_{P_i}B)\simeq\oplus
\Hom(j^*_{W_{S-J}x,P_i}A,j^*_{W_{S-J}x}B)
\simeq\oplus
\Hom(j_{W_{S-J}x,!}j^*_{W_{S-J}x,P_i}A,B)
$$
where the sum is taken over $x\in W_{S-J}\backslash W$
such that $W^{(j)}x\sub P_i$ for every $j\in J$.
It follows that we have an isomorphism
\begin{align*}
\Hom(A,j^*_{P_i}B)\simeq
&\ker(\Hom(\oplus_{j,x} j_{W^{(j)}x,!}j^*_{W^{(j)}x,P_i}A,B)\ra
\Hom(\oplus_{|J|=2,x}
j_{W_{S-J}x,!}j^*_{W_{S-J}x,P_i}A,B))\\
&\simeq \Hom(j_{P_i,!}A,B)
\end{align*}
where
$$j_{P_i,!}A=\coker(
j_{W_{S-J}x,!}j^*_{W_{S-J}x,P_i}A\ra
\oplus_{j,x} j_{W^{(j)}x,!}j^*_{W^{(j)}x,P_i}A).$$
                                                
Theorem \ref{complexPi} and the above construction work
almost literally for the subset $P_is_i\sub W$ instead
of $P_i$. Now we can consider the functors
$j^*_{P_is_i}j_{P_i,!}$ and $j^*_{P_i}j_{P_is_i,!}$
as gluing data on the pair of categories
$\CC(\Phi_{P_i})$ and $\CC(\Phi_{P_is_i})$. Theorem \ref{braverman}
implies that the corresponding glued category is equivalent
to $\CC(\Phi)$.

\subsection{Supports of simple objects}
Now we are going to apply the explicit construction of the adjoint functors
corresponding to cosets of parabolic subgroups in $W$ to the study of
simple objects in the glued category.  

\begin{prop}\label{supp1} 
Let $P\sub W$ be a subset, 
$S$ a simple object of $\CC(\Phi_P)$, 
and $x\in P$ an element such that $s_ix\in P$ for some simple 
reflection $s_i$. Assume that $S_x=0$ and $S_{s_ix}\neq0$. 
Then $S_w=0$ for every $w\in P\cap P_ix$. 
\end{prop} 

\Pf . Let $W_i\sub W$ be the subgroup generated by $s_i$. Our 
assumptions imply that  $j^*_{W_ix,P}S\neq 0$, hence the
adjunction morphism 
$$A=j_{W_ix,P;!}j^*_{W_ix,P}S\ra S$$
is surjective. Now since $S_x=0$, the explicit construction of
$j_{W_ix,P,!}$ tells us that the $A_w=0$ if $w$ is closer to $x$
than to $s_ix$. The latter condition is equivalent to $wx^{-1}\in P_i$.
By surjectivity for such $w$ we also have $S_w=0$.
\ed

For every object $A=(A_w,w\in W)$ of $\CC(\Phi)$ we 
denote by $\Supp(A)$ the set of $w\in W$ such that $A_w\neq 0$. 
The above proposition gives serious restrictions on a subset $\Supp(S)$
for a simple object $S\in\CC(\Phi)$.

\begin{thm}\label{supp2} Let $S$ be a simple object of $\CC(\Phi)$. 
Then either $\Supp(S)=W$ or $\Supp(S)$ is an intersection of 
subsets of the form $P_ix$. In particular, if $\Supp(S)\neq W$, 
then $\Supp(S)$ is convex. 
\end{thm}

\Pf . Note that an intersection of 
subsets of the form $P_ix$ is convex by Lemma \ref{convex}.
Thus, it suffices to prove that $W-\Supp(S)$ is a union of
subsets of the form $P_ix$. Let $w\not\in\Supp(S)$. Choose a 
geodesic path 
$$w\ra s_{i_1}w\ra\ldots s_{i_k}\ldots s_{i_1}w=x$$ 
of minimal length from $w$ to an element $x$ of $\Supp(S)$. Then we can 
apply Proposition \ref{supp1} to $x'=s_{i_k}x$ to conclude that
$P_{i_k}x'\sub W-\Supp(S)$. On the other hand,
clearly $w\in P_{i_k}x'$.
\ed

\section{Extensions in the glued category}
                   
\subsection{Adapted objects and finiteness of dimensions}
The proof of the following theorem is due to L.~Positselski.
       
\begin{thm}\label{adapted}
Let $k:V\hra X$ be an affine open subset
such that the projection $G\ra X$ splits over $V$.
Then the functors $F_{w,!}k_!$ and $F_{w,*}k_*$ are $t$-exact.
\end{thm}
       
\Pf . Since $k_!$ is $t$-exact and  $F_{w,!}$  is  $t$-exact
from the right, it is sufficient to prove that
$F_{w,!}k_!$ is $t$-exact from the left. Now by definition
the functor $F_{w,!}$ is given by the kernel
$\ov{K(w)}$ on $X\times X$, hence the functor
$F_{w,!}k_!$ is given by the kernel $\ov{K(w)}|_{V\times X}$
on $V\times X$. Since the projection $p_2:V\times X\ra X$
is affine, the functor $p_{2!}$ is left $t$-exact; hence
it is sufficient to prove that for any $A\in\Perv(V)$
the object $p_1^*A\otimes\ov{K(w)}|_{V\times X}$ is
a perverse sheaf on $V\times X$. Let $s:V\ra G$ be
a splitting of the projection $G\ra X$ over $V$.
Consider the isomorphism
$$\nu:V\times X\ra V\times X:(v,x)\mapsto (v,s(v)x).$$
Then it is sufficient to check that
$\nu^*(p_1^*A\otimes\ov{K(w)})\simeq p_1^*A\otimes
\nu^*\ov{K(w)}$ is perverse. The sheaf $\nu^*\ov{K(w)}$
is the Goresky---MacPherson extension of
$\nu^*K(w)$ on $\nu^{-1}(X(w))$. By definition
$$\nu^{-1}(X(w))=\{(v,x)\in V\times X|\ x\in X_w\}$$
where $X_w=Un_wT_wU/U\sub X$ is a locally closed subvariety
of $X$. Note that
the projection $\pr_w:X(w)\ra T_w$ factors as the composition
of projections $p_w:X(w)\ra X_w$ and $q_w:X_w\ra T_w$
where $p_w$ is smooth of relative dimension $\dim X$,
$q_w$ is smooth of relative dimension $l(w)$.
Also we have $K(w)=\pr_w^*L_w[l(w)]$ where $L_w$ is a perverse
sheaf on $T_w$. Hence, $K(w)=p_w^*(q_w^*L_w[l(w)])$ where
$q_w^*L_w[l(w)]$ is a perverse sheaf on $X_w$.
Now  we  have
$p_w\circ\nu|_{\nu^{-1}(X(w)}=p_2|_{\nu^{-1}(X(w))}$
where $p_2:V\times X\ra X$ is the projection, and hence
$$\nu^*K(w)\simeq p_2^*(q_w^*L_w[l(w)]).$$
It follows that $\nu^*\ov{K(w)}\simeq p_2^*M_w$
where $M_w$ is the Goresky---MacPherson extension of
$q_w^*L_w[l(w)]$ to $X$. Thus,
$$p_1^*A\otimes\nu^*\ov{K(w)}\simeq A\boxtimes M_w,$$
and the latter sheaf is perverse by \cite{BBD}, 4.2.8.
     
The exactness of functor $F_{w,*}k_*$ follows
from the isomorphism
$$F_{w,\psi,*}=D\circ F_{w,\psi^{-1},!}\circ D$$
where $D$ is the Verdier duality (see \cite{KL}, 2.6.2(i)).
\ed
       
\begin{thm}\label{findim} For any $A,B\in\AA$ the spaces
$\Ext^i_{\AA}(A,B)$ are finite-dimensional.
\end{thm}
          
\Pf . It is sufficient to prove that for every
$A\in\AA$ there exists a surjection $A'\ra A$
in $\AA$ such that $\Ext^i_{\AA}(A,B)$
is finite-dimensional for every $B\in\AA$.
We start with the surjection
$$\oplus_w j_{w,!}A_w\ra A$$
where $A_w=j_w^*A$.
Now we choose a finite covering $(U_k)$ of $X$
by affine open subsets such that the projection
$G\ra X$ splits over every $U_k$. Let
$j_k:U_k\hra X$ be the corresponding open embeddings.
Now we replace every $A_w$ by $j_{k,!}j_k^*A_w$
to get the surjection
$$A'=\oplus_{i,w} j_{w,!}j_{k,!}j_k^*A_w\ra A.$$
It remains to prove that
$\Ext^i_{\AA}(j_{w,!}j_{k,!}C,B)$ is finite-dimensional
for every $C\in\Perv(U_k)$. According to Theorem \ref{adapted}
the functor $j_{w,!}j_{k,!}:\Perv(U_k)\ra\AA$ is exact.
Since it is left adjoint to $j_k^*j_w^*:\AA\ra\Perv(U_k)$
we get the isomorphism
$$\Ext^i_{\AA}(j_{w,!}j_{k,!}C,B)\simeq
\Ext^i_{\Perv(U_k)}(C,j_k^*j_w^*B),$$
where the latter space is finite-dimensional, as follows from
Beilinson's theorem \cite{Be}.
\ed

\subsection{Vanishing of $Ext^1$}
Let $(W,S)$ be a finite Coxeter system,
$\Phi$ be a $W$-gluing data.  
        
\begin{thm}\label{Ext1van}
Assume that for every $w\in W$ and $s\in S$  the
following condition holds: for every object $A\in\CC_w$ such
that the canonical morphism $\Phi_s^2(A)\ra A$ is zero, one
has $\Phi_s(A)=0$. Let $S$ and $S'$ be simple objects in $\CC(\Phi)$
such that $\Supp(S)\cap\Supp(S')=\emptyset$.
Then $\Ext^1_{\CC(\Phi)}(S,S')=0$.
\end{thm}

\Pf . By assumption we have a partition of $W$ into two subsets 
$I$ and $I'$ such that $S_w=0$ for $w\in I'$ and $S'_{w}=0$ 
for $w\in I$. Then we have $E_w=S_w$ for $w\in I$, and $E_w=S'_w$ 
for $w\in I'$. Now we claim that the structural morphisms  
$\Phi_sE_w\ra  E_{sw}$ vanish unless $w$ and $sw$ belong to the 
same subset of this partition. By the definition of $W$-gluing 
this would imply that $E\simeq S\oplus S'$. So assume for 
example that $w\in I$ and $sw\in I'$. Then the morphism 
$\Phi_s^2  S_w\ra  S_w$  is  zero, and hence $\Phi_s E_w=\Phi_s S_w=0$. 
Similarly, if $w\in I'$ and $sw\in I$, then $\Phi_s E_w=\Phi_s S'_w=0$.
\ed

\begin{cor}\label{Ext1inj} 
Let $\Phi$ be a $W$-gluing data  of finite type.
Under the assumptions of Theorem \ref{Ext1van},
for every pair of simple objects
$S=(S_w)$ and $S'=(S'_w)$ of $\CC(\Phi)$
the natural map
$$\Ext^1_{\CC(\Phi)}(S,S')\ra
\oplus_{w\in W}\Ext^1_{\CC_w}(S_w,S'_w)$$
is injective.
\end{cor}
        
\Pf . If the supports of $S$ and $S'$ do not intersect then we are done
by Theorem \ref{Ext1van}. So we can assume that
there  exists  $w\in  W$  such that $S_w$ and $S'_w$ are both
non-zero. Then according
to Lemmas \ref{ressim} and \ref{gormac}
we have $S\simeq j_{w,!*}S_w$ and
$S'\simeq j_{w,!*}S'_w$. Now assume that we have an extension
\begin{equation}\label{extSS'}
0\ra S'\ra E\ra S\ra 0
\end{equation}
in $\CC(\Phi)$ that induces a trivial extension
\begin{equation}\label{extw}
0\ra S'_w\ra E_w\ra S_w\ra 0.
\end{equation}
We claim that the adjunction morphism
$j_{w,!}E_w\ra E$ is surjective. Indeed, this
follows immediately from the 
commutative diagram with exact rows:
\begin{equation}
\begin{array}{cccccc}
j_{w,!}S'_w & \lrar{} & j_{w,!}E_w & \lrar{} &
j_{w,!}S_w &\lrar{} 0\\
\ldar{}& &\ldar{}& &\ldar{}\\
S'& \lrar{} & E & \lrar{} &
S &\lrar{} 0\\
\end{array}
\end{equation}
Similarly, one proves that the adjunction
morphism $E\ra j_{w,*}E_w$ is injective. Thus,
$E\simeq j_{w,!*}E_w$, and a splitting of the
extension (\ref{extw})
induces a splitting of (\ref{extSS'}).     
\ed
        
\begin{rems} 

\noindent
1. The conditions of Theorem \ref{Ext1van} are satisfied for
the gluing on the basic affine space as one can see combining
Corollary \ref{zer} and Lemma \ref{mainl}.
       
\noindent 2. There is an analogue of Theorem \ref{Ext1van}
for arbitrary gluing data of finite type. In this case one
should impose the condition that $\Phi_{ij}(\CC_{ij})=0$
(in the notation of section \ref{groth}). For example,
this condition is satisfied for the usual geometric gluing
data associated with an open covering (see \cite{KL}).
\end{rems}

Theorem \ref{Ext1van} can be strengthened using the adjoint
functors. 

\begin{thm} Assume that we have a gluing data $\Phi$ as in Theorem
\ref{Ext1van}.
Let $P\sub W$ be a subset such that there exists left and right adjoint
functors $j_{P,!},j_{P,*}:\CC(\Phi_P)\ra\CC(\Phi)$ to the natural 
restriction functor
$j_P^*:\CC(\Phi)\ra\CC(\Phi_P)$.
Let $S$ and $S'$ be simple objects in $\CC(\Phi)$ such that
$\Supp(S)\cap P\neq \emptyset$, $\Supp(S')\cap P\neq \emptyset$ but 
$\Supp(S)\cap\Supp(S')\cap P=\emptyset$. Then
$\Ext^1_{\CC(\Phi)}(S,S')=0$.
\end{thm}

\Pf . The proof of Theorem \ref{Ext1van} shows that there are no 
non-trivial extensions between $j_P^*S$ and $j_P^*S'$.
Now given an extension between $S$ and $S'$ the same argument as in
Corollary \ref{Ext1inj} shows that it splits.
\ed

For example, let $W=S_3$, $\Phi$ be a $W$-gluing data of finite type.
Then the conditions of this theorem are satisfied for any proper connected
subgraph $P\sub S_3$. Together with Theorem \ref{supp2} this implies
that in this case $Ext^1_{\CC(\Phi)}(S,S')=0$ unless either
$\Supp(S)\sub\Supp(S')$ or $\Supp(S')\sub\Supp(S)$.

\section{Mixed glued category}
        
\subsection{Definition}
We are going to define a mixed analogue of the  gluing  on  the
basic affine space. Namely, we consider the situation when
the field $k$ is the algebraic closure of $\F_p$ ($p>2$),
and we fix the finite subfield $\F_q\sub k$ such that
$G$, $T$, and $B$ are defined over $\F_q$.
Following \cite{BBD} we denote by subscript
$0$ objects defined over this finite field.
Thus $X=G/U$ is  obtained  by extension of scalars from
$X_0$, the corresponding scheme over $\F_q$.
Let $\Fr:X\ra X$ be the corresponding geometric Frobenius.
Since the Fourier transform commutes with $\Fr$, there is a
well-defined  functor  $\Fr^*$   on   the   glued   category
$\AA=\AA_{\psi}$.
Let $\AA^{\Fr}$ be the category of objects $A\in\AA$ together
with isomomorphisms $\a:A\wt{\ra}\Fr^*A$.
Let $\AA_m$ be the subcategory in $\AA^{\Fr}$ consisting of
$A=(A_w)\in\AA$  such  that  all  $A_w$  are  mixed
perverse sheaves on $X_0$, with the canonical morphism $\a$.
Then $\AA_m$ is obtained by gluing from $|W|$ copies
of the category $\Perv_m(X_0)$ of mixed perverse sheaves on
$X_0$ (since the relevant functors between perverse
sheaves preserve mixedness). We will call an object of $\AA_m$
{\it pure} if all its components are pure perverse sheaves of the
same weight. 
        
\subsection{Action of Frobenius on extensions}
Let $\CC$ be a $\ov{\Q_l}$-linear abelian category, and
$\Fr^*:\CC\ra\CC$ a $\ov{\Q_l}$-linear exact functor. Let $\CC^{\Fr}$
be the category of pairs $(A,\a)$ where $A\in\CC$,
$\a:A\wt{\ra}\Fr^*A$ is an isomorphism. Assume that
$\CC_0$ is a full subcategory of $\CC^{\Fr}$.
Let us denote  by  letters  with  subscript  $0$
objects in $\CC_0$ and omit the subscript when considering
the corresponding objects of $\CC$. For every pair of objects
$A_0$ and $B_0$ of $\CC_0$ there is a  natural  automorphism
$\Fr$ on $\Ext^*_{\CC}(A,B)$ induced by $\Fr^*$. Let us assume that
for every $i$ the $\ov{\Q_l}$-space $\Ext^i_{\CC}(A,B)$
is finite-dimensional and the eigenvalues of the above endomorphism
are  $l$-adic  units.  Finally,  notice  that  if  $V$  is a
finite-dimensional $\ov{\Q_l}$-vector   space   with   an
automorphism $\phi$, then we can define a functor of twist with
$(V,\phi)$ on $\CC^{\Fr}$ by sending $(A,\a)$ to the object
$A'=A\otimes_{\ov{\Q_l}} V$ with the isomorphism
$\a'=\a\otimes\phi:A'\wt{\ra}\Fr^*A'$. Our  last  assumption
is that the subcategory $\CC_0\sub\CC^{\Fr}$ is stable under
twists by finite-dimensional $\ov{\Q_l}$-vector spaces with
continuous $\hat{\Z}$-action.
                                     
\begin{thm}\label{Frinv}
With  the  above  assumptions the following two
conditions are equivalent:
\begin{enumerate}
\item canonical maps
$\Ext^i_{\CC_0}(A_0,B_0)\ra\Ext^i_{\CC}(A,B)^{\Fr}$
are surjective for all $A_0,B_0\in\CC_0$, $i\ge 0$,
\item for every $A_0,B_0\in\CC_0$ and every element
$e\in\Ext^i_{\CC}(A,B)$ there exists a morphism
$f:B_0\ra B'_0$ such that the image of $e$ in
$\Ext^i_{\CC}(A,B')$ is zero.
\end{enumerate}
\end{thm}
        
\Pf. (1)$\implies$(2). Obviously it is sufficient to check the
required  condition  for  elements  $e$  belonging to one of the
generalized eigenspaces of $\Fr$. Twisting $A_0$ or $B_0$ with a
one-dimensional $\hat{\Z}$-representation we can assume that
$(\Fr-1)^n\cdot  e=0$ for some $n\ge 1$.
Consider the element $e_1=(\Fr-1)^{n-1}\cdot e$. Then
$e_1$ is invariant under $\Fr$; hence by assumption it is induced by some
extension in $\CC_0$ between $A_0$ and $B_0$. In particular,
there exists a morphism $f:B_0\ra B'_0$ such that the
image  of  $e_1$  in $\Ext^i(A,B')$ is zero. Let $e'$ be the
image of $e$ in $\Ext^i(A,B')$. Then $(\Fr-1)^{n-1}e'=0$,  so
we can apply induction to finish the proof.
        
\noindent
(2)$\implies$(1).  We use induction in $i$.
For  $i=0$  the  required  surjectivity  follows  from   the
assumption that $\CC_0$ is a full subcategory of $\CC^{\Fr}$.
Let $i\ge 1$, $e\in\Ext^i(A,B)$  be  an  element
invariant under $\Fr$. Then by assumption there exists
an embedding $B_0\hra B'_0$ killing $e$. In other words, $e$
is the image of some element $e'\in\Ext^{i-1}(A,B'/B)$
under the boundary map
$\de:\Ext^{i-1}(A,B'/B)\ra\Ext^i(A,B)$. Notice that the element
$\Fr(e')-e'$ lies in the kernel of $\de$; hence it comes
from some element $d\in\Ext^{i-1}(A,B')$. Applying our
assumption to this element we find an embedding
$B'_0\hra B''_0$ such that $d$ is killed in
$\Ext^{i-1}(A,B'')$. Let $e''$ be the image of $e'$ in
$\Ext^{i-1}(A,B''/B)$. Then $e''$ is invariant under $\Fr$
and maps to $e$ under the boundary map
$\de':\Ext^{i-1}(A,B''/B)\ra\Ext^i(A,B)$. Now the proof is
finished by applying the induction hypothesis to $e''$.
\ed
        
\subsection{Weights of Ext-groups}
In this section we will temporarily use the notation
$\Phi_w=\sideset{^p}{^0}{H}F_{w,!}$ for Kazhdan---Laumon
gluing functors (reserving symbols $F_{\cdot}$ for 
the filtration on mixed objects).

\begin{prop}\label{pure}
Every simple object of $\AA_m$ is pure.
Every object $A$ of $\AA_m$ has a canonical increasing
filtration such that $\gr_n(A)$ is pure of weight $n$.
Every morphism in $\AA_m$ is strictly compatible with
this filtration.
\end{prop}
        
\Pf . Since $\AA_m$ is obtained by gluing,
from Lemmas \ref{ressim} and \ref{gormac} we conclude that
every simple object of $\AA_m$ has form
$j_{w,!*}(A_w)$ for some simple mixed perverse sheaf $A_w$ on
$X_0$. Recall that by  \cite{BBD}, 5.3.4  $A_w$ is pure
of some weight $n$.
Now since the symplectic Fourier transform preserves weights
we obtain that $j_{w,!}(A_w)$ has weights $\le n$, while
$j_{w,*}(A_w)$ has weights $\ge n$. Hence, $j_{w,!*}$ is pure
of weight $n$.
        
Let $A=(A_w)$ be an object of $\AA_m$, and let $F_{\cdot}(A_w)$ be
the canonical filtration on $A_w$ such that
$F_n(A_w)/F_{n-1}(A_w)$ is pure of weight $n$. Then
the weights of $\Phi_s(F_n(A_w))$ are $\le n$; hence
the structural morphisms $\Phi_s(A_w)\ra A_{sw}$ induce
the morphisms $\Phi_s(F_n(A_w))\ra F_n(A_{sw})$. Thus,
$(F_n(A_w))$ is a subobject of $A$ for every $n$, and
the $(F_{\cdot}(A_w))$ is the filtration with required
properties.
\ed
        
\begin{rem} Another proof of existence of canonical
filtrations in $\AA_m$ can be obtained using Corollary
\ref{Ext1inj} and \cite{BBD}, 5.1.15 and 5.3.6.
\end{rem}
        
\begin{thm}\label{bound}
Let $A_0$ and $B_0$ be pure objects in $\AA_m$
of weights $a$ and $b$, respectively. Then
weights of Frobenius in $\Ext^i_{\AA}(A,B)$ are
$\ge i+b-a$.
\end{thm}
        
\Pf . Without loss of generality we can assume that $A_0$ and
$B_0$ are simple.
For $i=0$ the assertion is clear, since
$\Hom_{\AA}(A,B)\hra\oplus_{w\in W}\Hom(j_w^*A,j_w^*B)$
and the weights of Frobenius on the latter space are
$\ge b-a$.
Similarly, for $i=1$ the statement follows from Corollary \ref{Ext1inj}
and \cite{BBD}, 5.1.15.
So let us assume that $i>1$ and
that the assertion is true for $i-1$.
Using twist one can easily see that
it is sufficient to prove that for $i>a-b$ one
has $\Ext^i_{\AA}(A,B)^{\Fr}=0$.
        
We claim that equivalent conditions
of Theorem \ref{Frinv} are satisfied for the pair
of categories $\AA_m\sub\AA$. Indeed,
the action of Frobenius on $\Ext^i$-spaces extends
to continuous $\hat{\Z}$-action as follows from
the proof of Theorem \ref{findim}. Next,
the equivalent conditions of \ref{Frinv} are satisfied for the pair
$\Perv_m(X_0)\sub\Perv(X)$ since by Beilinson's theorem
the $\Ext$-groups in these categories can be computed
in derived categories of construcible sheaves, and
then the condition (1) of \ref{Frinv}
follows from ``local to global" spectral sequence
(see e.g., \cite{BBD} 5.1). Now we can check the condition (2) of
\ref{Frinv} for $\AA_m\sub\AA$ as follows. Let
$e\in\Ext^i_{\AA}(A,B)$ be an element. For every $w\in W$
consider the embedding $j^*_wB_0\hra C_{w,0}$ such that
$C_{w,0}$ is $j_{w,*}$-acyclic (such an embedding exists
by Theorem \ref{adapted}). Then
$\Ext^i(A,j_{w,*}C_w)\simeq\Ext^i(j^*_wA,C_w)$; hence we
can find an embedding $C_{w,0}\ra C'_{w,0}$ such that
the image of $e$ is killed in $\Ext^i(j^*_wA,C'_w)$. Then
$$B_0\hra B'_0=\oplus_w j_{w,*}C'_{w,0}$$
is the required embedding. Thus the morphism
$\Ext^i_{\AA_m}(A,B)\ra\Ext^i_{\AA}(A,B)^{\Fr}$
is surjective.
        
Now given an Ioneda extension class
in $\AA_m$
\begin{equation}\label{Ion}
0\ra B_0\ra E^1_0\ra E^2_0\ra\ldots E^i_0\ra A_0\ra 0
\end{equation}
we can replace $E^1_0$ and $E^2_0$ by their quotients
$E^1_0/F_{b-1}(E^1_0)$ and $E^2_0/F_{b-1}(E^1_0)$ (where
$F_{\cdot}(\cdot)$ is the weight filtration) to get an equivalent
extension class (\ref{Ion}) such that weights of $E_0^1$
are  $\ge  b$. Let $C_0=E^1_0/B_0$. Then
the class $e$ of (\ref{Ion}) is an image under the boundary
map
$$\Ext^{i-1}(A_0,C_0)\ra\Ext^i(A_0,B_0)$$
of some class $e_1\in\Ext^{i-1}(A_0,C_0)$.
Since weights of $C_0$ are $\ge b$ we have an exact sequence
$$\Ext^{i-1}(A,F_b(C))\ra\Ext^{i-1}(A,C)\ra
\Ext^{i-1}(A,C/F_b(C)).$$
For every extension class $e$ in $\AA_m$ let us denote
by $c(e)$ the corresponding extension class in $\AA$.
By the induction hypothesis the image of $c(e_1)$ in
$\Ext^{i-1}(A,C/F_b(C))$ is zero, hence, $c(e_1)$ comes
from some class $e'_1\in\Ext^{i-1}(A,F_b(C))$.
Now the class $c(e)$ is the image of $c(e_1)$
under the boundary map
$$\Ext^{i-1}(A,C)\ra\Ext^i(A,B);$$
hence it lies in the image of the boundary map
$$\Ext^{i-1}(A,F_b(C))\ra \Ext^i(A,B)$$
corresponding to the exact sequence
$$0\ra B_0\ra F_b(E^1_0)\ra F_b(C_0)\ra 0.$$
But the image of this exact sequence in $\AA$ splits,
hence $c(e)=0$.
\ed

\begin{cor}\label{findimwt}
For every $A_0,B_0\in\AA_m$ and every
$n\in\Z$ let $\Ext^*_{\AA}(A,B)_n$ denotes the weight-$n$
component of $\Ext^*_{\AA}(A,B)$. Then all the spaces
$\Ext^*_{\AA}(A,B)_n$ are finite-dimensional and
$\Ext^*_{\AA}(A,B)_n=0$ for $n<<0$.
\end{cor}

\section{Induction for representations of braid groups}
\label{induction}
   
\subsection{Formulation of the theorem}
Let  $(W,  S)$  be  a  finite  Coxeter  group,  $B$  the
corresponding braid group, and $B^+\sub  B$  the  monoid  of
positive braids.
We fix a subset $J\sub S$ of simple reflections. Let
$W_J\sub  W$ be the subgroup generated by simple reflections
in $J$. Then $(W_J, J)$ is a Coxeter group and we denote
by $B_J$ the corresponding braid group.
Let $\Mod_J-B$ be the category of representations of $B$
of the form $\oplus_{x\in W/W_J} V_x$ such that the action of
$b\in B$ sends $V_x$  to  $V_{\ov{b}x}$  and   the
following condition is satisfied: for  every  $s\in  S$  and
every $x\in W/W_J$
such that $sx\neq x$ one has $s^2|_{V_x}=\id_{V_x}$.
Morphisms in
$\Mod_J-B$   are   morphisms   of   $B$-modules   preserving
direct sum decompositions. Let $\Mod_J-B^+$ be the similar
category of $B^+$-representations.
Let $x_0\in W/W_J$ be the coset
containing the identity.
        
\begin{thm}\label{braidmain}
The functor $\oplus_{x\in W/W_J} V_x\mapsto V_{x_0}$
is an equivalence of $\Mod_J-B$ with the category $\Mod-B_J$
of $B_J$-representations. Similarly the category
$\Mod_J-B^+$ is equivalent to $\Mod-B_J^+$.
\end{thm}

\subsection{Arrangements of hyperplanes}
Let us realize $W$ as the group generated by reflections
in a Euclidean vector space $V$ over $\R$.
Let $\HH$ be the corresponding arrangement of hyperplanes in
$V$, $X$ be the complement in $V_{\C}$
to all the complex hyperplanes $H_{\C}$, $H\in \HH$. For
a subset $J\sub S$ of simple roots we denote by $X_J$
the similar space associated with $W_J$.
By a theorem of P.~Deligne \cite{D1} the spaces
$X/W$ and $X_J/W_J$ are $K(\pi,1)$ with fundamental groups
$B$ and $B_J$, respectively. Let $X'_J$ be the complement
in $V_{\C}$ to all the complex hyperplanes $H_{\C}$ such
that the corresponding reflection $r_H$ belongs to $W_J$.
Then $X'_J$ is a $W_J$-invariant open subset
containing $X$. Furthermore, there is a $W_J$-equivariant
retraction $X'_J\ra X_J$ which induces a homotopic equivalence
of $X'_J/W_J$ with $X_J/W_J$. Let us consider the cartesian
power $Y_J=(X'_J/W_J)^{W/W_J}$ with the natural action of
$W$ by permutations and let us denote by
$X\times_W Y_J$ the quotient of $X\times Y_J$ by the diagonal
action of $W$. Let
$f:X\times_W Y_J\ra X/W$ be the natural projection.
Then $f$ is a fibration with fiber $Y_J$. Let $x_0\in X$
be a fixed point. Then $y_0=(x_0,\ldots,x_0)\in Y_J$
is a point fixed by $W$ and hence the map
$$X\ra X\times Y_J:x\mapsto (x,y_0)$$
descends to a section
$\si_0:X/W\ra X\times_W Y_J$ of $f$. Let $\ov{x_0}$ be the
image of $x_0$ in $X/W$. Using $\si_0$ we obtain a canonical
identification
$$\pi_1(X\times_W Y_J,\si_0(\ov{x_0}))\simeq
\pi_1(Y_J,y_0)\rtimes\pi_1(X/W,\ov{x_0}).$$
Note that $\pi_1(X/W,\ov{x_0})\simeq B$ acts
on $\pi_1(Y_J,y_0)\simeq B_J^{W/W_J}$ via the
action of $W$ by permutations. In particular,
we have the canonical surjection
$$\pi_1(X\times_W Y_J,\si_0(\ov{x_0}))\ra
B_J^{W/W_J}\rtimes W.$$
Now let us consider a $W$-equivariant map
$$\wt{\si}:X\ra X\times Y_J: x\mapsto (x,\pi_J(\ov{w^{-1}x})_{w\in W/W_J})$$
where $\pi_J:X'_J\ra X'_J/W_J$ is the natural projection.
Then $\wt{\si}$ induces a section
$\si:X/W\ra X\times_W Y_J$ of $f$, and hence a homomorphism
$$\si_*:B\ra\pi_1(X\times_W Y_J,\si(\ov{x_0})).$$
To  identify  the   latter   group   with
$\pi_1(X\times_W Y_J,\si_0(\ov{x_0}))$
we   have to construct  a  path  from
$\si_0(x_0)$ to $\si(x)$. In other words for every  $wW_J\in
W/W_J$ we have to construct a path from $\pi_J(x_0)$ to
$\pi_J(w^{-1}x_0)$  in  $X'_J/W_J$.  To do this we
take the canonical representative $w\in W$ of every $W_J$-coset
and consider the corresponding path between $x_0$ and
$w^{-1}x_0$  in  $X$  (there  is  a  canonical  choice up to
homotopy,   corresponding   to   the   section    $\tau:W\ra
B=\pi_1(X/W,\ov{x_0})$),
then project it to $X'_J/W_J$.
This gives the required identification  of  the  fundamental
groups, hence we get a homomorphism
$$f_J:B\ra B_J^{W/W_J}\rtimes W.$$
It is easy to check that for every $s\in S$ one has
$$f_J(s)=((b_w)_{w\in W/W_J},\ov{s})$$
where $w$ runs over canonical representatives of $W/W_J$,
$b_w=\tau(w^{-1}sw)$ if $w^{-1}sw\in W_J$, and $b_w=1$
otherwise. Note that since $w$ is a canonical representative
the condition $w^{-1}sw\in W_J$ implies that
$l(sw)=l(w)+l(w^{-1}sw)$; hence $w^{-1}sw$ is in fact a simple
reflection.
      
\subsection{Proof of Theorem \ref{braidmain}}
The inverse functor
$\Mod-B_J\ra\Mod_J-B$ is constructed as follows.
Let  $V_0$  be  a  representation  of $B_J$. Then there is a
natural action of $B_J^{W/W_J}\rtimes W$ on
$V=\oplus_{x\in W/W_J}V_0$ such that $B_J^{W/W_J}$
acts component-wise while $W$ acts by permutations of
components in the direct sum. Using the homomorphism
$f_J$ we obtain the action of $B$ on $V$. We claim
that $V$ belongs to the subcategory
$\Mod_J-B$. Indeed, by the definition of $f_J$,
if $sx\neq x$ for $x\in W/W_J$ and $s\in S$; then
$s$ acts by the permutation of factors on
$V_x\oplus V_{sx}\sub V$.
      
It remains to observe that for any $V\in\Mod_J-B$ we can
identify all the components $V_x$ with $V_{x_0}$ using
the action of the canonical representative for $x$.
This gives an isomorphism of $V$ with the $B$-representation
associated with $V_{x_0}$ as above.

\subsection{Commutator subgroup of the pure braid group}
When $J$ consists of one element we have $B_J=\Z$ and
the homomorphism $f_J$ factors through the
canonical homomorphism
$f:B\ra \oplus_{t\in T}\Z\rtimes W$
where $T$ is the set of
all reflections in $W$ (=the set of elements that are conjugate to
some element of $S$).
By definition $f(s)=(e_{\ov{s}},\ov{s})$,
where $e_t$ is the standard basis of $\oplus_{t\in T}\Z$.
        
\begin{prop} Assume that $W$ is finite. Then
$f$ is surjective with the kernel
$[P,P]$, the commutator subgroup of the pure braid group
$P\sub B$.
\end{prop}
   
\Pf . In the above geometric picture $f|_P$ corresponds to taking
the link indices of a loop in $X$ with complex
hyperplanes $H_{\C}$ for all $H\in\HH$. Thus, the map $f|_P$
can be identified with the natural projection
$\pi_1(X)\ra H_1(X)=\pi_1(X)/[\pi_1(X),\pi_1(X)].$
\ed
        
\subsection{Induction for actions of positive braid monoids}
There  is  a  version  of Theorem \ref{braidmain} concerning
the actions of braid monoids (or groups) on categories.
Namely, assume that we have a $B$-action on an additive category
$\CC=\oplus_{x\in W/W_J}\CC_x$ such that
the functor $T(b)$ corresponding  to  $b\in  B$  sends  $\CC_x$  to
$\CC_{\ov{b}x}$.  Assume  also  that  for $s\in S$ and $x\in
W/W_J$ such that $sx\neq x$, the functor
$T(s):\CC_x\ra\CC_{sx}$ is an equivalence. Then we can
reconstruct $\CC$ and the action of $B$ on it from the
action  of  $B_J$  on  $\CC_{x_0}$  exactly  as  in  Theorem
\ref{braidmain}. More precisely, if
$w\in  W$ is a canonical representative of $x\in W/W_J$, then
the functor $T(w):\CC_{x_0}\ra \CC_x$ is an equivalence,
and the corresponding equivalence
$$\CC\simeq \oplus_{x\in W/W_J}\CC_{x_0}$$
is compatible with the $B$-actions, where the action of $B$
on the right-hand side is the composition of the natural
$B_J^{W/W_J}\rtimes W$-action with $f_J$.
      
This   is  reflected  in  the  following  result  concerning
$W$-gluing.
Let $(\CC_w,\Phi)$ be a $W$-gluing data, and $J\sub S$
a subset of simple reflections. Then it induces a gluing
data $\Phi_J$ for the categories $(\CC_w), w\in W_J$.
        
\begin{thm}
Assume that for every
$s\in S$ and every $w\in W$ such that $w^{-1}sw\not\in W_J$,
the morphism $\nu_s:\Phi_s^2|_{\CC_w}\ra\Id_{\CC_w}$
is an isomorphism. Then the glued category $\CC(\Phi)$
is equivalent to $\CC(\Phi_J)$.
\end{thm}
          
\Pf . Note that our condition on the gluing functors
implies that for every $w\in W$, which is the shortest
element of $wW_J$ and every $w'\in W_J$,
the functors $\Phi_w|_{\CC_{w'}}$ and
$\Phi_{w^{-1}}|_{\CC_{ww'}}$ are quasi-inverse to each other.
     
We claim that the functors
$j^*=j^*_{W_J}:\CC(\Phi)\ra\CC(\Phi_J)$
and $j_!=j_{W_J,!}:\CC(\Phi_J)\ra\CC(\Phi)$ are
quasi-inverse to each other. Indeed, we always have
$j^*j_!=\Id$. Now let $A=(A_w; \a_{w,w'})$ be an object of
$\CC(\Phi)$. The canonical adjunction morphism
$j_!j^*A\ra A$ has as components the morphisms
$$\a_{n(w),p(w)}:\Phi_{n(w)}A_{p(w)}\ra A_w.$$
Since the functors $\Phi_{n(w)}$ and $\Phi_{n(w)^{-1}}$
between $\CC_{p(w)}$ and $\CC_w$ are quasi-inverse to each
other, the associativity condition on $\a$ implies
that $\a_{n(w),p(w)}$ is an isomorphism.
\ed
   
\section{Good representations of braid groups}\label{goodsec}
   
\subsection{Some ideals in the group ring of the pure
braid group}\label{ideals}
         
Let $(W,S)$ be a Coxeter group, and $B$ and $P$  the
corresponding (generalized) braid group and pure braid groups, respectively.
Recall that $P$ is the kernel of the
natural homomorphism $B\ra W:b\mapsto\ov{b}$.
In other words, $P$ is the normal
closure  of the elements $\{ s^2, s\in S\}$ in  $B$.
Below we view $S$ as a subset in $B$ and denote by $\ov{s}$
($s\in S$) the corresponding elements in $W$.
        
\begin{thm} There exists a unique collection of right ideals
$(I_w, w\in W)$ in $\Z[P]$ such that $I_1=0$,
$$I_{\ov{s}}=(s^2-1)\Z[P]$$ for every $s\in S$, and
$$I_{ww'}=I_w+\tau(w)I_{w'}\tau(w)^{-1}$$
for every pair $w,w'\in W$ such that $l(ww')=l(w)+l(w')$.
\end{thm}
         
\Pf . We use the induction in $l(w)$, so we assume that
all $I(y)$ with $l(y)<l(w)$ are already constructed and
satisfy the property
$$I_{yy'}=I_y+\tau(y)I_{y'}\tau(y)^{-1}$$
for $y,y'\in W$ such that $l(yy')=l(y)+l(y')<l(w)$.
Choose a decomposition $w=\ov{s}w_1$ with $l(w)=l(w_1)+1$. Then
we must have $I_w=I_{\ov{s}}+sI_{w_1}s^{-1}$. It remains to show that the
right-hand side
does not depend on a choice of decomposition.
Let $w=\ov{s'}w'_1$ be another decomposition with $l(w)=l(w'_1)+1$.
Assume first that $s$ and $s'$ commute. Then $w_1=\ov{s'}y$ and
$w'_1=\ov{s}y$ where $l(y)=l(w)-2$, so that by the induction hypothesis
we have
$$I_{w_1}=I_{\ov{s'}}+s'I_{y}(s')^{-1},$$
$$I_{w'_1}=I_{\ov{s}}+sI_{y}s^{-1}.$$
It follows that
$$I_{\ov{s}}+sI_{w_1}s^{-1}=
I_{\ov{s'}}+s'I_{w'_1}(s')^{-1}=I_{s}+ I_{s'}+(ss')I_y(ss')^{-1}.$$
Now if $s$ and $s'$ do not commute there is a defining relation of the form
$$ss's\ldots=s'ss'\ldots$$
where both sides have the same length $m$.
Let us write this relation in the form
$$sr=s'r'$$
where $r=s's\ldots$ and $r'=ss'\ldots$ are the corresponding elements of length
$m-1$.
In this case we have $w_1=\ov{r}y$, $w'_1=\ov{r'}y$ where $l(y)=l(w)-m$,
so by the induction hypothesis we have
$$I_{w_1}=I_{\ov{r}}+rI_{y}r^{-1}.$$
Hence
$$I_{\ov{s}}+sI_{w_1}s^{-1}=
I_{\ov{s}}+sI_{\ov{r}}s^{-1}+(sr)I_y(sr)^{-1}.$$
Comparing this with the similar expression for
$I_{\ov{s'}}+s'I_{w'_1}(s')^{-1}$ we see that it is sufficient to
prove the equality
$$I_{\ov{s}}+sI_{\ov{r}}s^{-1}=
I_{\ov{s'}}+s'I_{\ov{r'}}(s')^{-1}.$$
In other words, we have to check that
\begin{equation}
I_{\ov{s}}+sI_{\ov{s'}}s^{-1}+ss'I_{\ov{s}}(ss')^{-1}+\ldots=
I_{\ov{s'}}+s'I_{\ov{s}}(s')^{-1}+s'sI_{\ov{s'}}(s's)^{-1}+\ldots
\end{equation}
where both sides contain $m$ terms. More precisely, we
claim that the terms of the right-hand side coincide
with those of the left-hand side in inverse order.
This follows immediately from the identity
$rI_{s''}r^{-1}=I_s$ where $s''=s$ if $m$ is even, $s''=s'$
if $m$ is odd.
\ed
         
Note that $I_w$ is generated by $l(w)$ elements and
$I_w\sub I_y$ if $y=ww'$ and $l(ww')=l(w)+l(w')$.
In the case when
$W$ is finite, this implies that $I_{w_0}$ contains all the ideals $I_w$,
where $w_0$ is the longest element in $W$.
         
\begin{cor} Assume that $W$ is finite. Then
$I_{w_0}=I$ where $I$ is the augmentation ideal in $\Z[P]$.
\end{cor}
         
\Pf. Let $s\in S$. Then we have
$$I_{w_0}=I_s+sI_{sw_0}s^{-1}.$$
Hence,
$$sI_{w_0}s^{-1}=I_s+s^2I_{sw_0}=I_s+I_{sw_0};$$
here the last equality follows from the definition of $I_s$.
In particular, $sI_{w_0}s^{-1}\sub I_{w_0}$ for every $s\in S$.
It follows that
$$s^{-1}I_{w_0}s=s^{-2}(sI_{w_0}s^{-1})\sub s^{-2}I_{w_0}\sub
I_{w_0}+I_s=I_{w_0}.$$
Thus, $I_{w_0}$ is stable under the action of $B$ by conjugation.
Since it also contains $(s^2-1)$ for every $s\in S$, it should be
equal to $I$.
\ed
                                                
\subsection{Main definition}
       
Let $V$ be a representation of the braid group $B$.
For every $w\in W$ let us denote
$$V_w=I_wV\sub V$$
where $I_w\sub\Z[P]$ is the ideal defined above.
Then $V_s=(s^2-1)V$ for $s\in S$, and
$$V_{w_1w_2}=V_{w_1}+\tau(w_1)V_{w_2}$$
if $l(w_1w_2)=l(w_1)+l(w_2)$.
Define
$K_W(V)$ to be the following subspace of $V^W$:
$$K_W(V)=\{(x_w), w\in W:  x_w\in V, x_{\ov{s}w}-sx_w\in V_s,
\all s\in S, w\in W\}.$$
    
\begin{prop} For every $(x_w)\in K_W(V)$ and every
$w_1,w_2\in W$, one has
$x_{w_1w_2}-\tau(w_1)x_{w_2}\in V_{w_2}$.
For every $y\in W$
there is a natural embedding
$i_y:V\ra K_W(V)$ given by
$$i_y(v)=(\tau(wy^{-1})v, w\in W).$$
\end{prop}
                    
The proof is straightforward.
    
\begin{ex} Let $R=\Z[u,u^{-1},(u^2-1)^{-1}]$ where $u$ is an indeterminate.
It follows from Theorem \ref{K_0A} and Corollary \ref{s2-1}
that the $R$-module $K_0(\AA)\otimes_{\Z[u,u^{-1}]}R$ is naturally
isomorphic to $K_W(K_0(G/U)\otimes_{\Z[u,u^{-1}]} R)$ where
$\AA$ is the Kazhdan---Laumon category. The embeddings
$i_y$ correspond to functors $Lj_{y,!}$ (left adjoint to restrictions).
\end{ex} 
     
\begin{defi} A $B$-representation $V$ is called {\it good}
if $K_W(V)$ is generated by the subspaces $i_y(V)$, $y\in W$.
\end{defi}

In the situation of the above example it would be very desirable
to prove that an appropriate localization of $K_0(G/U)$ is
a good representation of $B$ since this would imply that the objects
of finite projective dimension generate the corresponding localization
of $K_0(\AA)$. Then the following result could be applied to construct a
bilinear pairing on appropriate finite-dimensional quotient of $K_0(\AA)$.

\begin{prop}\label{form}
Let $V$ be a finite-dimensional representation of $B$, and
$V^*$  the dual representation.
Then both $V$ and $V^*$ are good if and only if there exists a
non-degenerate pairing $\chi:K_W(V)\otimes K_W(V^*)\ra\C$ 
such that $\chi(i_y(v_y),v')=\lan v_y, p_y v'\ran$
and $\chi(v,i_y(v'_y))=\lan p_y v, v'_y\ran$ for every
$y\in W$, $v_y\in V$, $v'_y\in V^*$, $v\in K_W(V)$,
$v'\in K_W(V^*)$.
\end{prop}

\Pf . If such a pairing $\chi$ exists, then
the orthogonal to the subspace in $V$ generated by
elements of the form $i_y(v_y)$, $y\in W$, $v_y\in V$,
is zero, and hence, $V$ is good. Similarly, $V^*$ is good. 
Now assume that both $V$ and $V^*$ are good.
Then we have the surjective map
$$
\pi:\oplus_{y\in W} V\ra K_W(V):
(v_y,y\in W)\mapsto \sum_y i_y(v_y)$$
and the similar map $\pi':\oplus_{y\in W} V^*\ra K_W(V^*)$.
It is easy to see that the kernel of $\pi$
consists of collections $(v_y,y\in W)$ such that
$\sum_{y\in W}\tau(wy^{-1})v_y$ for every $w\in W$. 
The kernel of $\pi'$ has a similar description.
Now we define a pairing
$$\wt{\chi}:(\oplus_{y\in W} V)\otimes(\oplus_{y\in W}V^*)\ra\C$$ 
by the formula
$$\wt{\chi}((v_y,y\in W),(v'_y,y\in W))=
\sum_{y,w\in W}\lan v_w,\tau(wy^{-1})v'_y\ran=
\sum_{y,w\in W}\lan \tau(yw^{-1})v_w,v'_y\ran.$$
Since the left and right kernels of $\wt{\chi}$ coincide with
$\ker(\pi)$ and $\ker(\pi')$, respectively, it descends
to the non-degenerate pairing $\chi$ between $K_W(V)$ and $K_W(V^*)$.
\ed 

Let  $Br_n$  be  the Artin braid  group  such  that the
corresponding Coxeter group is the symmetric group $S_n$.
    
\begin{prop}\label{B2} Any representation of $Br_2=\Z$ is good.
\end{prop}
    
\Pf . A representation of $Br_2$ is a space $V$ with  an
operator $\phi:V\ra V$. The corresponding space $K_{S_2}(V)$
consists of pairs $(v_1,v_2)\in V^2$ such that
$v_2-\phi v_1\in (\phi^2-1)V$. The maps $i_1$ and $i_s$
(where $s$ is the generator of $Br_2$) are given by
$i_1(v)=(v,\phi v)$, $i_s(v)=(\phi(v), v)$.
Now for any $(v_1,v_2)\in K_{S_2}(V)$ we have
$$(v_1,v_2)=i_1(v_1+\phi(v))+i_s(v)$$
where $v\in V$ is such that
$v_2-\phi v_1=(\phi^2-1)v$.
\ed
                            
\begin{thm}\label{B3}
Let  $V$  be  a  representation of $Br_3$ over a
field $k$, such that for generators $s_1$ and $s_2$, the
following identity is satisfied in $V$:
$$(s^2-1)(s-\la)=0$$
where $\la\in k$ is an element such that
$\la^6\neq 1$. Then $V$ is good.
\end{thm}
    
\subsection{Criterion}  We  need  an auxiliary result which allows one to
check that a subspace of $V^2$ is of the form $K_{S_2}(V)$
for some action of $Br_2=\Z$ on $V$. It is more natural
to generalize this construction as follows. Let
$V=V_0\oplus V_1$ be a super-vector space, and $\phi:V\ra V$
an odd operator. Then $\phi^2$ is even, and
we have the super-subspace $V_{\phi}=(\phi^2-1)V\sub V$.
Now we define the vector space
$K(\phi)=\{v\in V|\ \phi(v)-v\in V_{\phi}\}$.
Note that this is a non-homogeneous subspace of $V$,
equipped with two sections $i_0$ and $i_1$ of
the projections to $V_0$ and $V_1$, namely,
$i_{\a}(v_{\a})=v_{\a}+\phi(v_{\a})$ for $\a=0,1$. We
want to characterize non-homogeneous subspaces of
$V$ arising in this way.
    
\begin{prop}\label{gluecheck}
Let $K\sub V$ be a non-homogeneous subspace, and
$i_0:V_0\ra K$ and $i_1:V_1\ra K$  sections of the
projections of $K$ to $V_0$ and $V_1$. Let $K^h$ be
the maximal homogeneous subspace of $K$. Let
$\phi:V\ra V$ be the odd operator with components
$p_1i_0$ and $p_0i_1$ where $p_{\a}:K\ra V_{\a}$
are the natural projections. Then we have inclusions
$V_{\phi}\sub K^h$ and $K(\phi)\sub K$. The equality
$K=K(\phi)$ holds if and only if $V_{\phi}=K^h$.
\end{prop}
    
\subsection{Proof of Theorem \ref{B3}}. We will prove first that
$K_{S_3}(V)=K_{S_2}(V')$ where $V'$ is some representation of
$Br_2$. Namely, let $V'$ be the space of triples
$(x,y,z)\in V^3$ such that $y-s_1x\in V_{s_1}$ and
$y-s_2z\in V_{s_2}$. We have the natural embedding
$$\kappa:K_{S_3}(V)\hra V'\oplus V':
(x_w,w\in S_3)\mapsto ((x_1,x_{s_1},x_{s_2s_1}),
(x_{s_2},x_{s_1s_2},x_{s_2s_1s_2})).$$
To apply Proposition \ref{gluecheck} to the
image of $\kappa$ we have to define maps
$i_0,i_1:V'\ra \kappa K_{S_3}(V)$ such that
$p_{\a}i_{\a}=\id$ for $\a=0,1$. Let us set
\begin{align*}
&i_0(x,y,z)=((x,y,z),(s_2x,s_1s_2x-s_2s_1s_2y+s_2s_1z,s_1z)),\\
&i_1(x,y,z)=((s_2x,s_1s_2x-s_2s_1s_2y+s_2s_1z,s_1z),(x,y,z)).
\end{align*}
Then the operator $\phi=p_1i_0=p_0i_1:V'\ra V'$ has the following
form:
$$\phi(x,y,z)=(s_2x,s_1s_2x-s_1s_2s_1y+s_2s_1z,s_1z).$$
Now according to Proposition \ref{gluecheck} in order to
check that $\kappa K_{S_3}(V)$ is obtained by gluing from
$(V'\oplus V',\phi)$ we have to check that the
image of the operator $\phi^2-\id$ is precisely the
subspace $U=\{(x,y,z)\in V'|\ x\in V_{s_2}, z\in V_{s_1}\}$.
Let $u\in U$ be an arbitrary element. Then
$$u=(\phi^2-\id)(\la^2-1)^{-1}u+u'$$
with $u'=(0,y,0)\in U$. Then $y\in V_{s_1}\cap V_{s_2}$.
Hence, $\tau(w_0)$ acts on $y$ as multiplication by $\la^3$.
It follows that $u'=(\phi^2-\id)(\la^6-1)^{-1}u'$.
Hence, $u'$ and $u$ are in the image of $\phi^2-\id$
as required.
                
It     follows     from     Proposition     \ref{B2}    that
$K_{S_3}(V)=K_{S_2}(V')$ is generated by the images of
$i_0$ and $i_1$. Thus, we are reduced to show that $V'$ is
generated by elements of the form $(x,s_1x,s_2s_1x)$,
$(s_1y,y,s_2y)$ and $(s_1s_2z,s_2z,z)$ which is
straightforward.
\ed

\subsection{Representations of quadratic Hecke algebras are good} 
Let $H_q$ be the Hecke algebra of $(W,S)$ over
$\C$ with complex parameter $q\in\C$. Recall that
$H_q$ is the quotient of the group algebra $\C[B]$ where
$B$ is the corresponding braid group by the quadratic relations
$(s-q)(s+1)=0$ for every $s\in S$.
   
\begin{thm}\label{quHecke} Assume that $q$ is not a root of unity.
Then every representation of $H_q$ is good.
\end{thm}
   
\begin{cor}\label{qubraid} Let $V$ be a representation of $B$ such that
$(s-\lambda)(s-\mu)=0$   for    every    $s\in    S$,    where
$\lambda\in\C^*$,  $\mu\in\C$, and $\frac{\mu}{\la}$ is not a root of
unity. Then $V$ is good.
\end{cor}
   
Recall that we denote $\pi=\tau(w_0)^2\in B$
where $w_0\in W$ is the longest element in $W$.
   
\begin{lem}\label{mainHecke} Assume that $q$ is not a root of unity. Let
$V$  be an irreducible representation of $H_q$.
Then either $\pi-1$ acts by a non-zero constant on $V$
or $V$ is the one-dimensional representation such that
$s=-1$ on $V$ for every $s\in S$.
\end{lem}
   
\Pf . Let $E$  be an irreducible representation of $W$, $E(u)$ be the
corresponding irreducible representation of the Hecke algebra
$H$ of $(W,S)$  over  $\Q[u^{\frac{1}{2}},u^{-\frac{1}{2}}]$
defined by Lusztig. According to \cite{Lus} (5.12.2) one has
$$\pi=u^{2l(w_0)-a_E-A_E}$$
on $E(u)$ where the integers $a_E\ge0$ and $A_E\ge0$ are the lowest
and the highest degrees of $u$ appearing with non-zero
coefficient in the formal dimension $D_E(u)$ of $E$. Recall
that $D_E(u)$ is defined from the equation
$$D_E(u)\cdot\sum_{w\in W} u^{-l(w)}\Tr(\tau(w),E(u))^2=
\dim(E)\cdot\sum_{w\in W}u^{l(w)}.$$
It is known that $a_{E\ot sign}=l(w_0)-A_E$
where $sign$ is the sign representation of $W$.
In particular, $a_E\le A_E\le l(w_0)$.
Thus, we have $\pi=u^l$ on $E(u)$, where $l>0$ unless
$a_E=A_E=l(w_0)$. In the latter case
$D_E(u)=f_E^{-1}u^{l(w_0)}$ where $f_E>0$ is an integer.
Thus, $D_E(1)=\dim(E)=f_E^{-1}=1$ and $E=sign$.
\ed
   
\noindent
{\it   Proof  of  Theorem  \ref{quHecke}}.  Since  $H_q$  is
finite-dimensional it is sufficient to prove the assertion
for a finite-dimensional representation $V$ of $H_q$.
Clearly, we can assume that $V$ is irreducible. Assume first
that $V$  is one-dimensional and $s=-1$ on $V$ for every $s\in S$.
Then  $V_s=0$  for every $s$ so $K_W(V)\simeq V$ and all the
maps $i_y$  are  isomorphisms,  hence  $V$  is  good.  Thus,
according  to  Lemma  \ref{mainHecke}  we  can  assume  that
$\pi-1\neq 0$ on $V$. Now we are going to use
the  $K_0$-analogue  of  the  complex  defined  in   section
\ref{complex}. First of all for every coset $W_Jx\in W$
we have the map
$i_{W_Jx}:K_{W_J}(V)\ra K_W(V)$ which is a section of
the projection $p_{W_Jx}$ onto components $W_Jx\sub W$.
Namely, using notation of section \ref{complex}
the $w$-component of $i_{W_Jx}(v_y,y\in W_Jx)$
is $n_{W_Jx}(w)(v_{p_{W_Jx}(w)})$.
To prove that $V$ is good it is sufficient to show that
$K_W(V)$ is generated by images of all maps $i_{W_Jx}$ where
$J\sub  S$  is  a  proper  subset.  Now the proof of Theorem
\ref{homology} shows that for every $v\in K_W(V)$, we have the identity
$$\sum_{J\sub S,|J|<n,x\in W_J\backslash W}
(-1)^{|J|}i_{W_Jx}p_{W_Jx}(v)=\iota(v)+(-1)^{n-1}v,$$
where $W_{\emptyset}=\{ 1\}$, $\iota(v_w,w\in W)=(w_0v_{w_0w},w\in W)$.
One has $\iota^2(v_w,w\in W)=(\pi v_w,w\in W)$.
It follows that for every $v\in K_W(V)$ the element
$(\pi-1)v$ is a linear combination of elements of the form
$i_{W_Jx}(v')$ where $J\sub S$ is a proper subset, and hence
$v$ itself is such a linear combination.
\ed
   
\subsection{Good representations and cubic Hecke algebras}
Let   $H^c_q$  be  the  cubic  Hecke  algebra  with  complex
parameter $q$, i.e., the quotient of the group algebra $\C[B]$
of the braid group $B$ of $(W,S)$ by the relations
$(s-q)(s^2-1)=0$, $s\in S$,
where $q\in\C$ is a constant. An optimistic conjecture would be that
if $q$ is not a root of unity, then every representation of $B$
that factors through $H^c_q$ is good.
The particular cases are Proposition \ref{B2} and Theorems
\ref{B3}, \ref{quHecke}. Below we check some other particular
cases of this conjecture. However, it seems that in general one needs
to add some higher polynomial identities on the generators in order
for such a statement to be true (see \cite{Wenzl} for an example of
such identities).
   
\begin{thm}\label{thmAn} Let $(W,S)$ be of type $A_n$ for $n\le 3$. 
Then every representation of $H^c_q$ is good provided that
$q$ is not a root of unity.
\end{thm}
               
\begin{thm}\label{B_2}   Let  $(W,S)$  be  of  type  $B_2$.
Then every representation of $H^c_q$ is good
provided that $q^8\neq1$.
\end{thm}
               
The structure of the proof of these two theorems is the same
as of Theorem \ref{B3}.
We choose a simple reflection $s_i\in S$ ($s_2\in S_3$ in
Theorem \ref{B3}) and consider the
partition of $W$ into two pieces: $P_i$ and $P_is_i$.
This way we get an inclusion
$K_W(V)\hra V'_0\oplus V'_1$ where
$V'_0\sub  \oplus_{P_i}V$  is   the   space   of   collections
$(v_w,w\in P_i)$ such that $v_{sw}-sv_{w}\in V_s$ whenever
$w, sw\in P_i$; $V'_1$ is the similar space for $P_is_i\sub W$
instead  of  $P_i$. Let $p_0$ and $p_1$ be the corresponding
projections of $K_W(V)$ to $V'_0$ and $V'_1$.
We are going to construct sections $s_0$ and $s_1$ of these
projections.
For every coset $W_Jx\sub P_i$ (resp. $W_Jx\sub P_is_i$)
let us consider the map $i_{W_Jx,P_i}=p_0i_{W_Jx}:K_{W_J}\ra V'_0$
(resp. $i_{W_Jx,P_is_i}=p_1i_{W_Jx}:K_{W_J}\ra V'_1$).   
The   proof   of    Theorem
\ref{complexPi} shows that for every $v'\in V'_0$ we have the
identity
$$\sum_{J,x}
(-1)^{|J|-1}i_{W_{S-J}x,P_i}p_{W_{S-J}x}(v')=v',$$
where the sum is taken over non-empty subsets $J\sub S$ 
and $x\in W_{S-J}\backslash W$ such that
$W^{(j)}x\sub P_i$ for every $j\in J$.
Now we define $s_0$ to be the similar sum
$$s_0(v')=\sum_{J,x}(-1)^{|J|-1}i_{W_{S-J}x}p_{W_{S-J}x}(v').$$

\begin{lem}
For every $y\in P_i$ one has
$i_y=s_0i_{y,P_i}$. 
\end{lem}

\Pf . Note that for a coset $W_Jx\sub P_i$ we have
$p_{W_Jx}i_{y,P_i}=i_{x',W_Jx}$, where $x'$ is the element
of $W_Jx$ closest to $y$. Now the proof of Theorem \ref{complexPi}
shows that we have a formal equality
$$\sum_{J,x}(-1)^{|J|-1}x'(W_{S-J}x)=y,$$
where the sum is taken over non-empty subsets $J\sub S$ 
and $x\in W_{S-J}\backslash W$ such that
$W^{(j)}x\sub P_i$ for every $j\in J$,
$x'(W_{Kx})$ denotes the element in $W_Kx$ closest to $y$.
The assertion follows immediately, since
for every $x'\in W_Kx$ we have $i_{W_Kx}i_{x',W_Kx}=i_{x'}$.
\ed

Similarly, one constructs the section
$s_1:V'_1\ra K_W(V)$.
Thus, it is sufficient to show that the subspace
$K_W(V)\sub V'_0\oplus V'_1$ coincides with the subspace $K(\phi)$
where $\phi|_{V'_0}=p_1s_0:V'_0\ra V'_1$,
$\phi|_{V'_1}=p_0s_1:V'_1\ra V'_0$. By Proposition \ref{gluecheck}
this amounts to showing that the image of $\phi^2-\id$
surjects onto $\ov{V}_0\oplus \ov{V}_1$, where e.g.,
$\ov{V}_0=K_W(V)\cap V'_0$ consists of elements $(v_y, y\in P_i)$
such that 
$v_y\in V_{ys_iy^{-1}}$ whenever $ys_iy^{-1}\in S$.
One way to show it would be to consider some natural
filtration on $\ov{V}_0$ preserved by $\phi^2$ and considering the
induced map on the associated graded factors.
The filtration is defined as follows. The subspace
$F_k(V'_0)\sub \ov{V}_0$ for $k\ge 0$ consists
of $(v_y,y\in P_i)\in \ov{V}_0$ such that
$v_y=0$ whenever $l(ys_iy^{-1})\le 2k-1$. There is a similar
filtration on $\ov{V}_1$ (replace $P_i$ by $P_is_i$ everywhere).
We claim that $\phi$ is compatible with these filtrations.
As we will see below this is a consequence of the definition of $s_0$
and of the following result.

\begin{prop}\label{geod_yw} Let $y,w\in W$ be a pair of elements. Assume that
$l(ys_iy^{-1})>l(ws_iw^{-1})$. Then no geodesics from $y$ to $w$
pass through $ys_i$.
\end{prop}

\Pf . Let us denote $r=ws_iw^{-1}$, $r'=ys_iy^{-1}$, $w_1=yw^{-1}$.
The existence of a geodesic from $y$ to $w$ passing through $ys_i$
is equivalent to the existence of a geodesics from $w_1$ to $1$
passing through $ys_iw^{-1}=w_1r$. In other words this is equivalent
to the equality $l(w_1)=l(w_1r)+l(w_1rw_1^{-1})$, i.e.,
$l(w_1)=l(w_1r)+l(r')$ (since $w_1rw_1^{-1}=r'$).
But one has $l(w_1)\le l(w_1r)+l(r)<l(w_1r)+l(r')$, which contradicts
the above equality.
\ed

\begin{cor} The map $\phi:V'_0\ra V'_1$ sends $F_k(\ov{V}_0)$ to
$F_k(\ov{V}_1)$.
\end{cor}

\Pf . As we have seen in the proof of Theorem \ref{complexPi},
in the sum defining $p_ws_0$
all the terms cancel except those that correspond to $x\in P_i$,
such that the set $\{s\in S\ |\ l(wx^{-1}s)<l(wx^{-1})\}$ coincides
with the set of $s_j$ such that $W^{(j)}x\sub P_i$. By Lemma \ref{appear}
the latter set coincides with the set $S(xs_ix^{-1})$
of simple reflections appearing
in any reduced decomposition of $xs_ix^{-1}$. In particular,
this implies that for any $w_1\in W_{S(xs_ix^{-1})}$ there exists
a geodesic from $x$ to $w$ passing through $w_1x$. Taking $w_1=xs_ix^{-1}$
we see that for such $x$ there exists a geodesic from $x$ to $w$
passing though $xs_i$. According to Proposition \ref{geod_yw}
this implies that $l(xs_ix^{-1})\le l(ws_iw^{-1})$, hence the
assertion. 
\ed

\begin{lem} Let $w\in W$ be an element, $r\in W$  a reflection, and
$S(r)$ the set of simple reflections in a reduced decomposition
of $r$. Assume that $l(wrw^{-1})=l(r)$ and
the set $\{ s\in S\ |\ l(ws)<l(w)\}$ coincides with $S(r)$.
Assume also that $W$ is either of type $A_n$ or of rank 2. 
Then $w$ is the longest element of $W_{S(r)}$.
\end{lem}

\Pf . The rank 2 case is straightforward, so let us prove the statement
for $W=S_n$. Let $r=(ij)$ where $i<j$. Then the assumptions of the lemma are:
1) for any $k\in [1,n]$ one has $w(k)>w(k+1)$ if and only if
$k\in [i,j-1]$, 2) $j-i$=$w(i)-w(j)$. It follows that $w$ maps the
interval $[i,j]$ to the interval $[w(j),w(i)]$ reversing the order
of elements in this interval. One the other hand, $w$ preserves the
order of elements in $[1,i]$ and in $[j,n]$. Hence $w(i)=j$, $w(j)=i$
and $w$ stabilizes all elements outside $[i,j]$.
\ed

It follows from this lemma
that the map 
$$\gr\phi:\gr^F_k(\ov{V}_0)\ra \gr^F_k(\ov{V}_1)$$
sends $(v_y,y\in P_i,l(ys_iy^{-1})=2k+1)$ to
$(v_w,w\in P_is_i,l(ws_iw^{-1})=2k+1)$ where
$$v_w=(-1)^{|S(ws_iw^{-1})|-1}\tau(w_0(w))v_{w_0(w)w};$$ 
here $S(ws_iw^{-1})$ is
the set of simple reflections in any reduced decomposition of 
$ws_iw^{-1}$, and $w_0(w)$ is the longest element in $W_{S(ws_iw^{-1})}$.
It follows that $\gr\phi^2$ sends $(v_y)$ to $(\tau(w_0(y))^2v_y)$.
By induction it is easy to see that it would be enough to prove the
surjectivity of this map for the smallest filtration term
$F_l(\ov{V}_0)$ where $l$ is the maximal length of reflections in $W$.

Assume that $W$ is of classical type. 
Let $r_0$ be the reflection corresponding to the maximal positive root.
Then $r_0$ has maximal length in its conjugacy class.
We can choose $s_i\in S$ which is conjugate to $r_0$.
Then the space $F_l(\ov{V}_0)$ consists of collections 
$(v_y, y\in P_i, ys_iy^{-1}=r_0)$ such that $v_y\in V_{s}$ whenever
$sr_0s\neq r_0$, 
and $v_{sy}-sv_y\in V_{s}$ if $sr_0s=r_0$. 

If $W=S_n$ is the symmetric group, then $r_0$ is the transposition $(1,n)$.
The space $F_l(\ov{V}_0)$ consists of collections 
$(v_y, y\in P_i, ys_iy^{-1}=r_0)$ such that $v_y\in V_{s_1}\cap V_{s_{n-1}}$ 
and $v_{s_jy}-s_jv_y\in V_{s_j}$ for $1<j<n-1$. Note that
the set of $y\in P_i$ such that $ys_iy^{-1}=r_0$ constitutes one coset
for the subgroup $W_{[2,n-2]}\sub W$ generated by $s_j$ with $1<j<n-1$.

\vspace{2mm}

\noindent
{\it Proof of Theorem \ref{thmAn}}.
Let $W=S_4$. 
Then $s_jr_0\neq r_0s_j$ for $j\neq 2$
and $s_2r_0=r_0s_2$. We have an involution $(v_y)\mapsto (v_{s_2y})$
on $F_l(\ov{V}_0)$. Consider the corresponding decomposition
$F_l(\ov{V}_0)=F_l^+\oplus F_l^{-}$. Thus, $F_l^{+}$ consists 
of collections $(v_y, y\in P_i, ys_iy^{-1}=r_0)$ such that
$v_{s_jy}\in V_{s_j}$ for $j\neq 2$,
$v_{s_2y}=v_y$, $(s_2-1)v_y\in V_{s_2}$.
Thus, for any component $v_y$ we have $s_j(v_y)=qv_y$ for $j\neq 2$,
$(s_2-1)(s_2-q)v_y=0$. It is easy to deduce from these conditions
that the subrepresentation of $B$ generated by $v_y$ factors through
the quadratic Hecke algebra $H_q$. It remains to apply Theorem
\ref{quHecke}.
\ed

\subsection{Some linear algebra}  

\begin{lem}\label{findi}
Let $V$ be a representation of the cubic Hecke
$H^c_q$ corresponding to $(W,S)$ of type $B_2$, where $q\neq 0$. Assume that
$v\in V$ is such that $s_1v=q v$ and
$(s_2-q)(s_2-1)v=0$. Then the subspace $\C[B]v\sub V$
is finite-dimensional.
\end{lem}

\Pf . We claim that the subspace $V'$ spanned by elements
$$(x, s_2x, s_1s_2x, s_1^{-1}s_2x, s_2s_1s_2x, s_2s_1^{-1}s_2x,
s_1s_2s_1^{-1}s_2x)$$ 
is closed under $s_1$ and $s_2$.
To prove this first note that $s_1^{-1}s_2^{-1}s_1^{-1}$
commutes with $s_2$, and therefore, 
$$(s_2-q)(s_2-1)s_1^{-1}s_2^{-1}x=0.$$
Since $s_1^{-1}s_2x$ is a linear combination of $x$ and $s_1^{-1}s_2^{-1}x$
it follows that
$$(s_2-q)(s_2-1)s_1^{-1}s_2x=0.$$
Similarly since $s_1s_2s_1$ commutes with $s_2$ we get
$$(s_2-q)(s_2-1)s_1s_2x=0.$$ 
It follows that the subspace spanned by
$$(x, s_2x, s_1s_2x, s_1^{-1}s_2x, s_2s_1s_2x, s_2s_1^{-1}s_2x)$$
contains $s_2^ns_1^ks_2^l$ for any $n,k,l\in\Z$.
It remains to show that $s_2s_1s_2s_1^{-1}s_2x$ and
$s_1^2s_2s_1^{-1}s_2x$ belong to $V'$. Because of the cubic relation
we have $s_1^2s_2s_1^{-1}s_2x\equiv s_1^{-1}s_2s_1^{-1}s_2x\mod V'$.
Since $s_2s_1^{-1}s_2x$ is a linear combination of $s_1^{-1}s_2x$ and
$s_2^{-1}s_1^{-1}s_2x$, it follows that 
$s_1^{-1}s_2s_1^{-1}s_2x\in V'$. Finally,
$s_2s_1s_2s_1^{-1}s_2x=s_1^{-1}s_2s_1s_2^2x$ is a linear combination
of $s_1^{-1}s_2s_1s_2x=q^{-1}s_2s_1s_2x$ and
$s_1^{-1}s_2s_1x$, so we are done.
\ed

\begin{lem}\label{q8mu} 
Let $s_1$, $s_2$ be a pair of invertible operators
on a vector space $V$ such that $(s_1s_2)^2=(s_2s_1)^2$ and
$(s_1-q)(s_1^2-1)=0$, and let
$v\in V$ be a vector such that $s_1v= q v$,
$(s_2-q)(s_2-1)v=0$, and $(s_1s_2)^2v=\mu v$, where $q\neq 0$ and $\mu$
are constants.
Assume that $q^8\neq 1$. Then $\mu^2\neq 1$.   
\end{lem}

\Pf . 
First, we claim that the subspace of $V$ spanned
by $v$ and $s_2v$ is closed under $s_1$ and $s_2$.
Indeed,  the identity
$$\mu v=(s_2s_1)^2v=q s_2s_1s_2 v$$
implies that
$$s_1s_2v=q^{-1}\mu s_2^{-1}v=q^{-2}(q+1)\mu v-q^{-2}\mu s_2v.$$  
Thus, we can assume that $V$ is generated by $v$ and $s_2v$.
If $s_2v=\nu v$, then $\nu$ is equal to either $1$ or $q$
and $\mu=q^2\nu^2$, so the assertion follows. Otherwise,
$v$ and $s_2v$ constitute a basis of $V$. The matrix of $s_1$
with respect to this basis is
$\left( \matrix q & \frac{\mu(q+1)}{q^2}\\ 
0 & -\frac{\mu}{q^2} \endmatrix \right)$. Since $s_1$ is diagonalizable
with eigenvalues among $\{\pm 1,q\}$, it follows that
$-\frac{\mu}{q^2}=\pm 1$, i.e. $\mu=\pm q^2$.
\ed

\noindent
{\it Proof of Theorem \ref{B_2}}.
We have $r_0=s_2s_1s_2$, $s_i=s_1$. Thus,
the space $F_l(\ov{V}_0)$ consists of pairs
$(v_{s_2},v_{s_1s_2})\in (V_{s_2})^{\oplus 2}$ such that
$v_{s_1s_2}-s_1v_{s_2}\in V_{s_1}$. We have the natural involution
on $F_l(\ov{V}_0)$ interchanging the components, so that
$F_l^{\pm}$ consists of $(v,\pm v)$, such that $s_2v=qv$,
$(s_1-q)(s_1 \mp 1)v=0$. According to Lemma \ref{findi} 
we can assume that $V$ is finite-dimensional.
Also we can assume that $V$ is irreducible as a representation of $B$.
Then the central element $(s_1s_2)^2$ acts as a constant $\mu$ on $V$,
and we can apply Lemma \ref{q8mu} to finish the proof.  
\ed

\subsection{Good representations and parabolic induction} 
Let $J\sub S$ be a subset, $W_J\sub W$
 the   corresponding   parabolic   subgroup,  and $V_0$  a
representation of $B_J$, the braid group corresponding to
$(W_J,J)$.  In  section   \ref{induction} we associated with $V_0$
the representation of $B$ in
$$V=\oplus_{x\in W/W_J}V_x,$$ 
where $V_x=V_0$ such that
$b\in B$ sends $V_x$ to $V_{\ov{b}x}$.
   
\begin{prop}\label{goodind}
If $V_0$ is a good representation of $B_J$, then
$V$ is a good representation of $B$.
\end{prop}
   
\Pf . We have a direct sum decomposition
$$K_W(V)\simeq \oplus_{x\in W/W_J}K_{W,x}(V)$$
where
$$K_{W,x}(V)=\{(v_w,w\in W)\ |\ v_w\in V_{wx},
sv_w-v_{sw}\in (s^2-1)V_{wx} \}.$$
Furthermore, the map
$i_y:V\ra K_W(V)$ for $y\in W$
decomposes into the direct sum of maps
$i_{y,x}:V_{yx}\ra K_{W,x}(V)$ where $x\in W/W_J$.
Thus, it is sufficient to check that for every
$x\in W/W_J$ the images of $i_{y,x}$, $y\in W$ generate
$K_{W,x}(V)$. Let $\wt{x}\in W$ be a representative of $x$.
Then we have the isomorphism
$$K_{W,x}(V)\wt{\ra}K_{W,x_0}(V):
(v_w)\mapsto (v_{w\wt{x}^{-1}})$$
where $x_0\in W/W_J$ is the class containing $1$.
Under this isomorphism the map $i_{y,x}$ corresponds
to $i_{y\wt{x},x_0}$. Thus, we can assume that $x=x_0$.
Now we claim that the canonical projection
$$p_{W_J}:K_{W,x_0}(V)\ra K_{W_J}(V_0),$$
leaving only coordinates corresponding to $w\in W_J$, is an isomorphism.
Indeed, let $v=(v_w,w\in W)\in K_{W,x_0}(V)$ be an element.
For any $x=wx_0\in W/W_J$ and $s\in S$ such that $sx\neq x$,
we have $(s^2-1)V_x=0$, hence $v_{sw}=sv_w$. It follows that
for every $w\in W$ we have $v_w=n_{W_J}(w)v_{p_{W_J}(w)}$.
Conversely, the latter formula produces an element of
$K_{W,x_0}(V)$ from an arbitrary element of $K_{W_J}(V_0)$.
Note that  for  $y\in  W_J$, the map $i_{y,x_0}$ corresponds via
$p_J$ to the map $i_y:V_0\ra K_{W_J}(V_0)$. Since the
$B_J$-representation  $V_0$  is  good,   we   conclude   that
$K_{W,x_0}$ is generated by images of $i_{y,x_0}$, $y\in W_J$.
\ed

The  main  example  where the above proposition works is the
following. Consider the situation of gluing on the
basic affine space of a group $G$ where $G$ and $B$ are defined over
the finite field $\F_q$. Then we have the trace map
$\tr: K_0(\Perv(G/U)^{\Fr})\ra C( G/U(\F_q))$
where $\Perv(G/U)^{\Fr}$ is the category of perverse Weil sheaves on
$G/U$, and $C( G/U(\F_q))$ is the space of functions
$G/U(\F_q)\ra\ov{\Q}_l$.
The action of the braid group $B$ corresponding to $(W,S)$
on $\Perv(G/U)^{\Fr}$
is compatible with the trace map and induces the action of $B$
on $C( G/U(\F_q))$.
   
\begin{thm} Assume that the center of $G$ is connected.
Then the representation of the braid group $B$ on
$C(G/U(\F_q))$ is good.
\end{thm}
   
\Pf . We have an action of the finite torus $T(\F_q)$ on
$G/U(\F_q)$ by right translation. Let
$C_{\theta}\sub C(G/U(\F_q))$ be
the subspace on which $T(\F_q)$ acts through the character
$\theta:T(\F_q)\ra\ov{\Q}_l^*$. Then we have the direct
sum decomposition
$$C(G/U(\F_q))=\oplus_{\theta}C_{\theta}$$
and the action of $b\in B$ sends
$C_{\theta}$ to $C_{\ov{b}\theta}$.
We claim that if for some character $\theta$  and  a  simple
reflection $s$, one has $s\theta\neq\theta$; then
$s^2=1$       on       $C_{\theta}$.       Indeed,       let
$f:G/U(\F_q)\ra\ov{\Q}_l$ be a function such that
$f(xt)=\theta(x)f(x)$. Then
$$\sum_{\la\in\F_q^*}f(x\a_s(\la))=0,$$
hence, Proposition \ref{Fouriersq} implies that $s^2f=f$.
   
Now let  $O$  be  an  orbit  of  $W$ on the set of characters of
$T(\F_q)$. Clearly it is sufficient to prove that
the representation of $B$ on
$C_O=\oplus_{\theta\in O}C_{\theta}$
is good. Now we can choose a representative $\theta_0\in O$
in such a way that all reflections stabilizing $\theta_0$
belong to $S$. Then the stabilizer of $\theta_0$ in $W$
is the parabolic subgroup $W_J\sub W$ corresponding to
some subset $J\sub S$ (this follows from Theorem 5.13 of \cite{DL}).
As we have shown above the representation of $C_O$
belongs to  $\Mod_J-B$; hence  it  is  recovered  from  the
representation   of   $B_J$   on   $C_{\theta_0}$  via
the construction of section \ref{induction}. According
to Proposition \ref{goodind} it is sufficient to check
that the representation of $B_J$ on $C_{\theta}$ is good.
But on $C_{\theta}$ we have $(s+q^{-1})(s-1)=0$ for every
$s\in  J$  (this  follows  from  the analogue of Proposition
\ref{Fourproj} for the finite Fourier transorm);
hence we are done by Corollary \ref{qubraid}.
\ed

\subsection{Final remarks}
Recall that according to Theorem \ref{K_0A} and Corollary
\ref{s2-1} we have
$$K_0(\AA^{\Fr})\otimes_{\Z[u,u^{-1}]}\ov{\Q}_l\simeq
K_W(K_0(\Perv(G/U)^{\Fr})\otimes_{\Z[u,u^{-1}]}\ov{\Q}_l)$$ 
where $\AA^{\Fr}$ is the
analogue of Weil sheaves in the Kazhdan---Laumon category, 
and the homomorphism
$\Z[u,u^{-1}]\ra\ov{\Q}_l$ sends $u$ to $q$. Hence, we have
the natural trace map
\begin{equation}\label{tracemap}
\tr:K_0(\AA^{\Fr})\ra K_W(C( G/U(\F_q))),
\end{equation}
induced by the trace maps on every component.
If the center of $G$  is connected, then the previous theorem
implies the surjectivity of this map. 
Assuming that the cohomological dimension of $\AA$ is finite, Kazhdan
and Laumon defined a bilinear pairing on
$K_0(\AA^{\Fr})$ by taking traces of Frobenius on $\Ext$-spaces.
Without this assumption one can still define this pairing for
classes of objects of finite projective dimension.
It is easy to check that
the map (\ref{tracemap}) is compatible with this pairing
and the non-degenerate form on
$K_W(C( G/U(\F_q)))$ defined in Proposition \ref{form}
(where $C( G/U(\F_q))$ is equipped with the standard scalar product
in which delta-functions constitute an orthonormal basis).
If in addition we know that
$K_0(\Perv(G/U)^{\Fr})\otimes_{\Z[u,u^{-1}]}\ov{\Q}_l$ is
a good representation of $B$, then we can deduce that the map
(\ref{tracemap}) is the quotient of $K_0(\AA^{\Fr})$ by the kernel
of the pairing with $K_0(\AA^{\Fr})$.
Since the representation of the braid group on $K_0(G/U)$ factors
through the cubic Hecke algebra, we can check the latter implication 
in the cases where $(W,S)$ is of type $A_2$, $A_3$, or $B_2$.

On the other hand, using Corollary \ref{findimwt} we can always define
a bilinear pairing on $K_0(\AA_m)$ with values in the field of
Laurent series $\ov{\Q}_l((u))$ by looking at the action of Frobenius
on weight-$n$ components of the $\Ext$-spaces. In the cases when
the representation of the braid group on
$K_0(\Perv_m((G/U)_0))\otimes_{\Z[u,u^{-1}]}\Q(u)$
is good, the Laurent series obtained as
values of this pairing are actually rational.
We conjecture that in fact they are always rational. Furthermore,
a similar pairing can be defined in the case when the Frobenius
automorphism is twisted by the action of some $w\in W$. We conjecture
that it still takes values in rational Laurent series and that
the quotient of
$K_0(\AA_m^{w\Fr})\otimes_{\Z[u,u^{-1}]}\ov{\Q}_l((u))$
by the kernel of this pairing is finite-dimensional.

\section{Appendix. Counterexample to the conjecture of
Kazhdan and Laumon.}

In this appendix we will show that the Kazhdan---Laumon
category $\AA$ has infinite
cohomological dimension in the simplest non-trivial case
$G=\SL_3$.

For every $w\in W$ let us denote by 
$O_w$ the simple object in $\AA$, such that the 
corresponding gluing data is the constant sheaf
$\ov{\Q}_{l,X}[d]$ at the place $w\in W$,
and zero at all the others, where $d=\dim X$.
Recall that for every $w\in W$ the functor
$j_{w,!}:\Perv(X)\ra\AA$ has the left derived one
$Lj_{w,!}$
(see Proposition \ref{leftderived} and the
remark after it). So for every $w\in W$ we can introduce
the following object in the derived category of $\AA$:
$P_w=Lj_{w,!}(\ov{\Q}_{l,X}[d])$. 
Note that by Proposition \ref{leftderived}
the functor $Lj_{w,!}$ is left adjoint to the
corresponding restriction functor, hence, we have
$$V^*_{w,w'}:=\Ext^*(P_w,O_{w'})=\cases 0,\ w\neq w',\\ H^*(X)=H^*(G), w=w'.
\endcases$$
(here cohomology is taken with coefficients in $\ov{\Q}_l$).
We want to compare these spaces with the spaces
$$E^*_{w,w'}:=\Ext^*(O_{w}, O_{w'})$$
(in fact, $E^*_{w,w'}$ depends only on $w'w^{-1}$).
This is done with the help of the following lemma.

\begin{lem}\label{labelrom}
For every $i$ we have a canonical isomorphism in $\AA$
$$H^{-i}(P_w)= \oplus_{w', \ell(w')=i}  O_{w'w}(i).$$
\end{lem}

\begin{rem} Most of the results in this paper can be translated  into the
  parallel setting where algebraic $D$-modules on the complex
algebraic variety
$(G/U)_{\C}$  
  are used instead of 
$l$-adic sheaves over the corresponding variety in characteristic
$p$.  The main result of \cite{BBP} asserts that the corresponding
``glued''  category is equivalent to the category of modules over the
ring of global differential operators on $(G/U)_{\C}$.
The functor $Rj_*$ (the Verdier dual to the $D$-module
counterpart of the functor $Lj_!$ considered above; it is somewhat 
more natural to work with this Verdier dual functor in the $D$-module
setting) is then identified with the derived functor of global sections
from the category of $D$-modules to the category of modules over the
global differential operators. 

The $D$-module counterpart of  Lemma \ref{labelrom} is easily seen to be 
equivalent to the Borel-Weil-Bott Theorem (which computes cohomology
of an equivariant line bundle on the flag variety $(G/B)_{\C}$).
\end{rem}

\noindent 
{\it Proof of Lemma \ref{labelrom}}. 
Recall that by Proposition \ref{leftderived}
the composition of functors
$j^*_{w'w}\circ Lj_{w,!}$ coincides with the left derived functor of
$\sideset{^p}{^0}{H} F_{w',!}$. Furthermore,
Theorem \ref{adapted} implies that under the
identification of the derived category of $\Perv(X)$ with
$D^b_c(X,\ov{\Q}_l)$ the
latter derived functor coincides with $F_{w',!}$.
Therefore, we have
$$j^*_{w'w}P_w=j^*_{w'w}Lj_{w,!}(\ov{\Q}_{l,X}[d])\simeq
F_{w',!}(\ov{\Q}_{l,X}[d]).$$

According to Lemma \ref{mainl} for every simple reflection $s$ we have
$$F_{s,!}(\ov{\Q}_{l,X}[d])=(\ov{\Q}_{l,X}[d])[1](1).$$
Therefore, for every $w'\in W$ we have
$$F_{w',!}(\ov{\Q}_{l,X}[d])=(\ov{\Q}_{l,X}[d])[\ell(w')](\ell(w')).$$
Hence, for every $i$ we have
$$j^*_{w'w}H^{-i}(P_w)=\cases 0,\ \ell(w')\neq i,\\ \ov{\Q}_{l,X}[d](i),\
\ell(w')=i.\endcases$$
This immediately implies our statement.
\ed

Thus, we have a spectral sequence with the $E_2$-term
$$\oplus_{p\le 0,q\ge 0} \oplus_{\ell(w_1)=-p} E^q_{w_1w,w'}(p)$$
converging to $V^*_{w,w'}$.
Now let us assume that the spaces $E^*_{w,w'}$ are finite-dimensional
and lead this to contradiction.
Note that all our spaces carry a canonical (mixed) action of Frobenius
which is respected by this spectral sequence. We can encode some information
about Frobenius action on such a space
by considering Laurent polynomials in $u$, where
the coefficient with $u^n$ is the super-dimension of the weight-$n$
component. Let $e_{w,w'}$ (resp.
$v_{w,w'}$) be such a Laurent polynomial in $u$ corresponding to
$E^*_{w,w'}$ (resp. $V^*_{w,w'}$).
Then the above spectral sequence implies that
\begin{equation}\label{inconsist}
v_{w,w'}=\sum_{w_1\in W} e_{w_1w,w'}(-u)^{\ell(w_1)}
\end{equation}
Now recall that $v_{w,w'}=0$ for $w\neq w'$ while
$v_{w,w}$ is equal to the Poincare polynomial of $G$
$p_G(u)=\prod(1-u^{e_i})$, with $e_i$ running over the set
of exponents of $G$. 
Therefore, we can rewrite (\ref{inconsist}) in the matrix form
as follows.
Let us consider the matrix $E=(e_{w,w'})$ 
with rows and columns numbered by $W$
and with entries in $\Z[u,u^{-1}]$. Let us also define the matrix $M$
of similar type by setting
$$M_{w,w'}=(-u)^{\ell(w'w^{-1})}.$$
Then (\ref{inconsist}) is equivalent to the following equality of matrices:
\begin{equation}\label{matrixeq}
ME= p_G\cdot I
\end{equation}
where $I$ is the identity matrix.

Now we are going to show that $M$ does not divide $p_G\cdot I$ for 
$G=\SL(3)$
in $Mat_6[u,u^{-1}]$ (as $\det (M)$ vanishes at a primitive
$6$-th root of $1$, while $p_G$ does not).

Recall that for any finite group $G$ one can form a matrix
$M^G \in Mat_n (\Z[x_1, .., x_n])$,  $n=|G|$; here rows/columns of $M^G$
and variables $x$ of the polynomial ring are indexed by elements of $G$, and
we set $M^G_{g_1,g_2}=x_{g_1g_2}$.
It is well known that
$\det(M^G)$ is the product of factors indexed by irreducible representations
of $G$, and the degree of a factor, as well as the power in which
it enters the decomposition equals the dimension of the representation.
Moreover, the factor corresponding to a ($1$-dimensional) character
of $G$ is $\sum \chi(g) x_g$.

%Proof: let $A= k[x_g, g\in G]$; then $M^G$ acts on $A[G]$ and sends $A[G]$
%into the submodule generated by the element  $(x_{g_1}, ...,x_{g_n})$
%(under the action of both $A$ and $G$). Being  $G$-inv-t, this submodule
%decomposes into a direct sum over isotypic components.

Now our matrix $M$ is obtained from $M^G$ for $G=W$ by
$x_w\mapsto (-u)^{\ell (w)}$
and multiplication by a matrix of the permutation $w\mapsto w^{-1}$. 
So $\det(M)$ is divisible by $(\sum_{w\in W} u^{\ell(w)}) \cdot
(\sum_{w\in W} (-u)^{\ell(w)})$.
For $G=\SL(3)$ we get that $\det(M)$ is divisible by
$1+2u + 2 u^2 + u^3$ and $1 - 2u +2 u^2 -u^3$. 
The latter polynomial vanishes at a primitive $6$-th root of
$1$. However, $p_{\SL_3}=(1-u^2)(1-u^3)$,
hence the equality (\ref{matrixeq})
is impossible.

Note that this counterexample does not contradict to the results
of section \ref{goodsec} since in that section we considered the
localization of $K_0(\AA)$ on which $u$ does not act as a root of unity.

\vspace{2mm}

\noindent
Department of Mathematics, Boston University, Boston, MA 02215       

\noindent
{\it E-mail address}: apolish@@math.bu.edu

\vspace{3mm}

\noindent
Department of Mathematics, University of Chicago, Chicago, IL 60637

\noindent
{\it E-mail address}: roman@@math.uchicago.edu
        
\end{document}